\documentclass[12pt]{amsart}
\usepackage{amsmath, amssymb}
\usepackage[hypertexnames=false]{hyperref}
\usepackage{cleveref}
\usepackage{mathrsfs}
\usepackage{bm}
\usepackage{bbding}
\usepackage{textcomp}

\usepackage{fontawesome5}

\newcommand{\alphafootnote}[1]{%
	\begingroup
	\renewcommand{\thefootnote}{a}%
	\footnote{#1}%
	\endgroup
}

\usepackage[greek,english]{babel}

\usepackage[x11names]{xcolor}

\usepackage{graphicx}
\usepackage{amsthm,amscd}
\usepackage{calc}
\usepackage{mathtools}
\usepackage{CJK,fancyhdr,amscd}

\usepackage{anysize}
\usepackage{epsfig,enumerate}
\usepackage{latexsym}
\usepackage{indentfirst, latexsym}
\usepackage{graphics}
\usepackage{tikz-cd}
\usepackage{lineno}
\usepackage[all,poly,knot]{xy}

\usepackage{xpatch}
\usepackage[shortlabels]{enumitem}

\usepackage[T1]{fontenc}
\usepackage{textcomp}
\usepackage[utf8]{inputenc}

\newcommand{\zd}{\mathrm d}

\hypersetup{
	colorlinks=true,
	linkcolor=Green4,
	citecolor=DodgerBlue2,
	urlcolor=Magenta4}

\renewcommand{\eqref}[1]{\textcolor{Red2}{\hyperref[#1]{(\ref*{#1})}}}

\xpatchcmd{\step}{%
	\normalfont\scshape\centering}{%
	\normalfont\scshape}{\typeout{Success}}{\typeout{Failure}}%

\allowdisplaybreaks[2]

\providecommand{\U}[1]{\protect\rule{.1in}{.1in}}

\numberwithin{equation}{section}

\newtheorem{theorem} {Theorem} [section]
\newtheorem{bigthm}{Theorem}
 
\newtheorem{proposition}[theorem]{Proposition}
\newtheorem{corollary}  [theorem]     {Corollary}
\newtheorem{lemma}  [theorem]     {Lemma}
\theoremstyle{definition}
\newtheorem{example}  [theorem]     {Example}
\theoremstyle{plain}

\newtheorem*{questionno}{Question} 

\newtheorem{step}{Step}
\theoremstyle{definition}

\newtheorem{definition}  [theorem]     {Definition}
\newtheorem{notation}  [theorem]     {Notation}

\theoremstyle{plain}

\theoremstyle{definition}
\newtheorem{remark}  [theorem]     {Remark}

\newcounter{specialremarkcounter}[section]
\renewcommand{\thespecialremarkcounter}{\thesection.\arabic{specialremarkcounter}}

\newtheoremstyle{specialremarkstyle} 
{}                          
{}                          
{\normalfont}               
{}                          
{\itshape}                  
{.}                         
{ }                         
{\textit{Remark \thespecialremarkcounter}} 

\theoremstyle{specialremarkstyle}
\newtheorem{specialremark}[specialremarkcounter]{Remark}

\newcommand{\dz}{d \bar{z}}
\newcommand{\bfitDelta}{\bm{\mathit{\Delta}}}

\newcommand{\Upomega}{\boldsymbol{\Omega}}

\newcommand{\co}{\mathrm{co}}
\newcommand{\im}{\mathrm{Im}}

\newcommand{\cooo}{\mathpzc{co}}

\newcommand{\res}{\text{Res}}
\newcommand{\AAA}{\mathrm{\mathbf A}}

\newcommand{\dgbva}{DGBVA}
\newcommand{\dgla}{DGLA}

\newcommand{\Dt}{\Delta^{\prime\prime}}
\newcommand{\aaa}{a}
\newcommand{\qqq}{q}
\newcommand{\Do}{\Delta^{\prime}}
\newcommand{\ddd}{\mathrm{\mathbf d}}
\newcommand{\LLL}{\mathrm{\mathbf L}}
\newcommand{\Res}{\mathrm{Res}}

\newcommand{\nulll}{\mathrm{null}\,}
\newcommand{\btheorem}{\begin{theorem}}
	\newcommand{\etheorem}{\end{theorem}}
\newcommand{\bproposition}{\begin{proposition}}
	\newcommand{\eproposition}{\end{proposition}}
\newcommand{\bdefinition}{\begin{definition}}
	\newcommand{\edefinition}{\end{definition}}
\newcommand{\bcorollary}{\begin{corollary}}
	\newcommand{\ecorollary}{\end{corollary}}
\newcommand{\bproof}{\begin{proof}}
	\newcommand{\eproof}{\end{proof}}
\newcommand{\beq}{\begin{equation}}
	\newcommand{\eeq}{\end{equation}}
\newcommand{\ee}{\end{eqnarray*}}
\newcommand{\be}{\begin{eqnarray*}}

\newcommand{\elemma}{\end{lemma}}
\newcommand{\blemma}{\begin{lemma}}

\newcommand{\p}{\partial}
\newcommand{\mmm}{m}
\newcommand{\bp}{\overline\partial}

\newcommand{\DDD}{D^\prime_{h^*}}

\newcommand{\III}{\mathrm{I}}

\newcommand{\bd}{\begin{enumerate} }
	\newcommand{\ed}{\end{enumerate}}

\newcommand{\codl}{\mathpzc{co}\mathscr{D}_{L^*}}

\DeclareMathAlphabet{\mathpzc}{OT1}{pzc}{m}{it}

\newcommand{\subscript}[2]{$#1 _ #2$}

\DeclareSymbolFont{extraup}{U}{zavm}{m}{n}
\DeclareMathSymbol{\varheart}{\mathalpha}{extraup}{86}
\DeclareMathSymbol{\vardiamond}{\mathalpha}{extraup}{87}

\def\p{\partial}

\def\mb{\mathbb}

\def\l{\lrcorner}

\begin{document}
	\title{Logarithmic $\p\bp$-lemma and several geometric applications (with an Appendix joint with Sheng Rao)}
	\let\oldmaketitle\maketitle
	\renewcommand\maketitle{{\bfseries\boldmath\oldmaketitle}}	
	\begin{abstract}
		In this paper, we prove a $\p\bp$-type lemma on compact K\"ahler manifolds for logarithmic differential forms valued in the dual of a certain pseudo-effective line bundle, thereby confirming a conjecture proposed by X. Wan. 
		
		We then derive several applications, including strengthened results by   H. Esnault--E. Viehweg on the degeneracy of the spectral sequence at the $E_1$-stage  for projective manifolds associated with the logarithmic de Rham complex, as well as  by L. Katzarkov--M. Kontsevich--T. Pantev on the unobstructed locally trivial deformations of a projective generalized  log  Calabi--Yau pair  with some weights, both of  which are extended  to the broader context of compact K\"ahler manifolds.
		
		Furthermore, we establish the K\"ahler version of an injectivity theorem originally formulated by F. Ambro in the algebraic setting. Notably, while O. Fujino previously addressed the K\"ahler case, our proof takes a different approach by avoiding the reliance on mixed Hodge structures for cohomology with compact support.  
	\end{abstract}

\author{Runze Zhang}

	\address{School of Mathematics and statistics, Wuhan  University,
		Wuhan 430072, China; Universit\'e C\^ote d'Azur, CNRS, Laboratoire J.-A. Dieudonn\'e, Parc Valrose, F-06108 Nice Cedex 2, France}
	\email{runze.zhang@unice.fr}
\urladdr{\href{https://sites.google.com/view/runzezhang}{https://sites.google.com/view/runzezhang}}

	\dedicatory{A Tribute to Professor Kefeng Liu: On the Occasion of His 60th Birthday}

	\date{\today}
	

	\setcounter{tocdepth}{2}
	
		\maketitle
	
	\section{Introduction}
	\subsection{Backgrounds and main results}
	One of the strongest results in deformation theory is the \textit{Bogomolov--Tian--Todorov (BTT) theorem} \cite{Bog78,Tia87,Tod89}, which, using differential geometric methods, establishes that Calabi--Yau manifolds have unobstructed deformations.  This means that the Kuranishi space of deformations of complex structures on such manifolds is smooth. We should remark here that while many authors consider Calabi--Yau manifolds to be projective, we make no projectivity assumptions. Specifically, in this paper, a \textit{Calabi--Yau manifold} $X$ is defined as a \textit{compact K\"ahler} manifold  with  trivial \textit{canonical bundle}   $K_X \simeq \mathcal{O}_X.$
	Algebraic proofs of the BTT theorem have also been provided using the degeneration of the Hodge-to-de Rham spectral  sequence plus the nowadays called ${T}^1$-lifting technique  \cite{Ran92,Kaw92,FM99}, see also \cite{IM10}. In fact, the K\"ahler condition in this theorem can be relaxed. For instance, it applies to a compact complex manifold $X$  with  trivial canonical bundle that  satisfies the $\p\bp$-lemma, see e.g. \cite[the proof in Chapter 6 after minor changes]{Huy05} and \cite[Theorem 1.3]{Pop19}.  Recall that the  \textit{(standard) $\p\bp$-lemma} refers to: for every pure-type $d$-closed form on a compact complex manifold, the properties of $d$-exactness, $\p$-exactness, $\bp$-exactness, and $\p\bp$-exactness are equivalent. More generally, it also holds for such manifolds satisfying a weak version of the  $\partial\bar{\partial}$-lemma\footnote{The weak version of the standard  $\partial\bar{\partial}$-lemma was first introduced by J. Fu--S.-T. Yau in \cite{FY11} while studying deformations of balanced manifolds.}, as shown by the proofs in \cite{LRW19}, see also \cite{RZ21} for a CR-version. Moreover, the BTT theorem holds for such manifolds if the Fr\"olicher spectral sequence degenerates at the first page, see for example, \cite[the proof in Theorem 4.18]{KKP08} and \cite[Theorem 3.3]{ACRT18} (this case is even true when  $K_X$  is a torsion bundle, see  e.g. \cite[Corollary 4.1]{Iac17}).  
	
	Noteworthy to mention that the BTT theorem has various extensions, which include results on the unobstructedness of deformations of the following:
	
	\begin{enumerate}[(\subscript{D}{{\arabic*}})]
		\item \,\label{previous log CY}(generalized) log Calabi--Yau pairs
		 \cite{KKP08,San14,Iac15,Ran17,LRW19,Wan18};
		 
		 \vspace{0.3em}
		\item \,(weak) Fano varieties \cite{Ran92,Min01,San14,Ran17};
		
			 \vspace{0.3em}
		\item\,(weak) Poisson structures \cite{Got10,Hit11,FM12,Ran17,Ran19};
		
			 \vspace{0.3em}
		\item\,Landau--Ginzburg models \cite{KKP17};
		
			 \vspace{0.3em}
		\item  \,a pair $(X,L)$ (resp.  $(X,\mathcal{F})$), where $L$ (resp.  $\mathcal{F}$) is a line bundle (resp.  coherent sheaf) on a smooth projective variety $X$ over an algebraically closed field of characteristic zero, with a trivial canonical bundle \cite{IM19,IM21}.
	\end{enumerate}
	
		 \vspace{0.3em}
	Notice that direction \ref{previous log CY} is particularly interesting from a mirror symmetry perspective, as the deformations of a (generalized) log Calabi--Yau pair are expected to mirror to the deformations of the  corresponding complexified symplectic form on the mirror Landau--Ginzburg model, see \cite{Aur07,Aur09,KKP08,Kon} for further details. Additionally, the BTT theorem has applications in the context of logarithmic geometry, see for instance, \cite{Fri83,KN94,CLM23,FFR21,Fel22,FP22,Fel23,FPR23}.

	\vspace{0.6em}
	In this paper, we focus our attention on the direction \ref{previous log CY}, specifically considering the pair \( (X, D) \), where \( X \) is a compact K\"ahler manifold and \( D \) is a simple normal crossing divisor. One of our main theorems, stated below, strengthens all the corresponding results in this direction.
	
	\begin{bigthm}\label{main thm}
		Let $X$ be an $n$-dimensional  compact K\"ahler manifold,  and let $D=\sum_{i=1}^s D_i$ be a simple normal crossing divisor on $X$. Assume that there exist rational weights
		$\{\aaa_i\}_{1\le i\le s}\subset [0,1]\cap \mathbb{Q}$ such that
		\begin{equation}\label{KKP's setting}
			\sum_{i=1}^s\,\mathcal{O}_X(\aaa_i D_i)=-K_X\in\mathrm{Pic}(X)\otimes_\mb Z \mb Q=:\mathrm{Pic}_{\mb Q}(X),
		\end{equation}which means that 
		\begin{equation}\tag{1.1$^\prime$}\label{22222222222}
			\sum_{i=1}^s \,\mathcal{O}_X(N\aaa_iD_i)=-NK_X\in\mathrm{Pic}(X),
		\end{equation}
		for some positive integer $N$. 
		
		Then,  the  locally trivial (infinitesimal) deformations of the pair $(X,D)$ are  unobstructed.  That is, for such deformations, the pair admits a smooth  Kuranishi space  (see $\S$ \ref{subsection 5.4.1} for further details).
	\end{bigthm}
	
	\begin{specialremark}\label{remark 1.}\,\,\begin{enumerate}[{$\bullet$}]
			\item The \textit{locally trivial} deformations of a pair $(X, D)$ can be roughly understood as deformations of $X$ in which $D$ deforms along with it in a locally trivial manner, so in particular, keeping  the analytic singularity types. Therefore, considering locally trivial deformations is \textit{not} a restriction but rather has geometric significance.  Notice that if $D$ is smooth, then every deformation of the pair is locally trivial, see for example  \cite[Lemmata 4.3.4 \& 4.3.5 \& Theorem 4.4.3]{Man22}.  Consequently in this case, the pair $(X, D)$ in Theorem \ref{main thm} has unobstructed deformations.
			
			\medskip
			\item Note that if the \emph{locally trivial} assumption is dropped in the above theorem 
			(i.e., \ if partial smoothings of $D$ are allowed), the conclusion may fail, 
			as illustrated by the examples in \cite{FPR23} due to S. Felten--A. Petracci--S. Robin.
		\end{enumerate}
	\end{specialremark}
	
	\begin{specialremark}\label{remark 1.1}
		\begin{enumerate}[{($i$)}]
		\item \label{remark 1.1 (i)}Theorem \ref{main thm} for the case of projective manifolds was  proved  by L. Katzarkov--M. Kontsevich--T. Pantev  \cite[$\S$  4.3.3 \textbf{(iii)}]{KKP08} (see also \cite[$\S$ 3, Case 3]{Kon} for a special case). In the projective setting, the equality \eqref{KKP's setting} can be rephrased using the terminology and concepts from algebraic geometry as the $\mb Q$-divisor
		$$K_X+\sum_{i=1}^s\aaa_i D_i$$
		being \textit{$\mathbb{Q}$-trivial}.  Here we also denote by $K_X$ the \textit{canonical divisor} on a projective manifold $X$. 
		
		\medskip
	We give a brief overview of the proof in \cite{KKP08}, which proceeds in three main stages.
	\begin{enumerate}[{1. }]
		\item \textit{Degeneration property of a DGBVA.}\quad
		Katzarkov--Kontsevich--Pantev construct a differential graded Batalin--Vilkovisky algebra (\dgbva),
		together with an associated differential graded Lie algebra (DGLA)
		\((\mathfrak{g}, \mathrm d_{\mathfrak{g}})\).
		This DGLA contains a direct summand
		\((\mathfrak{g}', \mathrm d_{\mathfrak{g}'})\) that governs the locally trivial  deformations
		of the pair \((X,D)\).\footnote{Roughly speaking, the underlying philosophy is encapsulated in
			Deligne's principle from his 1986 letter to J. Millson:
			``\textit{in characteristic $0$, a deformation problem is controlled by a \dgla, and quasi-isomorphic \dgla s give rise to the same deformation theory}'', see \cite{Del86}.}
		They then prove that this \dgbva\ satisfies the \emph{degeneration property} via the theory of mixed Hodge structures, see \cite[Lemma 4.21]{KKP08}.
		
		\vspace{0.5em}
		\item \textit{Homotopy abelianity and formal smoothness.}\quad
		The degeneration property implies that \((\mathfrak{g}, \mathrm d_{\mathfrak{g}})\) is \emph{homotopy abelian},
		and hence the same holds for \((\mathfrak{g}', \mathrm d_{\mathfrak{g}'})\).
		As a consequence, the \emph{formal moduli space} associated with the \dgla\
		\((\mathfrak{g}', \mathrm d_{\mathfrak{g}'})\) is smooth.
		This establishes the \emph{formal smoothness} of the Kuranishi space, meaning that the completion
		of its local ring at the maximal ideal is isomorphic to a ring of formal power series.
		
		\vspace{0.5em}
		\item \textit{From formal to analytic smoothness via Artin's theorem.}\quad
		To conclude that the Kuranishi space \(B\) is smooth as an analytic germ,
		one can appeal to M. Artin's fundamental approximation theorem \cite{Art68}  on the existence of 	convergent solutions for analytic equations  from the existence of  formal solutions, 
		which implies that this formally smooth analytic germ is in fact analytically
		smooth, see e.g.  \cite[Theorem B.1]{RS20} for a detailed proof.
	\end{enumerate}

			\vspace{0.7em}
			\item  Notably, the setup \eqref{KKP's setting} encompasses two important  endpoint situations (the corresponding projective cases were also presented in \cite[$\S$ 4.3.3 \textbf{(i)} \& \textbf{(ii)}]{KKP08}, respectively), namely:
			\begin{enumerate}[{(a)}]
				\item\label{!1}   \textit{the log Calabi--Yau case}, where the logarithmic canonical  bundle $$\Omega_{X}^n(\log D)\simeq K_X\otimes\mathcal{O}_X(D)$$ is trivial (i.e., all $\aaa_i=1$  and $N=1$ in \eqref{22222222222});
			\medskip
				\item \label{!2}\textit{the Calabi--Yau  case with  a simple normal crossing divisor},  where $K_X$ is trivial (i.e., all  $\aaa_i=0$ and $N=1$).
			\end{enumerate}
		
			\vspace{0.6em}
			Therefore, we can  regard  the pair $(X,D)$ in Theorem \ref{main thm} as a \textit{K\"ahler generalized  log Calabi--Yau pair}.
			
			\vspace{0.6em}
			Based on the theory of   DGLA  and the Cartan homotopy construction, D. Iacono proved cases \ref{!1} and \ref{!2} in the algebraic setting  over an algebraically closed field of characteristic zero \cite[Corollaries 5.5 \& 5.8 \& Remark 5.6]{Iac15} (see also \cite[$\S$ 4.2]{Iac17} from the perspective of abstract BTT theorem), and later extended the case \ref{!1}  to the context of the compact complex manifolds where  the  logarithmic Hodge-to-de Rham spectral sequence degenerates at the $E_1$-level 
			\cite[Corollary 4.5]{Iac17}. Additionally, see K. Sano's work employing the $T^1$-lifting technique \cite[Remark 2.5]{San14} and Z. Ran's approach via Poisson deformations \cite[Theorem 12]{Ran17}. More recently,  K. Liu--S. Rao--X. Wan proved cases \ref{!1} and \ref{!2} using  a power series method, more faithfully following Tian--Todorov's approach  \cite[Theorems 0.8 \& 0.9]{LRW19}.

			\vspace{0.3em}
			Furthermore, Iacono proved the case (also in the algebraic setting) if any $\aaa_i$ is equal to $\frac{1}{N}$ for some positive integer $N$ in \eqref{KKP's setting}, via a cyclic covering construction, see \cite[Proposition 6.4]{Iac15}. Wan also contributed to this line of research by settling the case of $N=2$ (when $X$ is a compact K\"ahler manifold), with his proof also relying on a cyclic cover trick, as shown in \cite[Theorem 0.5]{Wan18}.
		\end{enumerate}
	\end{specialremark}
	
Indeed, inspired by the approach in \cite{KKP08} (as outlined in \textit{Remark \ref{remark 1.1}} \ref{remark 1.1 (i)}), 
we provide two distinct proofs of Theorem \ref{main thm}:

\begin{enumerate}
	\item[--] a proof based on DGLA theory, following the more \textit{algebraic and abstract}
	framework of \cite{KKP08}, presented in \S\,5.4.3;
	
	\medskip
	\item[--] a purely analytic proof using power series method in the spirit of \cite{LRW19},
	which is more \textit{visualizable and geometric}. The proof, jointly with Sheng Rao,  is given
	in Appendix~\ref{jw}.
\end{enumerate}
In both cases, the key step is to solve a certain $\bar{\partial}$-equation for logarithmic
forms twisted by the dual of a pseudo-effective line bundle, originally conjectured by
Wan \cite{Wan18}.

	\begin{bigthm}[{\cite[Conjecture 0.8]{Wan18}}] \label{n,q-main} 
			Let $(X,D)$ be as in Theorem \ref{main thm},  and let $L$ be a holomorphic line bundle over $X$. 
			Assume that there exist rational weights
			$\{\qqq_i\}_{1\le i\le s}\subset [0,1]\cap \mathbb{Q}$ such that
			\begin{equation}\label{rational number}
				L=\sum_{i=1}^s \mathcal{O}_X(\qqq_i D_i)\in\mathrm{Pic}_{\mathbb{Q}}(X).
			\end{equation}
			Then for any $
			\alpha\in A^{0,q}\bigl(X,\Omega_X^{p-1}(\log D)\otimes L^{*}\bigr)$
			satisfying
			\[
			\bar{\partial}D^\prime_{h^{*}}\alpha = 0
			\quad \text{pointwise on } X^\circ := X\setminus \mathrm{Supp}\, D,
			\]
			the following twisted logarithmic $\bar{\partial}$-equation
			\begin{equation}\tag*{$(\curlyvee)$}\label{solution 2}
				\bar{\partial}\chi
				=
				D^\prime_{h^{*}}\alpha
				\qquad \text{pointwise on } X^\circ
			\end{equation}
			admits a solution $
			\chi\in A^{0,q-1}\bigl(X,\Omega_X^{p}(\log D)\otimes L^{*}\bigr).$
			
				\medskip
			Here  $D^\prime_{h^*}$ represents the $(1,0)$-part of the integrable logarithmic connection $\nabla_{h^*}$ along $D$, induced by the singular metric $h^*$ on the dual bundle $L^*$ (see $\S$ \ref{Section 2} for more details).
	\end{bigthm}

	\begin{specialremark}\label{remark 1.2}\,\,Theorem \ref{n,q-main} revisits two extreme cases: \begin{equation*}
		 L\simeq\mathcal{O}_X\quad\text{and}\quad
		 L\simeq \mathcal{O}_X(D).
		\end{equation*}  Both of which were proved by Liu--Rao--Wan \cite[Theorem 0.1 \& 0.2]{LRW19}. It also revisits the case where all $\qqq_i=\frac{1}{N}$ for some positive integer $N$, and $D$ is a smooth divisor, as proved by Wan \cite[Theorem 0.1]{Wan18}. See also \cite[$\S$ 4]{RZ24} for  some discussions related to Wan's conjecture from the perspective of double complexes.
	\end{specialremark} 
	
	\begin{specialremark}\,\,
		 Strictly speaking, Theorem \ref{n,q-main} is weaker than the ``genuine'' logarithmic version of the (standard)  $\p\bp$-lemma aforementioned, as we cannot guarantee that the solution $\chi$ obtained in  \ref{solution 2} is moreover $D^\prime_{h^{*}}$-exact.\footnote{Following the terminology of \cite[Notation 3.5]{RZ18} (see also \cite{RZ22} for the foliated setting), we can analogously say that the pair $(X,D)$ satisfies $\mathbb S^{p,q}$ with respect to $(L^*,h^*)$.} It would be very interesting to know if this can be achieved.  
	\end{specialremark}
	
	\vspace{0.5em}
	As is widely recognized, the $E_1$-degeneration of certain spectral sequences is a useful tool in algebraic geometry and complex geometry, such as being used to imply several injectivity, vanishing and torsion-free theorems. For a more comprehensive discussion, we refer the reader to \cite{EV92, deFEM, Fuj-Book}, along with the references therein.

		\vspace{0.5em}
	Following the approach in \cite[Theorem 3.2]{LRW19} and utilizing the logarithmic counterpart of the general description of terms in the Fr\"olicher spectral sequence as presented in \cite[Theorems 1 \& 3]{CFGU97}, we obtain the expression expression $E_r^{p,q}\simeq \frac{Z_r^{p,q}}{B_r^{p,q}}$  with the differential map $ {\zd_r^{p,q}}:E_r^{p,q}\rightarrow E_r^{p+r,q-r+1}$. The subgroups  $Z_r^{p,q}$ and $B^r_{p,q}$ are detailed in $\S$ \ref{subsection:E1}. 
	As a direct corollary of Theorem \ref{n,q-main}, we deduce that $\zd_r^{p,q} = 0$ for all $r \geq 1$ and all $p,q$. This leads to the following result, which extends the work of H. Esnault--E. Viehweg on projective manifolds \cite[Theorem 3.2 \& Remarks 3.3]{EV92} to the broader setting of compact K\"ahler manifolds.
	\begin{bigthm}\label{Thm 1.4}
		With the same setting \eqref{rational number} as in Theorem \ref{n,q-main}, the following spectral sequence
		$$E_1^{p,q}=H^q(X,\Omega_{X}^p(\log D)\otimes L^*)\Longrightarrow\mb H^{p+q}(X,\Omega_{X}^\bullet(\log D)\otimes L^*)$$
		associated to the logarithmic de Rham complex (see $\S$ \ref{Section 2})
		$$(\Omega_X^\bullet(\log D)\otimes L^*, \nabla_{{h^*}_\bullet})$$
		degenerates at the $E_1$-level.  Here $\mb H^{p+q}(X,\Omega_{X}^\bullet(\log D)\otimes L^*)$ denotes the hypercohomology. 
	\end{bigthm}
	
	\vspace{1em}
	Building on this, we establish the K\"ahler version of Ambro's injectivity theorem \cite[Theorem 2.3 \& Remark 2.6]{Amb14}, originally proved in the algebraic setting for proper, non-singular varieties over an algebraically closed field of characteristic zero, see also the recent work of T. Murayama \cite{Mur24}. The proof of Theorem \ref{inj-} is detailed in $\S$   \ref{subsection:inj}.   Notably, the case of Theorem \ref{inj-} where  $X$  is a compact complex manifold in Fujiki's class \( \mathscr{C} \) (i.e., bimeromorphic to a compact K\"ahler manifold) was established by Fujino \cite[Theorem 1.2]{Fuj17}, utilizing the theory of mixed Hodge structures for cohomology with compact support. Our approach is quite different from Fujino's (see Remark \ref{difference}).
	\begin{bigthm}\label{inj-}
		Let $(X,D)$ be as in Theorem \ref{main thm},   and let $L$ be a holomorphic line bundle over $X$ such that  \begin{equation}\label{satis}
			L=K_X\otimes\sum_{i=1}^{s}\mathcal{O}_X(b_iD_i)\in\mathrm{Pic}_\mb Q(X)\,\,\,\,\,\textit{with}\,\,b_i\in(0,1]\cap\mb Q\,\,\textit{for all }i. 
		\end{equation}
		Then the restriction homomorphism
		\begin{equation}\label{VVV}
			H^q(X,L)\stackrel{{j}}{\longrightarrow} H^q(X^\circ,L|_{X^\circ})
		\end{equation}
		is  injective, for all $q$. Equivalently, for every effective Cartier divisor $\widehat D$ with $\textrm{Supp}\,\widehat D\subset\textrm{Supp}\,D,$  the natural homomorphism
		\begin{equation}\label{XXX}
			H^q(X,L)\stackrel{j^\prime}{\longrightarrow} H^q(X,L\otimes\mathcal{O}_X(\widehat D))
		\end{equation}
		induced by the inclusion $\mathcal O_X\subset\mathcal{O}_X(\widehat D)$  is injective, for all $q$.
	\end{bigthm}

	\vspace{0.8em}
	In addition to Theorems \ref{main thm}, \ref{Thm 1.4}, and \ref{inj-}, we present another application of Theorem \ref{n,q-main}: the closedness of twisted logarithmic forms, see $\S$ \ref{subsection:clo}.

	\subsection{Idea of the proofs}
	We briefly outline the basic strategies behind the proofs of Theorem \ref{n,q-main} and Theorem \ref{main thm}.
	\vspace{0.5em}
	
	\noindent\faHandPointRight[regular]\,\,\,\,For Theorem \ref{n,q-main}, let us  first consider the special case: \(0 < \qqq_i \leq 1\) for every \(1 \leq i \leq s\), which corresponds to Theorem \ref{Thm 1.3}. The first observation is that any \(\alpha \in A^{0,q}(X, \Omega_{X}^p(\log D) \otimes L^*)\), and consequently \(D^\prime_{h^*} \alpha\),  is in fact \textit{smooth in  conic sense} (Proposition \ref{Lemma:>0=>sc}). By utilizing the \textit{Hodge decomposition for forms that are smooth in  conic sense}, as established by J. Cao--M. P\u{a}un in \cite{CP23b}, we can deduce that \([D^\prime_{h^*} \alpha]_{\bp} = 0\) in the $L^*$-valued $(p,q)$-conic Dolbeault cohomology group (Notation \ref{notation co}). Finally, to ensure the existence of a solution  $\chi$  (with at most logarithmic poles) in equation \eqref{solution 1}, we require an \textit{acyclic resolution} of the sheaf  $\Omega_X^p(\log D) \otimes L^*$  by sheaves of germs of  $L^*$ -valued  $(p, \bullet)$-forms that are smooth in conic sense (Proposition \ref{prop-resolution a}).
	
	\vspace{0.67em}
	In the more general setting of Theorem \ref{n,q-main}, \(D^\prime_{h^*} \alpha\) is no longer smooth in  conic sense; however, the good news is that it remains a \textit{conic current with values in $(L^*,h^*)$}, denoted by \(T_{D^\prime_{h^*} \alpha}\). We then decompose \(D\) into \(E\) and \(F\) \eqref{D=E+F}, where, roughly speaking, \(E\) (resp.  \(F\)) represents the \(\qqq_i = 0\) part (resp.  \(\qqq_i > 0\) part). By using  the \textit{de Rham--Kodaira decomposition for conic currents} \cite{CP23b}, we can decompose \(T_{D^\prime_{h^*} \alpha}\) into two parts. One part falls in the image of \(\bp\), and the other one is a \textit{residue term}.
	
	\vspace{0.67em}
	Inspired by \cite[Lemma 2.3]{LRW19} and also \cite[Lemma 4.3]{CP23b}, the key Lemma \ref{key lemma} shows that the residue term can also be expressed as \(\bp \widehat{T}\), where \(\widehat{T}\) is a conic current, modulo the space of \textit{on-\(E\) conic currents valued in \((L^*,h^*)\)} on \(X\) (where, roughly speaking, any element in this space annihilates all test forms valued in \(L\) that are smooth in  conic sense and moreover vanish on \(E\); consequently, they will vanish on \(D\) as shown in \eqref{!@!}). This concept is motivated by the work of J. R. King \cite{Kin83}.
	The final step to find the solution \(\chi\) in equation \ref{solution 2} follows a similar approach. The main difference is that this time we need to replace the protagonists in the corresponding acyclic resolution (Proposition \ref{resolution-2}) with the \textit{sheaf of log-\(E\) conic currents with values in  $(L^*,h^*)$}, which is defined as the quotient sheaf of the sheaf of $(L^*,h^*)$-conic currents by the sheaf of on-\(E\) conic currents  valued in $(L^*,h^*)$, see Definition \ref{LoGG}. 
	
\begin{specialremark}\,\,
	\begin{enumerate}
	[{$\bullet$}]
	\item  Note that Theorem \ref{Thm 1.3} can also be derived using the conic current approach, see Corollary \ref{corooooo}. Moreover, our method provides part of a new proof of \cite[Theorem 0.1]{LRW19}, see Remark \ref{R} \ref{RR2}.
	
	\medskip
	\item   
	We would also like to emphasize that, much recently, a general logarithmic type
	$\p\bp$-lemma (for top-degree) proved by Cao--P\u{a}un \cite[Theorem 1.1]{CP23b} has played a
	significant role in proving Fujino's conjecture on a K\"ahler injectivity theorem \cite[Theorem 1.2]{CP23b}
	(see also the recent independent work of T. O. M. Chan--Y.-J. Choi--S. Matsumura
	\cite{CCM23}).
	Moreover, their logarithmic $\p\bp$-lemma provides new results concerning the	invariance of plurigenera in the K\"ahler setting \cite{CP23a}.  
	
Our goal is also to establish a logarithmic  $\p\bp$-lemma, but our approach differs from that of \cite{CP23b}, especially in the construction of the acyclic resolution, while still relying on the conic Hodge theory they developed.

	\end{enumerate}
\end{specialremark}

\vspace{0.5em}
\noindent\faHandPointRight[regular]\,\,\,\,For Theorem \ref{main thm}, 

\begin{enumerate}
	\item[--](\textit{Algebraic approach})\quad
	As noted in Remark \ref{remark 1.1}, the key requirement is the \textit{degeneration property} (Definition \ref{1212}) for a suitable DGBVA. Using Theorem \ref{n,q-main} (or more practically, Theorem \ref{Thm 1.4}), we can achieve this by considering the DGBVA \( (\AAA, \ddd, \bfitDelta) \), where
	\[
	\AAA := A^{0, \bullet}\bigl(X,\wedge^{\bullet} T_X^1(-\log D)\bigr), \qquad 
	\ddd := \bar\partial,
	\]
	and
	\[
	\bfitDelta := \mathbf{i}_{\Upomega}^{-1} \circ D^\prime_{h^*} \circ \mathbf{i}_{\Upomega}.
	\]
	Here \(T_X^1(-\log D)\) denotes the logarithmic tangent bundle (Definition \ref{tangent}). The map  
	\[
	\mathbf{i}_{\Upomega}: \wedge^{\bullet} T_X^1(-\log D) \longrightarrow 
	\Omega_X^{n-\bullet}(\log D) \otimes L^*
	\]
	is the isomorphism given by contraction with a nowhere‑vanishing section
	\[
	\Upomega \in A^{0,0}\bigl(X,\Omega^n_X(\log D)\otimes L^*\bigr),
	\]
	where the holomorphic line bundle
	$$L := \Omega_X^n(\log D) \simeq K_X \otimes \mathcal{O}_X(D)$$ satisfies 
	\[
	L = \sum_{i=1}^s \mathcal{O}_X( \qqq_i D_i) \in \operatorname{Pic}_{\mathbb Q}(X), \qquad 
	\qqq_i := 1 - \aaa_i \in [0,1].
	\]
	The essential point is that only the \textit{existence} of logarithmic solutions is required.
	
	\medskip
	\item[--](\textit{Analytic approach})\quad
	Following Y. Kawamata \cite{Kaw78}, the Kuranishi family \(\mathscr{K}\) of locally trivial deformations of the pair \((X,D)\) (Definition \ref{Kuranishi-}) can be realized concretely inside the space of logarithmic Beltrami differentials
	\[
	\Gamma_{\text{real analytic}}\bigl(X, T_X^1(-\log D)\otimes \Lambda^{0,1}T^*X\bigr)
	\subset A^{0,1}\bigl(X,T_X^1(-\log D)\bigr),
	\]
	whose elements satisfy the integrability condition
	\[
	\bar\partial\varphi = \tfrac12[\varphi,\varphi].
	\]
	By constructing \textit{explicit} solutions to the twisted logarithmic \(\bar\partial\)-equation (Lemma \ref{lemma-exp}), we  follow the iterative method in \cite{LRW19} to produce a \textit{family of integrable logarithmic Beltrami differentials} over a sufficiently small disk. This yields the unobstructedness in a manner that aligns closely with the classical Tian--Todorov philosophy.
\end{enumerate}
	\subsection{Overview and outlook}

	We prove Theorem \ref{n,q-main} via an analytic approach, and Theorem \ref{main thm} is an application of it, generalizing \cite[$\S$ 4.3.3 \textbf{(iii)}]{KKP08}, and consequently \cite[$\S$ 4.3.3 \textbf{(i)} \& \textbf{(ii)}]{KKP08}, to the context of compact K\"ahler manifolds. Indeed, by constructing again a  DGBVA and employing the mixed Hodge theory,  Katzarkov--Kontsevich--Pantev provided another generalization of the previous results, that is  \cite[$\S$ 4.3.3 \textbf{(iv)}]{KKP08}, specifically when $X$ is a \textit{projective normal-crossing  Calabi--Yau}. More precisely, assume that $X$ is a strict normal crossings variety with irreducible components $X=\bigcup_{j\in J}X_j$ equipped with a holomorphic volume form $\Omega_X$ on $X-X_{\text{sing}}$. This form $\Omega_X$ satisfies the condition that its restriction to each $X_j$ has a logarithmic pole along $X_j\cap(\cup_{k\not=j}X_k$), and the residues of these restricted forms cancel along each $X_j\cup X_k$. Motivated by this, the following question naturally arises, which we will consider in a forthcoming paper.
	
	\begin{questionno}
		Could further analytic methods be developed to generalize  the variety $X$ in Theorem \ref{n,q-main} to the setting of normal-crossing Calabi--Yau varieties that are K\"ahler but not necessarily projective? If so, it may also be promising to extend the variety $X$ in Theorem \ref{main thm} to this broader context without relying on  the mixed Hodge theory.
	\end{questionno}

	\noindent
	\textbf{Notations and conventions.} 
	Throughout this paper, we work over the field of complex numbers.
	\begin{itemize}
		\item[--]Any (compact) complex manifold $X$ in this paper is assumed to be connected. 
		\item[--] Denote by $\mathscr A _X$ (resp.  $\mathcal O_X$) the sheaf of germs of $\mathscr{C}^\infty$ differentiable functions (resp.  holomorphic functions) over $X$.
		\item[--]\textit{Locally free sheaves of $\mathscr A _X$-modules (resp.  $\mathcal O_X$-modules)} and  \textit{$\mathscr{C}^\infty$ complex (resp.  holomorphic) vector bundles}  are considered synonymous.
		\item[--]The terminology \textit{Cartier divisors}, \textit{invertible sheaves}, and \textit{holomorphic line bundles} are used interchangeably. 
		\item [--] A sheaf $\mathcal F$ is called \textit{flabby} if for every open subset $V$ of $X$, the restriction map $\mathcal{F}(X)\rightarrow \mathcal{F}(V)$ is onto, i.e., if every section of $\mathcal F$ on $V$ can be extended to $X.$ 
		\item [--] A flabby sheaf $\mathcal F$ is  \textit{acyclic} on all open sets $V\subset X$, meaning that $H^q(V,\mathcal{F})=0$ for any $q\geq 1$.
		\item[--] We  use additive notation for tensor products and powers of line bundles, and multiplicative notation (resp.  additive notation) for hermitian metrics (resp.  local weights) of line bundles.  For example, $(L_1,h_1,\varphi_1), (L_2,h_2,\varphi_2),$ and $(L_1 + L_2,  h_1\cdot h_2,\varphi_1+\varphi_2)$. Here $h_1=e^{-\varphi_1}$ (resp.  $h_2=e^{-\varphi_2}$) is a (possibly singular) metric on $L_1$ (resp.  $L_2$).
		\item[--] For $\aleph\in\mb Q$, $\lfloor\aleph\rfloor$ denotes the  \textit{integral part} of $\aleph$, defined as the only integer such that 
		$\lfloor\aleph\rfloor\leq\aleph\leq\lfloor\aleph\rfloor+1$.
	\end{itemize}
	
	\vspace{0.9em}
	This paper is organized as follows.
	
	\setcounter{tocdepth}{3}
	\tableofcontents
	\textbf{Acknowledgement}\quad This paper forms part of the author's PhD thesis.  He expresses his deepest gratitude to his advisors, Professors Junyan Cao and Sheng Rao, for their insightful guidance and constant encouragement. He is also sincerely thankful to Professors Tsz On Mario Chan, Ya Deng, Simon Felten, Andreas H\"{o}ring, Donatella Iacono, Kefeng Liu, Shin-ichi Matsumura, Mihai P\u{a}un, Helge Ruddat, and Xueyuan Wan for numerous helpful discussions and valuable comments on this work. Special thanks are due to Donatella Iacono and Simon Felten for their kind responses to his email inquiries, and to Xueyuan Wan for pointing out a serious mistake  in an earlier version of this manuscript. 
	
	\medskip
The author is further indebted to Professors Philippe Eyssidieux, Andreas H\"{o}ring,
Kefeng Liu, and Mihai P\u{a}un for their patient attention and perceptive questions
during the defense of his thesis. In particular, a question raised by Mihai P\u{a}un
provided the original motivation for Appendix~\ref{jw}. The author is especially
grateful to Professor Sheng Rao for many fruitful discussions, which ultimately
led to the joint development of Appendix~\ref{jw}. Part of this work was completed
during a visit to the Shanghai Center for Mathematical Sciences. The author warmly
thanks Professor Chen Jiang for his gracious hospitality and for providing an
excellent research environment.
	
	\medskip
	Finally, the author sincerely acknowledges the reviewers for their thoughtful comments and suggestions, which have greatly improved the quality of this paper. This work was supported by the China Scholarship Council (Grant No.~202306270252).
	\section{Logarithmic connection and logarithmic complex}\label{Section 2}
	In this section, we will recall some basic notions regarding the   logarithmic connection (induced by the singular metric on  a holomorphic line bundle) and the logarithmic de Rham complex. For more details refer to   \cite[Chapter 3.5]{GH78} and \cite[Chapter 2]{EV92}.

	Let $X$ be a compact complex manifold of dimension $n$ and $D=\sum_{i=1}^sD_i$ be a simple normal crossing divisor ( i.e., a divisor with non-singular components $D_i$ intersecting each other transversally) on $X$. Let
	\begin{equation*}\label{meromorphic sheaf}
		\Omega_{X}^p(\star D)=	\varinjlim_{v}\Omega_X^p({v}D)
	\end{equation*}
	be the sheaf of germs of $p$-meromorphic forms which are holomorphic on $X-D$ but possibly with arbitrary orders of poles along $D$. Obviously, $(\Omega_X^\bullet(\star D),d)$ is a complex. Deligne introduced the \textit{sheaf of germs of logarithmic $p$-forms} \cite{Del69}
	\begin{equation*}
		\Omega_{X}^p(\log D),
	\end{equation*}
	which is defined as the subsheaf of $\Omega_{X}^p(\star D)$ with logarithmic poles along $D$, i.e., if $V\subset  X$ is open, then
	$$\Gamma(V,\Omega_X^p(\log D))=\left\{\alpha\in\Gamma(V,\Omega_{X}^p(\star D))\,|\,\alpha\text{ and }d\alpha\text{ both have simple poles along } D\right\}.$$
	It turns out that there exists a subcomplex $(\Omega_{X}^\bullet(\log D),d)\subset(\Omega_{X}^\bullet(\star D),d)$, see for example \cite[II, 3.1-3.7]{Del70} or \cite[Properties 2.2]{EV92}. Furthermore, 
	$$\Omega^p_X(\log D)=\bigwedge^p\Omega_{X}^1(\log D)$$
	is locally free. More precisely, for any $z \in X$, suppose $z \in D_i$ for any $1 \leq i \leq d$ and $z \notin D_i$ for $d+1\leq i\leq  s$. We can then choose local coordinates $\{ z^1, \ldots, z^n \}$ in a small neighborhood $V$ of $z = (0, \ldots, 0)$ such that $D_i \cap V = \{ z^i = 0 \}$ for $1 \leq i \leq d$. One writes 
	\begin{equation*}
		\delta_j=\begin{cases}
			\frac{dz^j}{z^j}\quad\text{if}\quad j\leq d;\\
			dz^j\quad\text{if}\quad j>d,
		\end{cases}
	\end{equation*}
	and for $J=\left\{j_1,\ldots,j_p\right\}\subset\left\{1,\ldots,n\right\}$ with $j_1<j_2<\ldots<j_p,$
	$$\delta_J=\delta_{j_1}\wedge\ldots\wedge\delta_{j_p}.$$
	Then \begin{equation}\label{hhhhhh}
		\left\{\delta_J\,|\,\,\# J=p\right\}
	\end{equation} forms a basis of  $\Omega_X^p(\log D)$ as a free $\mathcal{O}_X$-module over $V$. Furthermore,  we denote by 
	$$\mathscr A ^{0,q}(\Omega_{X}^p(\log D))$$
	\textit{the sheaf of germs of $(0,q)$-forms valued in $\Omega_{X}^p(\log D)$}, which  is a locally free sheaf of $\mathscr A _X$-modules.
	Elements in 
	$$A^{0,q}(X,\Omega_{X}^p(\log D)),$$ i.e., the global sections of $\mathscr A ^{0,q}(\Omega_{X}^p(\log D))$, are called  \textit{logarithmic $(p,q)$-forms}. 
	
	Next, we recall the definition of (integrable) logarithmic connection along $D$ with respect to a holomorphic vector bundle.
	\begin{definition}[{\cite[Definition 2.4]{EV92}}] \label{log connection def}Let $\mathcal{E}$ be a locally free sheaf of $\mathcal{O}_X$-modules and let
		$$\nabla:\mathcal{E}\longrightarrow\Omega_X^1(\log D)\otimes\mathcal{E}$$
		be a $\mb C$-linear map satisfying
		$$\nabla(f\cdot e)=f\cdot \nabla(e)+df\otimes e.$$
		One defines
		$$\nabla_p:\Omega_{X}^p(\log D)\otimes\mathcal{E}\longrightarrow\Omega_X^{p+1}(\log D)\otimes\mathcal{E}$$
		by the rule
		$$\nabla_p(\omega\otimes e)=dw\otimes e+(-1)^p\omega\wedge\nabla(e).$$
		We assume that $$\nabla_{p+1}\circ\nabla_p=0.$$  Such $\nabla$ will be called an \textit{integrable logarithmic connection along $D$}, or just a connection. The complex
		$$(\Omega_{X}^\bullet(\log D)\otimes\mathcal{E},\nabla_\bullet)$$
		is called the \textit{logarithmic de Rham complex} of $(\mathcal{E},\nabla)$.
	\end{definition}
	
	Let $L$ be a holomorphic line bundle over $X$ satisfying 
	\[ L = \sum_{i=1}^s \mathcal{O}_X( \qqq_i D_i) \in \textrm{Pic}_{\mathbb{Q}}(X) \]
	with $\qqq_i \in \mathbb{Q}$ (note that here we do not restrict the sign of $\qqq_i$). Then there naturally exists a singular metric $h:=h_L =e^{-\varphi_L}$ on $L$, where $\varphi_L$ is a collection of functions defined on small open sets, called the \textit{local weight}, locally can be written as
	\[ \varphi_L =\sum_{i=1}^s \qqq_i \log|z^i|^2, \]
	with $z^i=0 \, (i=1,\ldots,s)$ representing the local equations of components of $D$, see  e.g.  \cite[$\S$ 2]{Dem92}. We then obtain that the curvature current with respect to $h$ is 
	\begin{equation*}\label{curvatureaa} 
		i\Theta_{h}(L)=2\pi i\sum_{i=1}^s \qqq_i[D_i],
	\end{equation*}
	thanks to the Lelong--Poincar\'e formula, where $[D_i]$ is the current of integration over the irreducible $(n-1)$-dimensional analytic set $D_i$ for any $i$. In particular, $L$ is \textit{pseudo-effective} if all $\qqq_i$ are non-negative. 
	One can  verify without difficulty that there exists a (global)  integrable logarithmic connection along $D$ induced by the metric $h$, denoted by $\nabla_{h}$, which has a decomposition
	\[ \nabla_{h} = D^\prime_{h} + \bar{\partial}. \]
	 Here $D^\prime_{h}$ is the $(1,0)$-part of $\nabla_{h}$, which  has the local expression
	\begin{equation*}\label{def of 1,0}
		D^\prime_{h} := \partial - \partial \varphi_L=\partial-\sum_{i=1}^s \qqq_i\frac{dz^i}{z^i}.
	\end{equation*}
	More explicitly, following the notations as in  Definition \ref{log connection def}, we have
	\begin{equation}\label{`}
		\nabla_{h}(\omega\otimes e)=d\omega\otimes e+(-1)^{p}\omega\wedge(\sum_{i=1}^s \qqq_i\frac{dz^i}{z^i})\otimes e.
	\end{equation}
	One then can check that
	$$\nabla_{h}^2=0.$$

	\section{Warm-up: all  {$\qqq_i>0$}} 
	As a warm-up, we first consider a simple case of Theorem \ref{n,q-main}. In this section, we are devoted to proving:
	\begin{theorem}\label{Thm 1.3}
		Under the same settings as in Theorem \ref{n,q-main}, assume further that all rational numbers $\qqq_i$ in \eqref{rational number} lie in  $(0,1]$.
		Then, for any $\alpha\in A^{0,q}(X,\Omega_X^{p}(\log D)\otimes L^{*})$  satisfying  $\bp D^\prime_{h^*}\alpha=0$ pointwise on $X^\circ,$ the logarithmic $\bp$-equation:
		\begin{equation}\label{solution 1}
			\bp \chi=D^\prime_{h^*}\alpha\,\,\,\,\text{pointwise on }X^\circ 
		\end{equation}
		has a solution $\chi\in A^{0,q-1}(X,\Omega_X^{p}(\log D)\otimes L^{*})$.
	\end{theorem}
	\subsection{Hodge decomposition: conic version}
	In this subsection, we will state the Hodge decomposition for metric with conic singularities, following \cite[$\S$ 2.2 \& 2.3]{CP23b}.

	The \textit{setting} of this subsection is as follows.  Let $X$ be a compact K\"ahler manifold with dimension $n$, and let $D=\sum_{i=1}^s D_i$ be a simple normal crossing  divisor on $X$. Let $L$  be a holomorphic line bundle over $X$ admitting  a singular metric $$h:=h_L=e^{-\varphi_L}$$ which has \textit{logarithmic poles}, i.e., its local weight can be written as
	\begin{equation}\tag{$\varheart$}\label{weight}
		\varphi_L=\sum_{i=1}^s \qqq_i \log|z^i|^2+\varphi_{L,0},
	\end{equation}
	where $\qqq_i\in\mb Q$, $z^i=0$ represent the local equations of components of $D$ and $\varphi_{L,0}$ is a smooth function. Note that the $\qqq_i$ are not necessarily positive.
	The condition \eqref{weight} implies the corresponding curvature current is given by
	\begin{equation}\label{curvature}
		i\Theta_{h}(L) = 2\pi i \sum_{i=1}^s \qqq_i [D_i] + \theta_L,
	\end{equation}
	where $\theta_L$ is a smooth function on $X$. Choose a positive integer $m$ such that 
	\begin{equation}\tag{$\spadesuit$}\label{restriction2}
		\textit{for each } \qqq_i\in\mb Q \setminus \mb Z, \,\,\mmm \qqq_i\in \mb Z \,\,\textit {  and  }\,\,\lfloor \qqq_i-\frac{1}{\mmm}\rfloor=\left\lfloor \qqq_i\right\rfloor \textit{ both hold true.}
	\end{equation} One then denotes by $\omega_{\mathrm c}$ a metric (on $X^\circ=X\setminus \text{Supp}\, D$) with \textit{conic singularities} along the $\mb Q$-divisor
	$$D_\mmm:=\sum_{i=1}^s(1-\frac{1}{\mmm})D_i.$$
	By this we mean that if $z^1\cdots z^s=0$ is the local equation of the divisor $D$, then 
	\begin{equation}\label{omega_C}
		\omega_{\mathrm c}\simeq \sum_{j=1}^s\frac{i dz^j\wedge d\bar{z}^j}{|z^j|^{2-\frac{2}{\mmm}}}+\sum_{j=s+1}^n i dz^j\wedge d\bar{z}^j=:\omega_{\textrm{model}},
	\end{equation}
	that is to say, $\omega_{\mathrm c}$ is \textit{quasi-isometric} with $\omega_{\mathrm{model}}$, i.e.,
	$$C^{-1} \cdot\omega_{\mathrm{model}}\leq\omega_{\mathrm c}\leq C\cdot\omega_{\mathrm{model}}$$
	for some constant $C>0$.
	Notice that $\omega_{\mathrm c}$ is a closed positive $(1,1)$-current (smooth away from the support of $D$) on $X$. For the existence and the explicit constructions of such metric, refer to e.g. \cite[Proposition 2.1]{Cla08}.
	Now let  
	\begin{equation*}
		\left(V_i;z_i^1,\ldots,z_i^n\right)_{i\in \mathfrak{I}}
	\end{equation*}
	be a finite cover with coordinate charts such that 
	\begin{equation*}\label{equation of D}
		z^1_i\cdots z_i^s=0
	\end{equation*}
	is the local equation of the divisor $D$ when restricted to the set $V_i$.  We next consider the local ramified maps
	\begin{equation*}
		\pi_i: U_i\rightarrow V_i,\,\,\,\,\,\,\, \pi_i(w_i^1,\ldots,w_i^n):=((w_i^1)^\mmm,\ldots,(w_i^s)^\mmm,w_i^{s+1},\ldots,w_i^n).
	\end{equation*}
	It defines the orbifold structure corresponding to $(X, D_\mmm)$.
	
	Then, introducing the following definition is natural in our context, as it accounts for the singularities of the metric $h$ on $L$. This can be seen as a generalization of the usual definition of ``orbifold differential forms'', see e.g. \cite[$\S$ 5.4]{MM07} and \cite[$\S$ 2]{Cam04}.
	\begin{definition}[{\cite[Definition 2.14]{CP23b}}] \label{conic def}
		Let $\phi$ be a smooth form of $(p,q)$-type with values in $L$ defined on the open set
		$X^\circ.$ We say that $\phi$ is \textit{smooth in conic sense}
		if the quotient of the local inverse images 
		\begin{equation} \label{conic defin expression}
			\widetilde{\phi_i}:=\frac{1}{w_i^{\qqq\mmm}}\pi_i^*(\phi|_{V_i})
		\end{equation}	
		admits a smooth extension to $U_i$. Here in \eqref{conic defin expression} we are using the notation
		$$w_i^{\qqq\mmm}:=(w_i^1)^{\qqq_1\mmm}\cdots (w_i^s)^{\qqq_s\mmm}$$ in order to simplify the writing.
	\end{definition}
	\begin{remark}
		Obviously,	the notion of ``smooth in conic sense'' from the above definition   can be applied to general (not necessarily compact) complex manifolds.
	\end{remark}
	The following proposition plays an important role,  as it builds a correspondence between the intrinsic differential operators and the local ones associated to the data $(\omega_{\mathrm c}, h)$.
	\begin{proposition}[{\cite[Proposition 2.17]{CP23b}}]\label{Prop-dic}
		Let $\phi$ be an $L$-valued $(p,q)$-form, smooth in conic sense. Then its ``natural'' derivatives $D^\prime_{h},D^{\prime *}_{h},\bp,\bp^*,$ are also smooth in conic sense. Furthermore, one has
		
		\begin{enumerate}[{$\mathrm(a)$}]
			\item $\sup_{X^{\circ}}|\phi|_{h,\omega_{\mathrm c}}<\infty$, \text{i.e., forms which are smooth in conic sense are bounded}.
			
			\vspace{0.5em}
			\item \text{The following equalities hold true}
			\begin{equation}\label{dictionary}
				\pi_i^*(D^\prime\phi)=w_i^{\qqq\mmm}D^\prime\widetilde{\phi_i},\quad	\pi_i^*(\bp\phi)=w_i^{\qqq\mmm}\bp\widetilde{\phi_i};
			\end{equation}
			together with
			\begin{equation*}
				\pi_i^*(D^{\prime*}\phi)=w_i^{\qqq\mmm}D^{\prime *}\widetilde{\phi_i},\quad
				\pi_i^*(\bp^*\phi)=w_i^{\qqq\mmm}\bp^*\widetilde{\phi_i}.
			\end{equation*}
			 Here the notation $D^\prime$ on the left-hand side (resp.  the right-hand side) of the first equality in \eqref{dictionary} refers to $D^\prime_{h}$ (resp.  is defined by $D^\prime\varsigma := \partial\varsigma - \partial(\varphi_{L,0}\circ\pi_i) \wedge\varsigma$). We recall that $\varphi_L=\sum_{i=1}^s \qqq_i \log|z^i|^2+\varphi_{L,0}$.
			
			\vspace{0.5em}
			\item Let $\Dt:=[\bp,\bp^*]$ be the Laplace operator with respect to $(\omega_{\mathrm c},h).$ Then one has
			\begin{equation}\label{6655}
				\pi_i^*\Dt\phi=w_i^{\qqq\mmm}\cdot\Delta^{\prime\prime}_\mathrm{sm}\widetilde{\phi_i},
			\end{equation}
			where $\Delta^{\prime\prime}_\mathrm{sm}$ is the Laplace operator for the local, non singular setting $(\pi_i^*\omega_{\mathrm c},\varphi_{L,0}\circ\pi_i)$.
		\end{enumerate}
	\end{proposition}
	
	As a result, we know that forms which are smooth in conic sense behave well by integrations by parts.
	\begin{proposition}[{\cite[Corollary 2.19]{CP23b}}]\label{prop-integration}
		Let $\alpha$ and $\beta$ be an $L$-valued $(p,q)$-form and an $L^*$-valued $(n-p-1,n-q)$-form, respectively, both of which are smooth in  conic sense. Then the usual integration by parts formula holds true
		\begin{equation*}
			\int_X D^\prime_{h}\alpha\wedge\beta={(-1)}^{p+q+1}\int_X \alpha\wedge D^\prime_{h^*}\beta.
		\end{equation*}
	\end{proposition}
	Furthermore, one has the following regularity theorem.
	\begin{proposition}[{\cite[Corollary 2.20]{CP23b}}]
		Suppose that $\zeta$ is  an $L$-valued $L^2$ form on $X$ such that $\Dt(\zeta)=\phi$ holds in the sense of currents (see Definition \ref{def-conic-current})  on $X$ for some $\phi,$ smooth in conic sense. Then $\zeta$ is also smooth in conic sense. In particular, any $L$-valued $L^2$ form which is  $\Dt$-harmonic  is smooth in conic sense.
		
	\end{proposition}

	Setting $\Do:=[D^\prime_{h},D^{\prime*}_{h}]$,\footnote{We leave out the subscript here for simplicity, as long as it does not cause confusion.} one gets the conic version of Bochner--Kodaira--Nakano formula.
	\begin{proposition}[{\cite[Proposition 2.21]{CP23b}}]\label{prop-Bochner-cor}
		Let $\phi$ be an $L$-valued $(p,q)$-form, which is smooth in  conic sense. Then the equality 
		\begin{equation}\label{7788}
			\Dt\phi=\Do\phi+[\theta_L,\Lambda_{\mathrm c}]\phi
		\end{equation}
		holds pointwise on $X^\circ$ (for the notation $\theta_L$ see \eqref{curvature}), where $\Lambda_{\mathrm c}$ is the adjoint of the Lefschetz operator $L_{\mathrm c}:=\omega_{\mathrm c}\wedge\bullet$. Furthermore, we have the following Bochner formula
		\begin{equation}\label{9988}
			\begin{aligned}
				\int_X|\bp\phi|^2 e^{-\varphi_L}dV_{\omega_{\mathrm c}} +\int_X |\bp^*\phi|^2e^{-\varphi_L}dV_{\omega_{\mathrm c}}&=\int_X|D^\prime_{h}\phi|^2e^{-\varphi_L}dV_{\omega_{\mathrm c}}+\int_X |D^{\prime *}_{h}\varphi|^2e^{-\varphi_L}dV_{\omega_{\mathrm c}}\\
				&+\int_X \langle [\theta_L,\Lambda_{\mathrm c}]\phi,\phi\rangle e^{-\varphi_L}dV_{\omega_{\mathrm c}}.
			\end{aligned}
		\end{equation}
	\end{proposition}
	
	Indeed, the special choice of the curvature of the line bundle  $L$  will broaden the range of validity for \eqref{7788}. We will use the following useful proposition later.
	\begin{proposition}\label{pppppp}
		Assume further that  $\qqq_i\geq 0$  for every $i$ 
		in \eqref{weight}. 	Let $\phi$ be an $L$-valued $(p,q)$-form, smooth in  conic sense. Then  the equality 
		$$\Delta^{\prime\prime}\phi=\Delta^{\prime}\phi+[\theta_L,\Lambda_{\mathrm c}]\phi$$
		holds pointwise on the whole of $X$.
	\end{proposition}
	\begin{proof}
		On each local ramified cover \(\pi_i: U_i \rightarrow V_i\),  the relationship \eqref{6655},  combined with the usual Bochner equality  
		\[
		\Delta^{\prime\prime}_\mathrm{sm} = \Delta^{\prime}_\mathrm{sm} + [i\Theta_{\varphi_{L,0}\circ\pi_i}, \Lambda_{\pi_i^*\omega_{\mathrm{c}}}],
		\]
		where \(\Lambda_{\pi_i^*\omega_{\mathrm{c}}}\) is the adjoint of the Lefschetz operator \(L_{\pi_i^*\omega_{\mathrm{c}}} := \pi_i^*\omega_{\mathrm{c}} \wedge \bullet\), implies that, on \(U_i\), we have the following as smooth forms (noting that \(w_i^{qm}\) does not introduce poles along the divisors since \(q_i \geq 0\) for every \(i\)):  
		\[
		\pi_i^*\Dt\phi = w_i^{\qqq\mmm} \cdot \Delta^{\prime\prime}_\mathrm{sm} \widetilde{\phi_i} = w_i^{\qqq\mmm} \cdot (\Delta^{\prime}_\mathrm{sm} \widetilde{\phi_i}+[i\Theta_{\varphi_{L,0}\circ\pi_i}, \Lambda_{\pi_i^*\omega_{\mathrm{c}}}]\widetilde{\phi_i})= \pi_i^* (\Delta^{\prime}\phi+[\theta_L,\Lambda_{\mathrm c}]\phi).
		\]
		This leads to the desired result.
	\end{proof}
	
	For certain bidegree types, we obtain the following straightforward corollary, which has already been established in \cite[Corollary 2.23]{CP23b}. This result will play a role later in establishing the closedness of twisted logarithmic forms of bidegree $(p,0)$, see Theorem \ref{closed}.
	\begin{corollary}\label{particular}
		Assume further that  $\qqq_i\geq 0$  for every $i$ and $\varphi_{L,0}$ is a smooth plurisubharmonic function	in \eqref{weight}.  Let  $\mathfrak{h}$ be  a $\Dt$-harmonic $(n-p,n)$-form  with values in $L$. Then $D^{\prime}_{h}\mathfrak{h}=D^{\prime *}_{h}\mathfrak{h}=0$.
	\end{corollary}
	Similarly to the smooth case, one defines the Hodge operators $*$ and $\sharp$ in our setting, namely, given  an $L$-valued $(p,q)$-form $t$, there exists a unique $L^*$-valued $(n-p,n-q)$-form, denoted by $\sharp t$, such that for any $L$-valued $(p,q)$-form $s$, we have
	\begin{equation}\label{hehe}
		\langle s,t\rangle_{h}dV_{\omega_{\mathrm c}}=s\wedge\sharp t,
	\end{equation}
	where $s\wedge\sharp t$ is calculated  via the natural pairing $L\otimes L^*\rightarrow\mb C$.  One then can derive:
	\begin{proposition}[{\cite[Propositions 2.24 \& 2.25]{CP23b}}]\label{prop-pairing}
		\begin{enumerate}
			\item Let $t$ be an $L$-valued $(p,q)$-form, smooth in conic sense. Then $\sharp t$ is an $L^*$-valued $(n-p,n-q)$-form, also smooth in conic sense (with respect to $(h^*,\omega_{\mathrm c})$).
			\item\label{harmonic good} Let $t$ be a $\Dt$-harmonic form with values in $L$. Then $\sharp t$ is also a 
			$\Dt$-harmonic form with values in $L^*$.
		\end{enumerate}
	\end{proposition}

	Cao--P$\breve{\textrm{a}}$un also obtained the conic version of \textit{G\aa rding}   and \textit{Sobolev inequalities} together with the \textit{Rellich embedding theorem}. As a consequence, they finally got the ``conic Hodge decomposition'' as follows.
	\begin{theorem}[{\cite[Theorem 2.28]{CP23b}}]\label{thm-conic hodge smooth}
		Let $(L,h)$ be a line bundle  on $X$ endowed with a metric $h$ such that the requirements \eqref{weight} and \eqref{restriction2} are satisfied. Let $\omega_{\mathrm c}$ be a K\"ahler metric with conic singularities as in  \eqref{omega_C}. Then we have the following Hodge decomposition 
		\begin{equation}\label{smooth-conic-decom}
			A_\mathrm{co}^{p,q}(X,L)=\mathrm{Ker}\,\Dt_{h}\oplus\mathrm{Im}\,\Dt_{h},
		\end{equation}
		and
		\begin{equation*}
			L^2_{p,q}(X,L)=\mathrm{Ker}\,\Delta^{\prime\prime}_{h}\oplus \mathrm{Im}\,\bp\oplus\mathrm{Im} \,\bp^*,
		\end{equation*}
		where $	A_\mathrm{co}^{p,q}(X,L)$ (resp.  $L^{2}_{p,q}(X,L)$) is the space of $L$-valued $(p,q)$-forms, which are smooth in conic sense (resp.  in the $L^2$ sense) with respect to $(h,\omega_{\mathrm c}).$ 
		
		One can also easily check that $A_\mathrm{co}^{p,q}(X,L)\subset 	L^2_{p,q}(X,L)$.
	\end{theorem}

	\subsection{Proof of Theorem \ref{Thm 1.3}}
	We first give the following simple but important observation:
	\begin{proposition}\label{Lemma:>0=>sc}
		Let $\alpha\in A^{0,q}(X,\Omega_X^p(\log D)\otimes L^* ).$ Suppose that  $\qqq_i>0$ (not necessarily equal to or less than 1) for every $i$ in \eqref{weight}. Then $\alpha$ is  smooth in conic sense.
		\begin{proof}
			This conclusion is based on a local computation. Without loss of generality, we only consider the one-variable case, where $w^{\qqq\mmm}\cdot\frac{d(w^\mmm)}{w^\mmm}$ is smooth. It's worth noting that here $\qqq\mmm$ is a positive integer.
		\end{proof}
	\end{proposition}
	
	\begin{remark}
		Noteworthy to mention that the converse of Proposition \ref{Lemma:>0=>sc} is in general false when the antiholomorphic degree $q\geq 1$. For example, let $X$ be the unit disc with coordinate $z$, $D$ be a simple normal crossing divisor on it defined by the equation $z=0$, and $L$ be the trivial bundle over $X$, endowed with the metric $\varphi_{L}=\frac{1}{2}\log|z|$. In this case, the ($L^*$-valued)-form  $\alpha = \frac{dz \wedge d\overline{z}}{|z|}$  is smooth in  conic sense  but does not admit logarithmic poles along  $D$. However, these two concepts coincide in terms of (Dolbeault type) cohomology, provided that \( 0 < \qqq_i \leq 1 \) for every  $i$  in \eqref{weight}, see Proposition \ref{prop-resolution a}.
		
	\end{remark}
	\begin{notation}\label{notation co}
		\begin{enumerate}
			\item Fix $(p,q)$, denote by $\mathpzc{co}\mathscr{A}^{p,q}_{L^*}$  the  sheaf of germs of $L^*$-valued $(p,q)$-forms that are smooth in conic sense.  The union of elements in  $\cooo\mathscr{A}_{L^*}^{p,q}(X)$, i.e., global sections of $\cooo\mathscr A _{L^*}^{p,q}$,  is then the space $	A_\mathrm{co}^{p,q}(X,L^*)$. \item  Define the \textit{$L^*$-valued $(p,q)$-conic Dolbeault cohomology group} as the following vector space
			$$H^{p,q}_\mathrm{co}(X,L^*):=\frac{\mathrm{Ker}\left\{\bp:{A}_\mathrm{co}^{p,q}(X,L^*)\rightarrow {A}_\mathrm{co}^{p,q+1}(X,L^*) \right\}}{\mathrm{Im} \left\{\bp: A_\mathrm{co}^{p,q-1}(X,L^*)\rightarrow {A}_\mathrm{co}^{p,q}(X,L^*)\right\}}.$$
		\end{enumerate}
	\end{notation}  Now we are going to show:
	\begin{proposition}\label{prop-resolution a}Suppose that  $0<\qqq_i\leq 1$  for every $i$  in \eqref{weight}. Then $(\cooo\mathscr A _{L^*}^{p,\bullet},\bp)$ is an acyclic resolution of $\Omega_X^p(\log D)\otimes L^*$, i.e., one has the following exact sequence of sheaves over $X$:
		\begin{equation}\label{resol}
			\begin{aligned}
				0 &\rightarrow \Omega_X^p(\log D)\otimes L^* \stackrel{\iota}{\rightarrow} \cooo\mathscr{A}_{L^*}^{p,0} \stackrel{\bp}{\rightarrow}  \cooo\mathscr{A}_{L^*}^{p,1} \rightarrow \cdots \\
				&\rightarrow  \cooo\mathscr{A}_{L^*}^{p,q} \stackrel{\bp}{\rightarrow}  \cooo\mathscr{A}_{L^*}^{p,q+1} \rightarrow \cdots \rightarrow  \cooo\mathscr{A}_{L^*}^{p,n} \rightarrow 0,
			\end{aligned}
		\end{equation}
		such that $\mathscr{A}^{p,q}_\mathrm{co}(L^*)$ is an acyclic sheaf on $X$ for any $0\leq p,q \leq n$.
		In particular,  we have the following isomorphism:
		\begin{equation}\label{iso-1}
			H^{q}(X,\Omega_{X}^p(\log D)\otimes L^*)\simeq H^{p,q}_\mathrm{co}(X,L^*).
		\end{equation}
		\begin{proof}
			This is a purely local statement. For any  $x \in X$, which $d$ of these $D_i$ pass, we may choose local holomorphic coordinates \( \{ z^1, \dots, z^n \} \) in a small neighborhood  $V$  around \(x = (0, \dots, 0) \) such that
			$D \cap U = \left\{ z^1 \cdots z^d = 0 \right\}$.
			One then considers the local ramified map
			\begin{equation*}
				\pi:U\rightarrow V,\qquad\pi(w^1,\cdots,w^n):=((w^1)^\mmm,\cdots,(w^d)^\mmm,w^{d+1},\cdots,w^n).
			\end{equation*}
			As all $\qqq_i > 0$, it follows that $\iota$ is indeed injective by Proposition \ref{Lemma:>0=>sc}. Suppose that $\alpha$ is a local $(p,0)$-form with values in $L^*$, smooth in  conic sense. Fix a local basis of $L^*$ on $V$, by Definition \ref{conic def} plus the assumption on $\qqq_i$, write locally
			\begin{equation}\label{alpha}
				\alpha=\alpha_{i_1\cdots i_p}\frac{dz^{i_1}}{(z^{i_1})^{a_{i_1}}}\wedge\cdots\wedge\frac{dz^{i_p}}{(z^{i_p})^{a_{i_p}}},
			\end{equation}
			where $a_{i_j}$ is an integer for any $j=1,\cdots,p$ and $\alpha_{i_1\cdots i_p}$ is a smooth function  that is not divisible by $z^{i_1}, \ldots, z^{i_p}$.  We then get  that 
			
			\begin{itemize}
				\item[--]
				if $i_j\in\left\{1,\ldots,d\right\}$, then 
				$$a_{i_j}\leq 1+\qqq_{i_j}-\frac{1}{\mmm}\leq 2-\frac{1}{\mmm},\,\,\text{so }a_{ij} \text{ is an integer no more than } 1;$$
				\item[--] 
				if $i_j\notin\left\{1,\ldots, d\right\}, $ then
				$$a_{ij}\leq 0.$$
			\end{itemize} 
			Furthermore, if $\bp\alpha=0$, then $\alpha_{i_1\cdots i_p}$ is a holomorphic function. Therefore, $\alpha$ in \eqref{alpha} is a logarithmic form.  Furthermore, as one can verify $\bp\circ\iota=0$ easily, the exactness of \eqref{resol} at level $0$ is thus proved.
			
			\vspace{0.8em}
			Let us now turn to the proof of the exactness of $(\mathpzc{co}\mathscr{A}^{p,q}_{L^*},\bp)$   at any level $q\geq1$. 	Thanks to Proposition \ref{Prop-dic}, we can conclude that $(\mathscr{A}_\mathrm{co}^{p,\bullet}(L^*),\bp)$ is a differential complex. Now, suppose that $\alpha$ is a local $(p,q)$-form with $q\geq 1$, smooth in  conic sense. If $\bp\alpha=0$, then so is $$\widetilde{\alpha}:= w^{\qqq\mmm}\cdot\pi^*\alpha.$$
			Here  
			$w^{\qqq\mmm}=(w^1)^{\qqq_1\mmm}\cdots (w^d)^{\qqq_d\mmm}.$  By definition, $\widetilde{\alpha}$ is a smooth form on $U$ (after smooth extension). Thus, via Dolbeault--Grothendieck lemma, one can find a local smooth form $\beta$ on $U$ such that
			\begin{equation*}
				\bp\beta=\widetilde\alpha.
			\end{equation*}
			Equivalently,
			\begin{equation}\label{pia}
				\frac{\bp\beta}{w^{\qqq\mmm}}=\pi^*\alpha.
			\end{equation}
			Recall that the $\mmm$-ramified cover $\pi:U\rightarrow V$ induces the Galois groups $\rho_1,\cdots\rho_\mmm$ which are the automorphisms $\rho_i:U\rightarrow U$ invariant over $V$. Since the right-hand side of \eqref{pia} is $\rho_i$-invariant, we have
			$$\bp\gamma=\pi^*\alpha,$$
			where $$\gamma:=\frac{1}{\mmm}\sum_{i}\rho_i^*\frac{\beta}{w^{\qqq m}}$$ is $\pi$-invariant and the multiplication $w^{\qqq\mmm}\cdot\gamma$ is smooth. Accordingly, $\gamma$ is the pull back of some form $\varphi$ on $V$, smooth in  conic sense. Therefore, we get, on $U$,
			$$\bp\pi^*\varphi=\pi^*\bp\varphi=\pi^*\alpha.$$
			As a result, on $V$,
			$$\bp\varphi=\alpha.$$
			
			Finally,   $\mathpzc{co}\mathscr A ^{p,q}_{L^*}$   is an acyclic sheaf on $X$ for any $0\leq p,q\leq n$ since one can check that it is a $\mathscr A _X$-module, 
			cf. \cite[Corollary 4.19]{Deme}.
			
			The proof of Proposition \ref{prop-resolution a} is thus completed, by virtue of the standard de Rham--Weil isomorphism theorem, see e.g. \cite[Chapter IV, $\S$ 6, (6.4)]{Deme}.
		\end{proof}
	\end{proposition}
	
	With these preparations, we now come to prove Theorem \ref{Thm 1.3}.
	\begin{proof}[{Proof of Theorem \ref{Thm 1.3}}] 
		Given $\alpha\in A^{0,q}(X, \Omega_{X}^p(\log D)\otimes L^*)$, we can deduce that $\alpha$ is smooth in conic sense according to Proposition \ref{Lemma:>0=>sc}; consequently, so is $\alpha^\prime:=D^\prime_{h^*}\alpha$, as shown in Proposition \ref{Prop-dic}. By assumption $\DDD\alpha$ is $\bp$-closed as a logarithmic form, one can then verify that it is also $\bp$-closed in the conic sense. Using Theorem \ref{thm-conic hodge smooth}, we obtain the following decomposition:
		\begin{equation*}
			\DDD\alpha=\mathfrak{h}+\bp \mathfrak{u},
		\end{equation*} where $\mathfrak{h}\in \mathrm{Ker} \Delta^{\prime\prime}_{h^*}$, $\mathfrak{u}$ is also smooth in conic sense. One also can derive that $\mathfrak{h}=0,$ thus $\DDD\alpha$ is $\bp$-exact. Indeed, denote by $(s,t)$ the inner product of two $L^*$-valued $(p,q)$-forms $s$ and $t$ associated to the pointwise norm with respect to the data $(\omega_{\mathrm c},h^*)$, i.e., $(s,t):=\int_X\langle s,t\rangle_{h^*}dV_{\omega_{\mathrm c}}$ (see \eqref{hehe}). We also  write $||s||=\sqrt{(s,s)}.$ Then, one gets
		\begin{equation*}
			\begin{aligned}
				(\mathfrak{h},\mathfrak{h})&=(\mathfrak{h},\DDD\alpha-\bp \mathfrak u)\\
				&=(\mathfrak{h}, \DDD\alpha)\\
				&=(D^\prime_{h}\sharp\mathfrak{h},\sharp\alpha)=0,
			\end{aligned}
		\end{equation*}
		where the third equality comes from Proposition \ref{prop-integration}, and the last equality holds because $|(D^\prime_{h}\sharp\mathfrak{h},\sharp\alpha)|$ is bounded by \(||D^\prime_{h}\sharp\mathfrak{h}||\cdot ||\sharp\alpha||\), which equals zero by the Bochner formula \eqref{9988} in Proposition \ref{prop-Bochner-cor} (notice that $\sharp h$ is also $\Dt$-harmonic with values in $L$ due to Proposition \ref{prop-pairing} \eqref{harmonic good}).

		Thanks to the  isomorphism \eqref{iso-1},  we can obtain  an $L^*$-valued logarithmic form  $\chi\in A^{0,q-1}(X,\Omega_{X}^{p+1}(\log D)\otimes L^*)$ that satisfies \eqref{solution 1}.	This concludes the proof of Theorem \ref{Thm 1.3}.
	\end{proof}

	\section{General case:  not all  {$\qqq_i>0$}}
	This section aims to prove  Theorem \ref{n,q-main}.  Notice that  in this general case (i.e., not all $\qqq_i>0$ in \eqref{weight}), $\alpha$ may not be smooth in conic sense. For example, 
	let $X$ be the unit bidisc $\mathbb{D}^2=\left\{(z^1,z^2)\in\mb C^2\,:\,|z^i|<1,\,i=1,2\right\}$, $D$ be a simple normal crossing divisor defined by the equation $z^1z^2=0,$ and $L$ be the trivial bundle over $X$, endowed with the metric $\varphi_L=\frac{1}{2}\log|z^1|.$ Then the form $\alpha=\frac{dz^1}{z^1}\wedge\frac{dz^2}{z^2}\in A^{0,0}(X,\Omega_{X}^2(\log D)\otimes L^*)$ defined over $\mathbb{D}^2\setminus{\left\{z^1=0\right\}\cup\left\{z^2=0\right\}}$ is not smooth in conic sense.
	However, it is now indeed a $(L^*,h^*)$-conic current  (see Definition \ref{def-conic-current} and \eqref{mmmm} in Proposition \ref{resolution-2}). Therefore, it should be no surprise that the (conic) current theory will play a crucial role in this section.

	In this section, we \textit{always} let $(L,h)$ be a line bundle  on $X$ endowed with a metric $h$ such that the requirements \eqref{weight} and \eqref{restriction2} are satisfied. Let $\omega_{\mathrm c}$ be a K\"ahler metric with conic singularities as in  \eqref{omega_C}. 
	
In $\S$ 4.1, we introduce conic currents valued in $(L,h)$ and recall the
de Rham--Kodaira decomposition from \cite{CP23b}.

In $\S$ 4.2--4.3, under the assumptions of Theorem~\ref{n,q-main}, we decompose
$D=E+F$ according to whether $\qqq_i=0$ or $\qqq_i>0$. Applying the
de Rham--Kodaira decomposition, the associated conic current
$T_{D^\prime_{h^*}\alpha}$ splits into an $\bp$-exact part and a residue term.
Lemma \ref{key lemma}, inspired by \cite{LRW19,CP23b}, shows that the residue
term is $\bar\partial$-exact modulo on-$E$ conic currents valued in $(L^*,h^*)$.

Finally, we solve equation \ref{solution 2} via an acyclic resolution of
$\Omega_X^p(\log D)\otimes L^*$ (Proposition~\ref{resolution-2}) by  the
sheaf of log-$E$ conic currents valued in $(L^*,h^*)$ 
(Definition \ref{LoGG}), which extends King's result \cite{Kin83} to the
present twisted setting.

	\subsection{de Rham--Kodaira decomposition for conic currents}\label{;;;;}

	\begin{definition}[{\cite[Definition 2.30 \& Remark 2.31]{CP23b}}]\label{def-conic-current}
		A \textit{$(p,q)$-conic current $T $ with values in $(L,h)$} on $X$  is a ``$L$-valued current''  such that there exist a constant $C>0$ and a positive integer $\mathfrak{s}>0$ such that the inequality 
		\begin{equation}\label{conic current expression}
			\left|\int_X T\wedge\phi\right|\leq C\sum _{j=0}^\mathfrak{s}\sup_{X\setminus\mathrm{Supp}\, D}|\nabla^j \phi|_{h^*,\omega_{\mathrm c}}
		\end{equation}
		holds for any  $L^*$-valued $(n-p,n-q)$ form $\phi$ which is moreover smooth in conic sense.   Here $\int_X T\wedge\phi$ is a \textit{formal expression} denoting  the pairing between a conic current $T$ and a test form $\phi$.
	\end{definition}
	The condition \eqref{conic current expression} is equivalent to the following.
	\begin{proposition}[{\cite[Proposition 2.33]{CP23b}}]\label{equivvvvvv}
		A $(p,q)$-conic current \( T \) with values in  $(L, h)$  is represented by a collection of $T^{\mathrm{inv}}_i$ on $U_i$ (where $i\in \mathfrak{I}$) such that for each compact subset $K\subset U_i,$ there exists a constant $C_K>0$ and a positive integer $\mathfrak{s}>0$ such that 
		\begin{equation*}
			\left|\int_{U_i}\frac{T_i^{\mathrm{inv}}}{w_i^{\qqq\mmm}}\wedge \widetilde{\phi_i}\right|\leq C_K \sum_{j=0}^\mathfrak{s}\sup_{U_i}|\nabla^j\widetilde{\phi_i}|
		\end{equation*}
		holds for any $(n-p,n-q)$-form $\phi,$ which is $L^*$-valued, with compact support in $V_i$ and smooth in conic sense. Here $\widetilde{\phi_i}=w_i^{\qqq\mmm}\cdot\pi_i^*(\phi|_{V_i}).$ 
		\begin{proof}
			The equivalence is evident upon observing that locally we have  $T|_{V_i}=\pi_{i*} T_i^{\mathrm{inv}},$ and the relationship
			$$\sum_{j=0}^\mathfrak{s}\sup_{V_i\setminus \mathrm{Supp}\, D}|\nabla^j\phi_i|_{h^*,\omega_{\mathrm c}}\simeq \sum_{j=0}^\mathfrak{s}\sup_{U_i}|\nabla^j\widetilde{\phi_i}|.$$
		\end{proof}
	\end{proposition}

	\vspace{0.4em}
	Thanks to the  ``stabilities''  of forms which are smooth in conic sense under the  ``natural''  operators (cf. Proposition \ref{Prop-dic}), we  can define operators $D_{h}^\prime$, $D_{h}^{\prime *}$, $\bp$, and $\bp^*$ acting on an $(L,h)$-valued $(p,q)$-conic current $T $ by the following way:
	\begin{equation}\label{D prime T def}
		\int_X D^\prime_{h}(T )\wedge \beta:=(-1)^{p+q+1}\int_X T \wedge D^\prime_{h^*}(\beta),\quad
		\int_X D^{\prime *}_{h}(T )\wedge \beta:=(-1)^{p+q}\int_X T \wedge D^{\prime *}_{h^*}(\beta),
	\end{equation}
	and
	\begin{equation}\label{dbar T def}
		\int_X\bp (T )\wedge\beta:=(-1)^{p+q+1}\int_X T \wedge \bp(\beta),\quad
		\int_X\bp^* (T) \wedge \beta:=(-1)^{p+q}\int_X T \wedge\bp^* (\beta).
	\end{equation}
	Via \eqref{D prime T def} and \eqref{dbar T def}, we derive
	\begin{equation}\label{555}
		\int_X \Do (T) \wedge \beta=\int_X T \wedge\Do (\beta),\quad	\int_X \Delta^{\prime\prime}(T )\wedge\beta=\int_X T \wedge\Delta^{\prime\prime}(\beta).
	\end{equation}

	Cao--P$\breve{\textrm{a}}$un showed that the results of de Rham--Kodaira in \cite{deRK50} concerning the Hodge decomposition for currents on compact manifolds have a complete analogue in the conic setting.
	\begin{theorem}[{\cite[p. 23]{CP23b}}]With the same settings for $(X,\omega_{\mathrm c})$ and $(L,h)$ as in the beginning of this section.
		Let $T $ be a conic current of $(p,q)$-type with values in  $(L,h)$. Then, there exists a unique operator, called the \textit{Green operator} $\mathcal{G}$, acting on $T$ by duality and maintaining the degree,  defined as 
		\begin{equation}\label{def-Green current}
			\int_X\mathcal{G}(T)\wedge\beta:=\int_X T\wedge \mathcal{G}(\beta),
		\end{equation}
		where $\beta$ is any  $L^*$-valued form of $(n-p,n-q)$-type  on $X$, which moreover is smooth in conic sense with respect to the data $(\omega_{\mathrm c},h^*)$.
		The Green operator satisfies that
		\begin{equation*}
			\bp\mathcal{G}T=\mathcal{G}\bp T, \,\,\,\,\,\bp^*\mathcal{G}T=\mathcal{G}\bp^*T.
		\end{equation*}
		We furthermore have the following identities:
		\begin{equation*}
			T  - \mathcal{H}T  = \Delta^{\prime\prime}\mathcal{G}T  = \mathcal{G}\Delta^{\prime\prime}T , \quad \mathcal{H}\mathcal{G}T  = \mathcal{G}\mathcal{H}T ,
		\end{equation*}
		where $\mathcal{H}$ is the harmonic projection, defined as
		\begin{equation}\label{Harmonic proj}
			\mathcal{H}(T ) := \sum_{i} \left<T, \zeta_i\right>\cdot\zeta_i,
		\end{equation}
		where $\left\{\zeta_i\right\}_{i}$ is a basis of $L^2$ $\Dt$-harmonic forms of  $(p,q)$-type with values in $L$ and where the notation in \eqref{Harmonic proj} is 
		\begin{equation*}
			\left<T, \zeta_i\right>:=\int_X T\wedge \sharp\zeta_i.
		\end{equation*}
		Here  $\left\{\sharp\zeta_i\right\}_{i}$ is a basis of $L^2$ $\Dt$-harmonic forms of  $(n-p,n-q)$-type with values in $L  ^*$,  see Proposition \ref{prop-pairing} \eqref{harmonic good}.
	\end{theorem}
	
	For the specific restrictions on $(L, h)$ of interest, the following useful properties hold:
	\begin{proposition}
		Suppose further that  $\qqq_i \geq 0$ for every $i$ and that $\varphi_{L,0} = 0$ in \eqref{weight}. Let $T$ be a conic current valued in $(L^*,h^*)$. Then, there exists a Green operator $\mathcal{G}$ acting on $T$ as in \eqref{def-Green current} that satisfies the following relations:
		\begin{equation}\label{commuuuuuuu G}
			\bp\mathcal{G}T = \mathcal{G}\bp T, \quad \bp^*\mathcal{G}T = \mathcal{G}\bp^*T, \quad D^\prime_{h^*}\mathcal{G}T = \mathcal{G}D^\prime_{h^*}T, \quad D^{\prime *}_{h^*}\mathcal{G}T = \mathcal{G}D^{\prime *}_{h^*}T.
		\end{equation}
		Moreover, we have the following identities:
		\begin{equation}\label{dK conic curent iden}
			\Delta^{\prime}\mathcal{G}T = \mathcal{G}\Delta^{\prime}T = T - \mathcal{H}T = \Delta^{\prime\prime}\mathcal{G}T = \mathcal{G}\Delta^{\prime\prime}T, \quad \mathcal{H}\mathcal{G}T = \mathcal{G}\mathcal{H}T,
		\end{equation}
		where the harmonic projection $\mathcal{H}$ is defined as in \eqref{Harmonic proj}.
	\end{proposition}
	\begin{proof}
		This proposition follows from $\Delta^{\prime}T = \Delta^{\prime\prime}T$, which is a consequence of Proposition \ref{pppppp} along with the equalities in \eqref{555}.
	\end{proof}
	
	\subsection{Log conic current}\label{///}
	In this subsection, following  King \cite{Kin83}, we extend the notion of log currents to the conic setting, incorporating the metric on the twisted line bundle.
	
	Note that throughout \textit{this subsection and the next $\S$ \ref{&&&}}, we consistently impose additional assumptions on  $(L,h)$ as specified in the setting of Theorem \ref{n,q-main}.  Write
	\begin{equation*}
		i\Theta_{h}(L)=\sum_{\substack{1\leq i \leq r;\\ \qqq_i=0}}\qqq_i [D_i]+\sum_{\substack{r+1\leq i\leq s;\\ 0<\qqq_i\leq 1}}\qqq_i[D_i],
	\end{equation*}
	and consider the decomposition
	\begin{equation}\tag{\text{Dec.}}\label{D=E+F}
		D=E+F.
	\end{equation} 
	 Here
	\begin{equation*}
		E:= \sum_{1\leq i\leq r}E_i \quad\text{with}\quad E_i=D_i\,\text{ for }1\leq i\leq r;
	\end{equation*}
	and
	\begin{equation*}
		F:= \sum_{r+1\leq i\leq s}F_i \quad\text{with}\quad F_i=D_i\,\text{ for }r+1\leq i\leq s,
	\end{equation*}
	both being simple normal crossing divisors, do not intersect with each other.
	
	\vspace{0.4em}
	Let $Y$ be a simple normal crossing divisor on a complex manifold $X$.
	\begin{definition}[{\cite[Definitions 1.1.3 \& 1.1.7]{Kin83}}]
		\begin{enumerate}
			\item Let 	$\Omega_{X}^p(\mathrm{null}\, Y)$ be the subsheaf of $\Omega^p_X$ consisting of forms that vanish on $Y$. More precisely, if $\iota:Y\rightarrow X$ is the inclusion and $\lambda$ is a holomorphic $p$-form on an open set $V\subset X$, then $\lambda\in \Omega_{X}^p(\mathrm{null}\, Y)(V)$ if and only if $\iota^*\lambda=0$ on $\mathrm{Reg} Y\cap V$. 
			\item 	Set  $$\mathscr A_X^{p,q}(\nulll Y):=\Omega_{X}^p(\nulll Y)\wedge \mathscr{A}_X^{0,q}.$$
			This sheaf is a subsheaf of $\mathscr{A}_X^{p,q}$. The usual wedge product of forms gives $\mathscr A_X^{p,q}(\nulll Y)$ the structure of an $\mathscr{A}_X$-module. In particular, it is an acyclic sheaf on $X$.
		\end{enumerate}
	\end{definition}
	The relationships between these sheaves and the logarithmic sheaf are described  below.
	\begin{proposition}[{\cite[Propositions 1.1.4 \& 1.1.6 \& 1.1.8]{Kin83}}]\label{null sheaf}
		\begin{enumerate}[{$\mathrm(1)$}]
			\item $\Omega_{X}^p(\mathrm{null}\, Y)$ is an  $\mathcal O_X$-coherent sheaf.
			\item \label{null sheaf 2}If $Y$ is
			locally defined  as $z^{1}\cdots z^ k=0$, then \begin{equation}\label{868686}
				z^{1}\cdots z^ k\cdot\Omega^p_X(\log Y)=\Omega_{X}^p(\mathrm{null}\, Y).
			\end{equation} Consequently, 
			$\Omega_{X}^\bullet(\mathrm{null}\, Y)$ is a graded ideal of $\Omega_{X}^\bullet (\log Y).$ Similarly, the bigraded algebra $\mathscr A_X^{\bullet,\bullet}(\nulll Y)$ is a bigraded ideal in $\mathscr A ^{0,\bullet}(\Omega_{X}^\bullet(\log Y)).$
		\end{enumerate}
	\end{proposition}
	\begin{remark}
		In fact, King defined the subsheaf and proved the above propositions for a general divisor, i.e., a complex analytic subset of codimension one. In our setting, however, we focus exclusively on the normal crossing case.
	\end{remark}
	King defined the sheaf of currents which  can annihilate the null-$Y$ forms.
	\begin{definition}[{\cite[Definition 1.3.9]{Kin83}}]
		The \textit{sheaf of on-}$Y$ $(p,q)$-\textit{currents}, $\mathscr{D}^{\prime p,q}_X (\mathrm{on}\, Y),$ is defined as the subsheaf of $\mathscr D_X^{\prime p,q},$ the \textit{sheaf of germs of $(p,q)$-currents on $X$}, that is obtained by imposing the condition that \begin{equation}\label{,,,,,}
			T\in \mathscr{D}^{\prime p,q}_X (\mathrm{on}\, Y, V) \Longleftrightarrow \int_V T\wedge\phi =0
		\end{equation} for all  $\phi \in \Gamma_c(V, \mathscr{A}_X^{n-p,n-q}(\nulll Y)),$ the forms with values in $\mathscr{A}_X^{n-p,n-q}(\nulll Y)$  that have compact support in $V$.
		 Here $\int_V T\wedge\phi$  is still a  formal expression denoting  the pairing between a  current $T$ and a (special) test form $\phi$.
		It follows that  the on-$Y$ currents form a bigraded  subcomplex of $\mathscr{D}_X^{\prime \bullet,\bullet}$.
	\end{definition}
	\begin{remark}Let $T\in\mathscr D_X^{\prime p,q}(V), \mu\in\mathscr{A}_X^{r,s}(V)$, and let $\varphi$ be a test form. Recall that the exterior product of $T$ and $\mu$ is defined by
		$$T\wedge\mu\in \mathscr D_X^{\prime p+r,q+s}(V):\varphi\mapsto \int_V T\wedge (\mu\wedge\varphi).$$
		Another equivalent formulation of  \eqref{,,,,,} is then given by 
		\begin{equation}\label{kkk}
			T\in \mathscr{D}^{\prime p,q}_X (\mathrm{on}\, Y, V) \Longleftrightarrow\,   T\wedge\xi =0\text{ as currents for all } \xi\in\Gamma(V,\mathscr A_X^{n-p,n-q}(\nulll Y)).
		\end{equation}
		\eqref{,,,,,} implies \eqref{kkk} is obvious. For the converse direction, we should notice that 
		$$\int_V T\wedge \phi=\int_V (T\wedge \phi)\wedge\lambda,$$
		where $\lambda$ is any test function equal to 1 on $\mathrm{Supp}\,\phi$.
	\end{remark}
	\begin{remark}
		Every on-$Y$ current has support in $Y$, however, the converse does not generally hold. For example, let $X$ be the unit disc with coordinate  $z$ and let $T$ be  a non-zero on-$Y$ current of  $(1,1)$-type, where $Y$  is non-singular and given by $\left\{z=0\right\}.$ Then  $zT=0$ as currents. But
		$\frac{\p}{\p z}T$ is a current supported in $Y$ such that 
		$$z(\frac{\p}{\p z}T)=T-\frac{\p}{\p z}(zT)=T\not=0,$$ so $\frac{\p}{\p z}T$ is not an on-$Y$ current.
		
		Nonetheless, this property does hold when restricted to rectifiable currents, such as integral currents. This result follows from a flatness theorem, see \cite[Theorem 2.1.8]{Kin71}.
	\end{remark}
	The notation for log currents was also introduced.
	\begin{definition}[{\cite[Definition 1.3.10]{Kin83}}]
		The \textit{sheaf of log-$Y$ currents of bidegree $(p,q)$} is the quotient sheaf
		$$\mathscr{D}^{\prime p,q}_X (\log Y)=\mathscr{D}^{\prime p,q}_X /\mathscr{D}^{\prime p,q}_X (\mathrm{on}\, Y).$$ 
	\end{definition}
	
	\begin{remark}\label{rmk 1}
		\begin{enumerate}[{(1)}]
			\item\label{p=0} When taking $p=0$, one has $\mathscr{D}^{\prime 0,q}_X (\mathrm{on}\, Y)=\left\{0\right\}$, since $$\mathscr A_X^{n,n-q}(\nulll Y)\simeq \mathscr{A}_X^{n,n-q}.$$ As a result, $\mathscr{D}^{\prime 0,q}_X (\log Y)=\mathscr{D}^{\prime 0,q}_X.$ 
			
			\vspace{0.4em}
			\item \label{rmk 1 (1)}Observe that $\mathscr{A}^{0,q} (\Omega_{X}^p(\log Y))$ is canonically isomorphic onto a  subsheaf of the sheaf of log-$Y$ currents of bidegree $(p,q),$ via the quotient map $$\gamma: \mathscr{D}^{\prime p,q}_X\rightarrow \mathscr{D}^{\prime p,q}_X (\log Y).$$ This holds because $$\mathscr{A}^{0,q} (\Omega_{X}^p(\log Y))\cap \mathscr{D}^{\prime p,q}_X (\mathrm{on} \, Y)=\left\{0\right\},$$ given that any  logarithmic form with support contained in $Y$ must be zero.
			
			\vspace{0.2em}
			\item The log current sheaf is flabby because it is the quotient of two flabby sheaves.
			
			\vspace{0.2em}
			\item Sometimes it is convenient to consider the log currents as elements of a function space and not just as a quotient. More precisely, $\mathscr{D}^{\prime p,q}_X (\log Y,  V)$ is canonically isomorphic to the space of linear functionals on $\Gamma_c (V, \mathscr{A}_X^{n-p,n-q}(\nulll Y))$, continuous in the relative topology from $\Gamma_c (V, \mathscr{A}_X^{n-p,n-q}).$ See \cite[Proposition 1.3.12]{Kin83} for a detailed proof.
		\end{enumerate}
	\end{remark}
	\begin{example}
		Let  $f$ be  holomorphic function on $X=\mathbb{C}$ with coordinate $z$.  Set $Y=\left\{0\right\}$. The equivalence  class of $\bp(f\frac{dz}{z})$ in $\mathscr{D}^{\prime 1,0}_{\mathbb{C}} (\log Y),$ denoted by $\left[\bp(f\frac{dz}{z}) \right],$ is equal to $\left[f\delta_0\right]$ and is therefore indeed the zero class. Here $\delta_0$ denotes the Dirac mass at $0$. 
	\end{example}
	The following crucial properties hold, as shown by King, with the quasi-isomorphisms \eqref{KING} playing an important role in our proof of Proposition \ref{resolution-2}.  For the convenience of the readers, a detailed proof of Proposition \ref{PPP} is provided in Appendix \ref{KINGKING}.
	\begin{proposition}[{\cite[Theorems 1.3.11 \& 2.1.2]{Kin83}}]\label{PPP}
		The sheaf $\mathscr{D}^{\prime \bullet,\bullet}_X (\log Y)$ is a bigraded  $\mathscr{A}^{0,\bullet} (\Omega_{X}^\bullet(\log Y))$-module, and this module structure defines a canonical isomorphism $$\Omega_{X}^p(\log Y)\otimes_{\mathcal{O}_X} \mathscr{D}_X^{\prime 0,q}\stackrel{\simeq}{\longrightarrow} \mathscr{D}_X^{\prime p,q}(\log Y).$$ We furthermore have the quasi-isomorphisms of complexes of sheaves:
		\begin{equation}\label{KING}
			\Omega_{X}^p(\log Y)\hookrightarrow\mathscr{A}^{0,\bullet}(\Omega_{X}^p(\log Y))\hookrightarrow \mathscr{D}^{\prime p,\bullet}_X (\log Y).
		\end{equation}
		So,  $( \mathscr{D}^{\prime p,\bullet}_X (\log Y),\bp)$ is an acyclic resolution of $ \Omega_{X}^p(\log Y).$
	\end{proposition}
	
	Motivated by King's work, we will explore the more general twisted case, taking the bundle's metric into account, see construction \eqref{D=E+F}. From now on we assume further that $X$ is a compact K\"ahler manifold.
	\begin{proposition}\label{prop-vanish}
		We have the following isomorphism for the two aforementioned sheaves:
		$$\mathpzc{co}\mathscr{A}^{p,q}_{L}\simeq \mathscr A_X^{p,q}(\nulll F)\otimes L.$$
		In particular, for any \(\beta \in A^{p,q}_{\co}(X, L)\) (see Notation \ref{notation co}), $\beta$ vanishes on each component of F. This property can be seen as the dual version of Proposition \ref{Lemma:>0=>sc}.
	\end{proposition}
	\begin{proof}The proof is nearly identical to that in Proposition \ref{prop-resolution a}. Nevertheless, we provide a proof here for the sake of completeness. 
		Assume there are \( d \) specific \( D_i \) that intersect at a chosen point \( x \in X \). We can select a small neighborhood \( V \) around \( x \) such that \( E \cap V \) is defined by \( \{ z^1 \cdots z^c = 0 \} \) and \( F \cap V \) by \( \{ z^{c+1} \cdots z^d = 0 \} \). In this context, consider the local ramified map given by:
		\[
		\pi: U \rightarrow V, \quad \pi(w^1, \ldots, w^n) = ((w^1)^\mmm, \ldots, (w^d)^\mmm, w^{d+1}, \ldots, w^n).
		\]

		Now, suppose that $\alpha$ is a local $(p,q)$-form with values in $L$, smooth in  conic sense. Fix a local basis of $L$ on $V$, by Definition \ref{conic def} plus the assumption on $\qqq_i$,  write locally
		\begin{equation*}
			\alpha=\alpha_{i_1\cdots i_pj_1\cdots j_q}{dz^{i_1}}\wedge\cdots\wedge{dz^{i_p}}\wedge\dz^{j_1}\wedge\cdots\wedge{\dz^{j_q}},
		\end{equation*}
		where $\alpha_{i_1\cdots i_pj_1\cdots j_q}$ is a smooth function.
		We then obtain 	the following two cases:
		\begin{itemize}
			\item[--] 
			if  each $i_j\in\left\{1,\ldots, c, d+1,\ldots,n\right\}$, then    $\alpha_{i_1\cdots i_pj_1\cdots j_q} \text{ is divisible by } (z^{k}) ^{a_{k}}$, for every $k=c+1,\ldots,d,$ where 
			$$a_{k}\geq {\qqq_{k}},\,\,\text{so }a_{k} \text{ is an integer no less than } 1.$$
			This means that $\alpha_{i_1\cdots i_pj_1\cdots j_q} \text{ is divisible by } z^{c+1}\cdots z^d.$
			\item[--]
			if there exists any $i_j\in\left\{c+1,\ldots,d\right\}$, say without loss of generality,  $i_1,\ldots,i_{l}$, then we set $$\left\{k_1,\ldots,k_{d-c-l}\right\}=\left\{c+1,\ldots,d\right\} \setminus \left\{i_1\ldots,i_l\right\}.$$ We can deduce that
			$\alpha_{i_1\cdots i_pj_1\cdots j_q} \text{ is divisible by } (z^{i_j}) ^{a_{i_j}}$ for $1\leq j\leq l$ and also  by $z^{k_1}\cdots z^{k_{d-c-l}}$ (as shown by the previous case). Here for $1\leq j\leq l$,
			$$ a_{i_j}\geq \qqq_{i_j}+\frac{1}{\mmm}-1,\,\,\text{so }a_{i_j} \text{ is an integer no less than } 0.$$
			Thus, $\alpha_{i_1\cdots i_pj_1\cdots j_q} \text{ is divisible by } z^{k_{1}}\cdots z^{k_{d-c-l}}.$
		\end{itemize} 
		In summary, the above  two cases can imply that the natural sheaf morphism $$z^{c+1}\cdots z^d\cdot\mathscr{A}^{0,q}(\Omega_{X}^p(\log F))\otimes L\longrightarrow\mathpzc{co}\mathscr{A}^{p,q}_{L}$$ is an isomorphism. Thus, the proof is complete,  thanks to Proposition \ref{null sheaf} \ref{null sheaf 2}.
	\end{proof}
	\begin{definition}
		Let $\mathpzc{co}\mathscr{A}_L^{p,q}(\nulll E)$ be the subsheaf of $\mathpzc{co}\mathscr{A}^{p,q}_{L}$ consisting of forms that vanish on $E.$ Then analogous to Proposition \ref{prop-vanish}, we can show that
		\begin{equation}\label{!@!}
			\mathpzc{co}\mathscr{A}^{p,q}_{L}(\nulll E)\simeq  \mathscr A_X^{p,q}(\nulll D)\otimes L.
		\end{equation}
	\end{definition}
	\begin{notation}
		For fixed $(p,q)$, denote by $\codl^{\prime p,q}$ the sheaf of germs of $(L^*,h^*)$-valued $(p,q)$-conic current.  Elements in  $\codl^{\prime p,q} (X)$, i.e., global sections of $\codl^{\prime p,q}$, are then $(L^*,h^*)$-valued  $(p,q)$-conic current on $X$.
	\end{notation}
	\begin{definition}
		The sheaf of on-$E$ $(p,q)$-conic currents valued in $(L^*,h^*)$, denoted by $\codl^{\prime p,q}(\mathrm{on}\,E),$ is defined as a subsheaf of $\codl^{\prime p,q}.$ Specifically, for $T\in \codl^{\prime p,q}(V)$, we say \textit{$T$ is an on-$E$ conic current on $V$}, written as
		$$T \in \codl^{\prime p,q}(\mathrm{on}\,E, V),$$
		if and only if
		$$\int_V T \wedge \phi = 0,$$
		for all $\phi \in \Gamma_c(V, \mathpzc{co}\mathscr{A}_L^{n-p,n-q}(\nulll E)),$
		the forms with values in $\mathpzc{co}\mathscr{A}_L^{n-p,n-q}(\nulll E)$ that have compact support in $V$.
	\end{definition}
	\begin{remark}\label{666}
		Every on-$E$  conic currents  has support in $E$, since the restriction of \(\mathpzc{co}\mathscr{A}_L^{n-p,n-q}(\nulll E)\) to $X - E$ is simply \(\mathpzc{co}\mathscr{A}_L^{n-p,n-q}\).
	\end{remark}
	We proceed to introduce the notation for log-$E$ conic currents valued in $(L^*,h^*)$.
	\begin{definition}\label{LoGG}
		The \textit{sheaf of log-$E$ conic currents of bidegree $(p,q)$ valued in $(L^*,h^*)$} is the quotient sheaf
		\begin{equation}\label{def log current}
			\codl^{\prime p,q}(\log E)=\codl^{\prime p,q}/\codl^{\prime p,q}(\mathrm{on}\, E).
		\end{equation}
	\end{definition}
	\begin{remark}\label{remark differn co}
		\begin{enumerate}[{$\mathrm(1)$}]
			\item  \label{remark differ co (1)}It is straightforward to verify that both $D^\prime_{h}$ and $\bp$ preserve the form valued in $\mathpzc{co}\mathscr{A}_L^{n-p,n-q}(\nulll E)$, and consequently, they also preserve  $\codl^{\prime p,q}(\mathrm{on}\,E).$   Therefore, the operator $\nabla_{h^*}=D^\prime_{h^*}+\bp$  on conic currents makes $\codl^{\prime p,q}(\log E)$ into a double complex of sheaves. By abuse of notation, this operator will also be referred to as $\nabla_{h^*}$ on the quotient. 
			\item \label{<>>>>>}By directly checking the definitions, we find that the two sheaves on the right-hand side of \eqref{def log current} are both flabby via the push-forward operator, which implies that the sheaf of log-$E$ conic currents valued in $(L^*,h^*)$ is also flabby. Consequently, all these sheaves are acyclic on any open set $V$ in $X$. One can then derive that
			$$\codl^{\prime p,q} (\log E, V) \simeq \frac{\codl^{\prime p,q} (V)}{\codl^{\prime p,q}(\textrm{on}\,E, V)}.$$
		\end{enumerate}
	\end{remark}
	
	Next we are going to state the ``singular'' version of Proposition \ref{prop-resolution a}. 
	\begin{proposition}\label{resolution-2}
		One has the following acyclic resolution of $ \Omega_{X}^p(\log D)\otimes L^* $,
		\begin{equation}\label{@2@}
			\begin{aligned}
				0 &\rightarrow \Omega_{X}^p(\log D)\otimes L^* \stackrel{\iota}{\rightarrow} \codl^{\prime p,0} (\log E) \stackrel{\bp}{\rightarrow} \codl^{\prime p,1} (\log E)
				\rightarrow \cdots \\
				&\rightarrow \codl^{\prime p,q} (\log E) \stackrel{\bp}{\rightarrow} \codl^{\prime p,q+1} (\log E) \rightarrow \cdots \rightarrow \codl^{\prime p,n} (\log E) \rightarrow 0.
			\end{aligned}
		\end{equation}
		In particular,  we have the following isomorphism:
		\begin{equation}\label{!!!!!!!!!}
			H^{q}(X,\Omega_{X}^p(\log D)\otimes L^*)\simeq H^{p,q}_\mathrm{co;\mathscr{D}^\prime_{\log E}}(X,L^*):=\frac{\mathrm{Ker}\left\{\bp:\codl^{\prime p,q} (\log E, X)\rightarrow \codl^{\prime p,q+1} (\log E,X)\right\}}{\mathrm{Im} \left\{\bp:\codl^{\prime p,q-1} (\log E,X)\rightarrow\codl^{\prime p,q} (\log E,X)\right\}}.
		\end{equation}
	\end{proposition}	
	\begin{proof}Consider the same setup as in the beginning of the proof of Proposition \ref{prop-vanish}, with the $7$-tuple $(c, d, E, F, \pi, U, V)$.
		
		Assume that $\alpha$ is a local  $L^*$-valued $(p,q)$-form with logarithmic poles along $D$ defined on the open set $V$. Since $\alpha$ has $L^1_{\mathrm{loc}}$-coefficients, it can be associated with a local conic current $T_\alpha$ defined by 
		\begin{equation}\label{mmmm}
			\int_V T_\alpha \wedge \beta := \int_V \alpha \wedge \beta, \quad \beta \in \Gamma_c(V, \mathpzc{co} \mathscr{A}_L^{n-p,n-q}).
		\end{equation}
		The proof of this follows the same argument as in \cite[$\S$ 2.A, \text{Example (2.5)}]{Deme}. This association shows that $\Omega_{X}^p(\log D) \otimes L^*$ can be viewed as a subsheaf of $\codl^{\prime p,0}$.

		Additionally, the map $\iota$ is injective, by the same reasoning as given in Remark \ref{rmk 1} \ref{rmk 1 (1)}: indeed, suppose  \(\alpha \in \codl^{\prime p,q}(\mathrm{on}\, E, V)\). Then, the coefficients of the $L^*$-valued logarithmic form $\alpha$ must vanish everywhere on $X$, owing to the condition that $\mathrm{Supp}(\alpha) \subseteq E$ as a conic current, see Remark \ref{666}.

		Furthermore, we have \(\bp \circ \iota = 0\), by virtue of  Proposition \ref{prop-vanish} and based on the definition of the sheaf \(\codl^{\prime p,q}(\mathrm{on}\, E)\).

		To prove $ \mathrm{Ker}\,\bp\subset \im\,\iota,$ we will need the following property:
		\vspace{1em}
		
		\noindent\textbf{Claim:} On $V$, any local $(p,q)$-type conic current \( T \) valued in $(L^*, h^*)$ can be represented by a usual (untwisted) current $\widetilde{T}_{\mathrm{usu}}$ on $U$, such that $\widetilde{T}_{\mathrm{usu}} = w^{\qqq\mmm} \cdot T_{\mathrm{inv}}$, where  $\pi_* T_{\mathrm{inv}} = T$.
		
		\begin{proof}[Proof of the Claim]
			For any $(n-p, n-q)$-smooth form $\eta$ with compact support $K$ in $U$, since the local current $T_{\mathrm{inv}}$ is $\pi$-invariant, we have
			\begin{equation*}
				\int_U \widetilde{T}_{\mathrm{usu}} \wedge \eta = \frac{1}{\mmm} \sum_{i=1}^{\mmm} \int_U \widetilde{T}_{\mathrm{usu}} \wedge \frac{\rho_{i}^*(w^{\qqq\mmm} \cdot \eta)}{w^{\qqq\mmm}}.
			\end{equation*}
			 Here the Galois group elements $\rho_{i}$ for $i = 1, \ldots, \mmm$ are the automorphisms $\rho_{i} : U \rightarrow U$ that are invariant over $V$ and induced by the $\mmm$-ramified covering $\pi : U \rightarrow V.$
			
			By a similar argument as in Proposition \ref{prop-resolution a}, there exists an $L$-valued form $\phi$ on $V$, smooth in conic sense, such that
			$$
			\sum_{i=1}^\mmm \rho_{i}^*(w^{\qqq\mmm} \cdot \eta) = \pi^* \phi.
			$$
			Invoking Proposition \ref{equivvvvvv} and the compact support condition, there exist two positive constants $C_K^1$ and $C_K^2$ and a positive integer $\mathfrak{s} > 0$ such that 
			\begin{equation*}
				\left|\int_U \widetilde{T}_{\mathrm{usu}} \wedge \eta\right| \leq C_K^1 \sum_{j=0}^{\mathfrak{s}} \sup_{U} |\nabla^j \widetilde{\phi}| \leq C_K^2 \sum_{j=0}^{\mathfrak{s}} \sup_{U} |\nabla^j \eta|,
			\end{equation*}
			where $\widetilde{\phi} = \frac{\pi^* \phi}{w^{\qqq\mmm}}$. This implies that $\widetilde{T}_{\mathrm{usu}}$ is indeed a local $(p,q)$-current. This claim is thus proved.
		\end{proof}
		
		\vspace{0.6em}
		We then let \( [T] \) denote the equivalence class of \( T \in \codl^{\prime p,0}(V) \) in $\codl^{\prime p,0} (\log E, V)$ such that \( \bp[T] = [0] \), i.e., \( \bp T \in \codl^{\prime p,1}(\text{on}\,E,V)\). By the above claim and the quasi-isomorphisms \eqref{KING} given by King, we have
		\[
		\widetilde{T}_{\mathrm{usu}} = w^{\qqq \mmm} \cdot T_{\mathrm{inv}} \in\Gamma(U, \Omega_{X}^p(\log E)).
		\footnote{Here we abuse notation by still denoting by \( E \) (resp.  \( F \)) the inverse image of \( E \) (resp.  \( F \)).}
		\]
		This means that \( T \) is a (local) \( \bar{\partial} \)-closed \( L^* \)-valued \textit{conic logarithmic form}\footnote{This notion naturally generalizes the concept of a form being smooth in  conic sense, as defined in Definition \ref{conic def}.} of bidegree \( (p,0) \), i.e., \( T \) is smooth outside \( D \), and \( \widetilde{T}_{\mathrm{usu}} \) satisfies both of the following:
		\begin{enumerate}[(i)]
			\item it extends smoothly across \( F \); 
			\item it admits at most logarithmic poles along \( E \).
		\end{enumerate}
		Therefore, by similar arguments as in Proposition \ref{prop-resolution a} near \eqref{alpha}, we can deduce that \( T \in \Gamma(V,\Omega_{X}^p(\log D) \otimes L^*) \). Thus, the exactness of \eqref{@2@} at level zero is proved. 
		
		\vspace{0.8em}
		We now proceed to the proof of the exactness of \eqref{@2@} at higher levels. 	By Remark \ref{remark differn co} \ref{remark differ co (1)}, we know that \((\codl^{\prime p, \bullet}(\log E), \bp)\) forms a differential complex.  Now for \( q \geq 1 \), denote by  \( [T] \)  the equivalence class of \( T \in \codl^{\prime p,q}(V) \) in \( \codl^{\prime p,q}(\log E, V) \), such that $\bp[T]=[0]$, i.e., $\bp T \in \codl^{\prime p,q+1}(\text{on}\,E, V).$
		
		By the above claim and King's result again, we deduce that there exists a local untwisted current \( T_1 \in \mathscr{D}_X^{\prime p,q-1}(U) \) such that  
		$$
		w^{\qqq \mmm} \cdot T_{\mathrm{inv}} - \bp T_1 \in \mathscr{D}_X^{\prime p,q}(\text{on}\,E, U).
		$$  
		It is worth noting a classical result by L. Schwartz \cite{Sch55} states that, for any given current \( T \) and holomorphic function \( g \) on an arbitrary complex manifold \( Z \), there always  exists a (not necessarily unique) current \( S \) on $Z$ such that $gS = T$.  
		Using  this, we find that  
		\begin{equation}\label{86}
			T_{\mathrm{inv}} - \bp S \in \frac{\mathscr{D}_X^{\prime p,q}(\text{on}\,E, U)}{w^{\qqq \mmm}},
		\end{equation}
		where \( S \) is a local current on $U$ satisfying $w^{\qqq \mmm} \cdot S = T_1$. Considering the actions of the Galois group elements \( \rho_i \) for \( i = 1, \ldots, \mmm \) on \eqref{86}, we have  
		$$
		T_{\mathrm{inv}} - \bp S_{\mathrm{inv}} \in \frac{\mathscr{D}_X^{\prime p,q}(\text{on}\,E, U)}{w^{\qqq \mmm}},
		$$  
		where 
		$
		S_{\mathrm{inv}} = \frac{1}{\mmm} \sum_{i=1}^{\mmm} {\rho_i}_* S
		$  
		is \(\pi\)-invariant.   It follows that  
		$$
		w^{\qqq \mmm}\cdot\big(T_{\mathrm{inv}} - \bp S_{\mathrm{inv}}\big)=w^{qm}\cdot T_{\mathrm{inv}}-\bp(w^{qm}\cdot S_{\mathrm{inv}}) \in \mathscr{D}_X^{\prime p,q}(\text{on } E, U).
		$$  
		The above claim then implies the existence of a local \((p,q-1)\)-type conic current \( S^\prime \), valued in \( (L^*, h^*) \), such that 
		$
		\pi_* S_{\mathrm{inv}} = S^\prime$  and
		$$
		T - \bp S^\prime \in \codl^{\prime p,q}(\text{on } E, V).
		$$  
		In other words,
		$$[T]=[\bp S^\prime]$$
		in \( \codl^{\prime p,q}(\log E, V) \),
		thus completing the proof.
	\end{proof}
	\begin{remark}
		The construction in \eqref{mmmm} depends on the condition that $\qqq_i \in [0,1] \cap \mathbb{Q}$ in \eqref{weight} (or, more generally, that all $q_i \geq 0$). In the most general case, for any $L$-valued $(p,q)$-form $\lambda$ on $X$ with logarithmic poles along $D$, if these rational numbers are chosen arbitrarily (i.e., $\lambda$ may not have sufficient integrability), there still exists a conic current $T_\lambda$ associated with $\lambda$, as shown in \cite[Examples 2.3.1, second bullet]{CP23b}.
		$$\int_X T _\lambda\wedge\beta:=\lim_{\epsilon\rightarrow 0}\sum_{i} \frac{1}{\delta_i}\int_{U_i}\mu_\epsilon \frac{\theta_i}{w_i^{\qqq\mmm}}\pi_i^*\lambda\wedge\overline{{\widetilde{\beta_i}}}, \quad\beta\in  A _\co^{n-q,n-p}(X,L),$$
		where $\widetilde{\beta_i}=\frac{\pi_i^*(\beta|_{V_i})}{w_i^{\qqq\mmm}},$ $\left\{\theta_i\right\}$ is a partition of unity with respect to the local ramified maps $\left\{\pi_i:U_i\rightarrow V_i\right\}$, the integer $\delta_i$ is the degree of the map $\pi_i$, and $\mu_\epsilon$ is  a family of truncation functions corresponding to the divisor $D$ (see \cite[Lemma 2.1]{CP23b} for the explicit constructions).
		Notice that \eqref{mmmm} is equivalent to the above construction by Cao--P$\breve{\textrm{a}}$un when our restriction is added,  due to the natural isomorphism $L^*\simeq \overline{L}$ and the Lebesgue's dominated convergence theorem.
	\end{remark}
	\begin{corollary}\label{corooooo}
		When $E=\emptyset$, i.e., $D=F$, we obtain the following quasi-isomorphisms of  complexes of sheaves using the above proposition in conjunction with Proposition \ref{prop-resolution a}:
		$$\Omega_{X}^p(\log F)\otimes L^*\hookrightarrow \mathscr{A}^{0,\bullet} (\Omega_X^p(\log F)\otimes L^*)\hookrightarrow \cooo\mathscr{A}_{L^*}^{p,\bullet}\hookrightarrow\codl^{\prime p,\bullet}.$$
		
		This result can be seen as the conic analogue of the well-known quasi-isomorphisms from the Dolbeault--Grothendieck lemma:
		$$\Omega_{X}^p\hookrightarrow\mathscr{A}^{p,\bullet}_X\hookrightarrow\mathscr{D}^{\prime p,\bullet}_X.$$
		
	\end{corollary}
	
	\subsection{Proof of Theorem \ref{n,q-main}}\label{&&&}
	In this subsection, we aim to prove Theorem \ref{n,q-main}. Aside from the de Rham--Kodaira decomposition for conic currents and Proposition \ref{resolution-2} from the previous subsections $\S$\ref{;;;;} and $\S$\ref{///}, the key Lemma \ref{key lemma}, inspired by \cite[Lemma 2.3]{LRW19} and also \cite[Lemma 4.3]{CP23b}, serves as a critical component.
	
	Before proceeding, we need to introduce the concept of  \textit{residue}.
	Given any $\alpha\in A^{0,q}(X,\Omega_X^p(\log D)\otimes L^*)$, locally on  a coordinate subset $V\subset X$ such that $D_1\cap V=\left\{z^1=0\right\}$, one writes
	$$\alpha=\frac{dz^1}{z^1}\wedge\alpha_1+\alpha_2,$$
	where $\alpha_1$ has at most poles along components $D_i, i\not=1$, and $\alpha_2$ is not divisible by the form $\frac{dz^1}{z^1}$. Set
	$$\Res_{D_1}(\alpha)=\alpha_1|_{D_1}$$
	on $V\cap D_1$. It can be checked that $\Res_{D_1}(\alpha)$ is globally well-defined.
	According to our decomposition \eqref{D=E+F}, we set
	$$\Res_D(\alpha):=\Res_E(\alpha)+\Res_F(\alpha),$$
	where $$\Res_E(\alpha):=\sum_{i=1}^{r}\Res_{E_i}(\alpha)\quad\text{ and }\quad\Res_F(\alpha):=\sum_{i=r+1}^{s}\Res_{F_i}(\alpha).$$
	In the sequel, we will always view $\Res_D(\alpha)$ (similarly for $\Res_E(\alpha)$ or $\Res_F(\alpha)$) as a conic current, defined as follows:
	\begin{equation}\label{hmmm}
		\int_X \Res_D(\alpha)\wedge \beta:=\sum_{i=1}^{s}\int_{D_i}\Res_{D_i}\alpha\wedge\iota_{D_i}^*\beta 
	\end{equation}
	for $\beta\in{A}_\co^{n-p,n-q-1}(X,L)$.  Here $\iota_{D_i}: D_i\hookrightarrow X$ represents the natural embedding. For simplicity, we will omit the pull-back operator $\iota_{D_i}^*$ in \eqref{hmmm} whenever it does not lead to ambiguity.
	
	We now present several properties that will be used later.  Notice that \ref{residue formula D} and \ref{residue formula}  in Proposition \ref{properties} can be regarded as the conic version of the usual residue formulae (see for example  \cite[formula (2.2)]{Nog95} and \cite[pp. 10-11]{LRW19}).
	\begin{proposition}\label{properties}
		For any $\alpha\in A^{0,q}(X,\Omega_{X}^p(\log D)\otimes L^*)$, one has
		\begin{enumerate} [{$\mathrm(1)$}]
			\item \label{residue formula D}$D^{\prime}_{h^*}T_\alpha=T_{	D^{\prime}_{h^*}\alpha}\,\,\text{in the sense of conic current};$
			\item \label{residue formula}
			$\bp T_\alpha-T_{\bp\alpha}=2\pi i \,\Res_E(\alpha)\,\,\text{in the sense of conic current};$
			\item\label{commutati}
			$D^{\prime}_{h^*}(\Res_{E_i}\alpha)+\Res_{E_i}(	D^{\prime}_{h^*}\alpha)=0\quad\text{on }E_i;$
			\item \label{commutati db}
			$\bp(\Res_{E_i}\alpha)+\Res_{E_i}(\bp\alpha)=0\quad\text{on }E_i.$
		\end{enumerate}
	\end{proposition}
	\begin{proof}
		We first prove \ref{residue formula} and the proof of \ref{residue formula D} is similar.
		For any $\beta\in{A}_\co^{n-p,n-q-1}(X,L),$
		\begingroup
		\allowdisplaybreaks
		\begin{align*}
			\int_X (\bp T _\alpha-T _{\bp\alpha})\wedge \beta&=\int_X(-1)^{p+q+1}\alpha\wedge\bp\beta-\bp\alpha\wedge\beta\\
			&=-\int_X \bp(\alpha\wedge\beta)\\
			&=-\int_Xd(\alpha\wedge\beta) \\
			&=2\pi i\sum_{i=1}^{s}\int_{D_i} \Res_{D_i}(\alpha)\wedge\beta\\
			&=2\pi i \left(\sum_{i=1}^{r}\int_{E_i} \Res_{E_i}(\alpha)\wedge\beta+\sum_{i=r+1}^{s}\int_{F_i} \Res_{F_i}(\alpha)\wedge\beta\right)\\
			&=2\pi i \sum_{i=1}^{r}\int_{E_i} \Res_{E_i}(\alpha)\wedge\beta.
		\end{align*}
		\endgroup
		 Here the third equality holds due to type considerations, the second-to-last equality follows from the Stokes formula for the regularized integral (cf. \cite[Theorem 2.14]{LZ21} for the one-dimensional case and \cite[Theorem B.4]{CMW23} for arbitrary-dimensional cases), and the last equality arises from Proposition \ref{prop-vanish}.
		
		We then prove \ref{commutati} and the proof of \ref{commutati db} is similar. Suppose that $E_i$ is defined by $z^i=0$ and $h^*=e^\varphi$ on  a coordinate subset $V\subset X$. Write
		$\alpha=\frac{dz^i}{z^i}\wedge\alpha_1+\alpha_2$
		for two logarithmic forms $\alpha_1, \alpha_2$ on $V.$ Then,
		$$D^{\prime}_{h^*}(\Res_{E_i}\alpha)=(\partial \alpha_1+\partial\varphi\wedge\alpha_1)|_{E_i}.$$
		Also,
		\begin{equation*}
			\begin{aligned}
				\Res_{E_i}(D^{\prime}_{h^*}\alpha)&=\Res_{E_i}(D^\prime_{h^*}(\frac{dz^i}{z^i}\wedge\alpha_1))\\&=\Res_{E_i}(-\frac{dz^i}{z^i}\wedge\partial\alpha_1+\partial\varphi\wedge\frac{dz^i}{z^i}\wedge \alpha_1)=-(\partial \alpha_1+\partial\varphi\wedge\alpha_1)|_{E_i}\\
			\end{aligned}
		\end{equation*}
		These yield \ref{commutati}.
	\end{proof}
	As applications, we get
	\begin{corollary}\label{coosjcs}Let $\alpha\in A^{0,q}(X,\Omega_X^{p}(\log D)\otimes L^{*})$ satisfy that $\bp D^\prime_{h^*}\alpha=0$ pointwise on $X^\circ.$ Then as conic currents, 
		\begin{equation}\label{778878}
			\mathcal{H}(T_{D^\prime_{h^*}\alpha})=0,
		\end{equation}
		\begin{equation}\label{9999}
			\bp^*\bp\mathcal G T_{D^\prime_{h^*}\alpha}=2\pi i\sum_{i=1}^r\bp^*\mathcal G\Res_{E_i}(D^\prime_{h^*}\alpha),
		\end{equation}
		and 
		\begin{equation}\label{7777}
			\begin{aligned}
				T_{D^\prime_{h^*}\alpha}=\bp\bp^*\mathcal G T_{D^\prime_{h^*}\alpha}+2\pi i\sum_{i=1}^r\bp^*\mathcal G\Res_{E_i}(D^\prime_{h^*}\alpha).
			\end{aligned}
		\end{equation}
	\end{corollary}
	\begin{proof}
		For \eqref{778878}, 
		\begin{equation*}
			\begin{aligned}
				\mathcal{H}(T_{D^\prime_{h^*}\alpha})&=\sum_i \left(\int_X {D^\prime_{h^*}\alpha}\wedge \sharp\zeta_i\right)\cdot\zeta_i\\&=(-1)^{p+q+1}\sum_i  \left(\int_X {\alpha}\wedge D^\prime_{h}\sharp\zeta_i\right)\cdot\zeta_i=0,\\
			\end{aligned}
		\end{equation*}	
		where the first equality follows from  \eqref{Harmonic proj} (here  $\left\{\zeta_i\right\}_{i}$ is a basis of $L^2$ $\Dt$-harmonic forms valued in $L^*$, which is in particular smooth in conic sense)  and \eqref{mmmm}, the second equality holds by Proposition \ref{properties} \ref{residue formula D} and \eqref{D prime T def}, and the last equality is a consequence of  Proposition \ref{pppppp}.
		
		For \eqref{9999},  it holds due to   the commutativity of $\mathcal{G}$ and $\bp$ (\ref{commuuuuuuu G}),  Proposition \ref{properties} \ref{residue formula} plus the assumption  that $\bp D^\prime_{h^*}\alpha=0$ on $X^\circ$. 
		
		Finally, the de Rham--Kodaira decomposition for conic currents as in \eqref{dK conic curent iden} implies that 
		\begin{equation*}
			\begin{aligned}
				T_{D^\prime_{h^*}\alpha}&=\mathcal{H}T_{D^\prime_{h^*}\alpha}+\bp\bp^*\mathcal G T_{D^\prime_{h^*}\alpha}+\bp^*\bp\mathcal G T_{D^\prime_{h^*}\alpha}.
			\end{aligned}
		\end{equation*}
		Therefore, \eqref{778878} and \eqref{9999} yield \eqref{7777}.
	\end{proof}
	
	\vspace{0.5em}
	We are then particularly interested in the \textit{second term} on the right-hand side of \eqref{7777} and wonder whether it can be expressed in a \textit{more refined shape} in the sense of conic current. 
	
	To gain some insight, we  first examine the \textit{simplest case where $r=1$}, i.e.,  $E=E_1$ is a smooth divisor.  For any $L$-valued $(n-p-1, n-q)$ form $\beta$ on $X$, which is moreover smooth in  conic sense, the pairing of that term and $\beta$ is given by
	\begin{equation}\label{A}
		(**):=2\pi i(-1)^{p+q}\int_E \Res_E (D^\prime_{h^*}\alpha)\wedge \Big(\iota_E^*\bp^*\mathcal{G}\beta\Big).
	\end{equation}
	
	Applying Corollary \ref{coosjcs} to 	$\Res_E (D^\prime_{h^*}\alpha)$, Proposition \ref{properties} \ref{commutati}  and  using the assumption that $E$ has only one component, we obtain
	\begin{equation}\label{wc}
		(**)=2\pi i(-1)^{p+q}\int_E \Big(\bp\bp^*_E\mathcal{G}_E\Res_E (D^\prime_{h^*}\alpha)\Big)\wedge \Big(\iota_E^*\bp^*\mathcal{G}\beta\Big),
	\end{equation}
	where $\mathcal G_E,$ $\bp^*_E$  are the  operators on  $E$ with respect to the induced metrics from $(h,\omega_{\mathrm c})$ on $X$. Then,
	\begin{equation}\label{B}
		\begin{aligned}
			\frac{1}{2\pi i}(**)&=\int_E \Big(\bp^*_E\mathcal{G}_E\Res_E (D^\prime_{h^*}\alpha)\Big)\wedge \Big(\iota_E^*\bp\bp^*\mathcal{G}\beta\Big),\\
			&=\int_E \Big(\bp^*_E\mathcal{G}_E\Res_E (D^\prime_{h^*}\alpha)\Big)\wedge \iota_E^*\beta-\int_E \Big(\bp^*_E\mathcal{G}_E\Res_E (D^\prime_{h^*}\alpha)\Big)\wedge \iota_E^*\mathfrak h\\&\,\,\,\,\,\,\,\,\,\,\,\,\,\,\,\,\quad-\int_E \Big(\bp^*_E\mathcal{G}_E\Res_E (D^\prime_{h^*}\alpha)\Big)\wedge \Big(\iota_E^*\bp^*\mathcal G\bp\beta\Big),\\
		\end{aligned}
	\end{equation}
	where the first equality holds by \eqref{dbar T def},  and the second equality comes from the Hodge decomposition \eqref{smooth-conic-decom} (where $\mathfrak{h}$ is a $\Delta^{\prime\prime}$-harmonic form with values in $L$ on $X$).
	Furthermore, we have
	\begin{equation}\label{C}
		\begin{aligned}
			\int_E 
			\Bigl(
			\bp_E^* \mathcal G_E \Res_E \bigl(D^\prime_{h^*}\alpha\bigr)
			\Bigr)
			\wedge \iota_E^*\mathfrak h
			&= - \int_E
			\Bigl(
			D^\prime_{h^*}\,
			\bp_E^* \mathcal G_E \Res_E \alpha
			\Bigr)
			\wedge \iota_E^*\mathfrak h
			\\
			&= (-1)^{p+q}
			\int_E
			\Bigl(
			\bp_E^* \mathcal G_E \Res_E \alpha
			\Bigr)
			\wedge
			\Bigl(
			\iota_E^* D^\prime_h \mathfrak h
			\Bigr)
			= 0, 
		\end{aligned}
	\end{equation}
	where  the first equality holds because $D^\prime_{h^*}$  commutes with $\Res_E$, $\mathcal{G}_E$ and $\bp_E^*$ (cf. Proposition \ref{properties} \ref{commutati} and \eqref{commuuuuuuu G}), the second equality follows from Proposition \ref{properties} \ref{residue formula D}  plus the commutativity of $\iota_E^*$ and $D^\prime_{h}$, and the last equality is a consequence of Proposition \ref{pppppp}.
	
	Combining  \eqref{A}, \eqref{B} and \eqref{C} together, one gets, as conic currents on $X$,
	\begin{equation}\label{D}
		\bp^*\mathcal G\Res_{E}(D^\prime_{h^*}\alpha)-\bp\bp^*\mathcal{G}{\iota_{E}}_*\bp^*_E\mathcal{G}_E\Res_E (D^\prime_{h^*}\alpha)=\bp^*_E\mathcal{G}_E\Res_E (D^\prime_{h^*}\alpha),
	\end{equation}
	since our  choice of $\beta$ is arbitrary. Here in \eqref{D},  ${\iota_{E}}_*$ is the push-forward operator with respect to the embedding $E\hookrightarrow X$. It is easy to verify that the right-hand side of \eqref{D} is an on-$E$ conic current (i.e., it annihilates all $L$-valued test form that  are smooth in conic sense and moreover vanish on $E$). Let us define:
	$$\widehat{T}:=\bp^*\mathcal{G}\iota_{E*}\bp^*_E\mathcal{G}_E\Res_E (D^\prime_{h^*}\alpha).$$
	Then, we have
	\begin{equation}\label{11111111}
		\bp^*\mathcal G\Res_{E}(D^\prime_{h^*}\alpha)-\bp\widehat{T}\in \codl^{\prime p,q}(\mathrm{on}\, E, X).
	\end{equation}

	\vspace{1.8em}
	In fact, a similar result to \eqref{11111111} also holds when $E$ has several components, i.e., $r > 1$,  as shown in Lemma \ref{key lemma}, which plays a crucial role in the proof of Theorem \ref{n,q-main}. The underlying idea behind the proof remains the same, though it requires a little more involved inductive arguments.

	Before presenting the lemma, we first introduce some notations. For any positive integer $\mathpzc{v} \leq r$, let $\mathrm{I}_\mathpzc v\subset \left\{1,\ldots, r\right\}$ with $|\III_\mathpzc v|=\mathpzc v$. Set $E_{\mathrm{I}_\mathpzc v}:=\bigcap\limits_{i\in \III_\mathpzc v}E_{i}$ and $E_{\mathrm{I}_0}:=X$. We also set $E_{\mathrm{I}_\mathpzc v}:=\emptyset$ for any integer $\mathpzc v$ bigger than $r$. Let $\phi$ be an $L$-valued form, defined on a space containing $E_{\mathrm{I}{\mathpzc v-1}}$, and smooth in  conic sense with respect to the induced metrics from $(h, \omega_{\mathrm c})$.  Define
	$$ \mathpzc{J}_{\mathrm{I}_{\mathpzc v}}(\phi):=\bp^*_{\mathrm{I}_{\mathpzc v-1}} \mathcal G_{\mathrm I_{\mathpzc v-1}}(\phi|_{E_{\mathrm I_{\mathpzc v-1}}}),$$
	where $\mathcal G_{\mathrm{I}_\mathpzc v}$ (resp.  $\bp^*_{\mathrm{I}_\mathpzc v}$)   is the Green operator (resp.   the $\bp^*$-operator) on  $E_{\mathrm{I}_\mathpzc v}$ with respect to the induced metrics from $(h,\omega_{\mathrm c})$ on $X$.  For any fixed sequence $\III_1\subset\cdots\subset \III_{k}$, we define the residue $\res_{E_{\mathrm{I}_k}}$ on $\DDD\alpha$ as follows:
	\begin{equation}\label{Res def}
		\textrm{Res}_{E_{\mathrm{I}_k}}(\DDD\alpha):=\res_{E_{\mathrm{I}_{k}}}\left(\cdots\res_{E_{\mathrm{I}_1}}(\DDD\alpha)\right).
	\end{equation}
	We will also interpret   \eqref{Res def} as a conic current.

	\begin{lemma}\label{key lemma}
		Let $\alpha \in A^{0,q}(X, \Omega_{X}^{p}(\log D) \otimes {L^*})$ satisfy $\bar{\partial} D’_{h^*} \alpha = 0$ pointwise on $X^\circ.$   
		Then for any $k=1,\ldots, r$, there exists an $(L^*,h^*)$-valued conic current $\widehat{T}_k$ on $X$ of $(p+1,q-1)$-type  such that, as conic currents,  
		\begin{equation*}
			T_k:=\sum_{\mathrm I_1\subset\cdots\subset \III_{k}}\bp^*\mathcal{G}{\iota_{E_{\mathrm{I}_1}}}_*\circ\cdots\circ\bp^*_{\III_{k-1}}\mathcal{G}_{\III_{k-1}}	\mathrm{Res}_{E_{\mathrm{I}_k}}(\DDD\alpha)\equiv\bp \widehat{T}_k\quad\mathrm{mod}\,\,\mathpzc{co}\mathscr{D}^{\prime p+1,q}_{L^*}(\mathrm{on}\, E,X),
		\end{equation*}
		where  ${\iota_{E_{\mathrm{I}_k}}}_*$ is the push-forward operator.
		In particular, taking $k=1$, we have
		\begin{equation}\label{important}
			\sum_{i=1}^{r}\bp^*\mathcal G\Res_{E_i}(D^\prime_{h^*}\alpha)\equiv\bp \widehat{T}_1\quad\mathrm{mod}\,\,\mathpzc{co}\mathscr{D}^{\prime p+1,q}_{L^*}(\mathrm{on}\, E,X),
		\end{equation}
		meaning that  on $X$, the difference of two conic currents 
		$$\sum_{i=1}^{r}\bp^*\mathcal G\Res_{E_i}(D^\prime_{h^*}\alpha)-\bp\widehat{T}_1$$
		is an  on-$E$ conic current.
		\begin{proof} For any $L$-valued $(n-p-1, n-q)$ form $\beta$ on $X$, which is moreover smooth in  conic sense, the pairing of $T_k$ and $\beta$ is given by
			\begin{equation*}
				\bm{\mathit{P}_k}:=	\langle T_k,\beta\rangle= m_k\sum_{\mathrm I_1\subset\cdots\subset \III_{k}}\int_{E_{\mathrm{I}_k}}\mathrm{Res}_{E_{\mathrm{I}_k}}(\DDD\alpha) \wedge \left(\mathpzc{J}_{\mathrm{I}_{k}}\circ\cdots\circ\mathpzc{J}_{{\mathrm{I}_1}}\beta\right),
			\end{equation*}
			where $m_k=(-1)^{k(p+q)+\frac{k(k-1)}{2}}$.
			Applying     Corollary \ref{coosjcs} to	$\Res_{E_{\III_{k}} }(D^\prime_{h^*}\alpha)$ and using Proposition \ref{properties} \ref{commutati},  we deduce that 
			\begin{equation}\label{9}
				\begin{aligned}
				\quad	\bm{\mathit{P}_k}=m_k\sum_{\mathrm I_1\subset\cdots\subset \III_{k}}\int_{E_{\mathrm{I}_k}}
					&\Big(\bp^*_{\mathrm{I}_{k}}\bp\mathcal{G}_{\mathrm{I}_{k}}\mathrm{Res}_{E_{\mathrm{I}_k}}(\DDD\alpha) \\&+\bp\,\bp^*_{\mathrm{I}_{k}}\mathcal{G}_{\mathrm{I}_{k}} \mathrm{Res}_{E_{\mathrm{I}_k}}(\DDD\alpha)\Big)\wedge
					\left(\mathpzc{J}_{\mathrm{I}_{k}}\circ\cdots\circ\mathpzc{J}_{\mathrm{I}_1}\beta\right)=m_k(\bm{\mathit{I}}+\bm{\mathit{II}}),\\
				\end{aligned}
			\end{equation}
			\footnote{Compare  \eqref{9} with \eqref{wc}.}
			where 
			\begin{equation*}
				\bm{\mathit{I}}:=\sum_{\mathrm I_1\subset\cdots\subset \III_{k}}
				\int_{E_{\mathrm{I}_k}}\Big(\bp^*_{\mathrm{I}_{k}}\bp\mathcal{G}_{\mathrm{I}_{k}}\mathrm{Res}_{E_{\mathrm{I}_k}}(\DDD\alpha)\Big)\wedge
				\left(\mathpzc{J}_{\mathrm{I}_{k}}\circ\cdots\circ\mathpzc{J}_{\mathrm{I}_1}\beta\right),
			\end{equation*}
			and
			\begin{equation*}
				\bm{\mathit{II}}:=\sum_{\mathrm I_1\subset\cdots\subset \III_{k}}
				\int_{E_{\mathrm{I}_k}}	\Big(\bp\,\bp^*_{\mathrm{I}_{k}}\mathcal{G}_{\mathrm{I}_{k}}\mathrm{Res}_{E_{\mathrm{I}_k}}(\DDD\alpha)\Big)\wedge
				\left(\mathpzc{J}_{\mathrm{I}_{k}}\circ\cdots\circ\mathpzc{J}_{\mathrm{I}_1}\beta\right).
			\end{equation*}
			
			For $\bm{\mathit{I}}$, with the aid of  Proposition \ref{properties} \ref{residue formula}  and  the assumption $\bp D^\prime_{h^*}\alpha=0$  on $X^\circ$,   we arrive at
			\begin{equation*}
				\begin{aligned}
					\bm{\mathit{I}}&=2\pi i\cdot m_k\,\sum_{\mathrm I_1\subset\cdots\subset \III_{k}\subset \III_{k+1}} (-1)^{p+q-k-1}  \int_{E_{\mathrm{I}_{k+1}}}	\mathrm{Res}_{E_{\mathrm{I}_{k+1}}}(\DDD\alpha)\wedge	\left(\bp^*_{\mathrm{I}_{k}}\mathcal{G}_{\mathrm{I}_{k}}\mathpzc{J}_{\mathrm{I}_{k}}\circ\cdots\circ\mathpzc{J}_{\mathrm{I}_1}\beta\right)\\
					&=2\pi i\cdot(-1)^{p+q-k-1}  \bm{\mathit{P}_{k+1}}\\
					&=2\pi i\cdot(-1)^{p+q-k-1} \langle T_{k+1}, \beta\rangle.\\
				\end{aligned}
			\end{equation*}
			
			For $\bm{\mathit{II}}$, by \eqref{dbar T def} one obtains
			\begin{equation*}
				\bm{\mathit{II}}=(-1)^{p+q-k-1}
				\sum_{\mathrm I_1\subset\cdots\subset \III_{k}}
				\int_{E_{\mathrm{I}_k}}	\Big(\bp^*_{\mathrm{I}_{k}}\mathcal{G}_{\mathrm{I}_{k}}\mathrm{Res}_{E_{\mathrm{I}_k}}(\DDD\alpha)\Big) \wedge
				\Big(\bp\,(\mathpzc{J}_{\mathrm{I}_{k}}\circ\cdots\circ\mathpzc{J}_{\mathrm{I}_1}\beta)\Big).
			\end{equation*}
			Then we apply the Hodge decomposition \eqref{smooth-conic-decom} to 
			$(\mathpzc{J}_{\mathrm{I}_{k-1}}\circ\cdots\circ\mathpzc{J}_{\mathrm{I}_1}\beta)|_{E_{\III_{k-1}}}$ to get 
			\begin{equation*}
				\bp(\mathpzc{J}_{\mathrm{I}_{k}}\circ\cdots\circ\mathpzc{J}_{\mathrm{I}_1}\beta)=\mathfrak{h}+(\mathpzc{J}_{\mathrm{I}_{k-1}}\circ\cdots\circ\mathpzc{J}_{\mathrm{I}_1}\beta)|_{E_{\III_{k-1}}}-\mathpzc{J}_{\mathrm I_k}\circ\bp\,(\mathpzc{J}_{\mathrm{I}_{k-1}}\circ\cdots\circ\mathpzc{J}_{\mathrm{I}_1}\beta),
			\end{equation*}
			where $\mathfrak h$ is an $L$-valued $\Dt_{E_{\III_{k-1}}}$-harmonic form on ${E_{\III_{k-1}}}$. By applying the same reasoning as in \eqref{C}, we conclude that
			\begin{equation*}
				\begin{aligned}
					\bm{\mathit{II}}=(-1)^{p+q-k-1}&\sum_{\mathrm I_1\subset\cdots\subset \III_{k}}\int_{E_{\mathrm{I}_k}}
					\Big(\bp^*_{\mathrm{I}_{k}}\mathcal{G}_{\mathrm{I}_{k}}\mathrm{Res}_{E_{\mathrm{I}_k}}(\DDD\alpha)\Big)\\
					&\quad\qquad\wedge\Big((\mathpzc{J}_{\mathrm{I}_{k-1}}\circ\cdots\circ\mathpzc{J}_{\mathrm{I}_1}\beta)-\mathpzc{J}_{\mathrm I_k}\circ\bp\,(\mathpzc{J}_{\mathrm{I}_{k-1}}\circ\cdots\circ\mathpzc{J}_{\mathrm{I}_1}\beta)\Big)=\bm{\mathit{III}}+\bm{\mathit{IV}},
				\end{aligned}	
			\end{equation*}
			where 
			\begin{equation*}
				\bm{\mathit{III}}:=(-1)^{p+q-k-1}\sum_{\mathrm I_1\subset\cdots\subset \III_{k}}
				\int_{E_{\mathrm{I}_k}}\Big(\bp^*_{\mathrm{I}_{k}}\mathcal{G}_{\mathrm{I}_{k}}\mathrm{Res}_{E_{\mathrm{I}_k}}(\DDD\alpha)\Big) \wedge(\mathpzc{J}_{\mathrm{I}_{k-1}}\circ\cdots\circ\mathpzc{J}_{\mathrm{I}_1}\beta),
			\end{equation*}
			and
			\begin{equation*}
				\bm{\mathit{IV}}:=(-1)^{p+q-k}\sum_{\mathrm I_1\subset\cdots\subset \III_{k}}\int_{E_{\mathrm{I}_k}}
				\Big(\bp^*_{\mathrm{I}_{k}}\mathcal{G}_{\mathrm{I}_{k}}\mathrm{Res}_{E_{\mathrm{I}_k}}(\DDD\alpha)\Big)\wedge \Big(\mathpzc{J}_{\mathrm I_k}\circ\bp\,(\mathpzc{J}_{\mathrm{I}_{k-1}}\circ\cdots\circ\mathpzc{J}_{\mathrm{I}_1}\beta)\Big).
			\end{equation*}
			
			For $\bm{\mathit{III}}$, one should notice that the involved term equals to zero as soon as $k\geq 2$. Indeed, there only exists two elements for every fixed pair $\III_{k}\setminus\III_{k-2}$. Then the alternate sum for every $\III_{k-1}$ with $\III_{k-2}\subset\III_{k-1}\subset\III_{k}$ equals to zero by our construction of $\Res_{E_{\mathrm{I}_k}}$ \eqref{Res def}. Therefore,
			\begin{equation*}
				\begin{aligned}
					\bm{\mathit{III}}&=(-1)^{p+q} \sum_{i=1}^{r}\int_{E_i}\Big(\bp^*_{E_i}\mathcal{G}_{E_i}\Res_{E_i}(\DDD\alpha)\Big)\wedge\beta\\
					&=\langle T_{\mathrm{on}\,E}, \beta\rangle
				\end{aligned}
			\end{equation*}
			where
			$$ T_{\mathrm{on}\,E}:=(-1)^{p+q} \sum_{i=1}^{r} \bp^*_{E_i}\mathcal{G}_{E_i}\Res_{E_i}(\DDD\alpha).$$
			This trick is indeed inspired by \cite[(2.3)]{Nog95}.
			We then can easily check that
			$$ T_{\mathrm{on}\,E}\in \mathpzc{co}\mathscr{D}^{\prime p+1,q}_{L^*}(\mathrm{on}\, E,X).$$

			For $\bm{\mathit{IV}}$, repeating the aforementioned arguments in succession, we get
			\begin{equation*}
				\begin{aligned}
					\bm{\mathit{IV}}=(-1)^{p+q-1}\sum_{\mathrm I_1\subset\cdots\subset \III_{k}}\int_{E_{\mathrm{I}_k}}
					\left(\bp^*_{\mathrm{I}_{k}}\mathcal{G}_{\mathrm{I}_{k}}\mathrm{Res}_{E_{\mathrm{I}_k}}(\DDD\alpha)\right)\wedge \left(\mathpzc{J}_{\mathrm{I}_k}\circ\mathpzc{J}_{\mathrm{I}_{k-1}}\circ\cdots\circ\mathpzc{J}_{\mathrm{I}_1}\circ\bp\beta\right).
				\end{aligned}
			\end{equation*}
			Therefore,
			\begin{equation*}
				\begin{aligned}
					\bm{\mathit{IV}}&=(-1)^{\frac{k(k+1)}{2}}\sum_{\mathrm I_1\subset\cdots\subset \III_{k}}\int_X \Big(\bp\bp^*\mathcal{G}{\iota_{E_{\mathrm{I}_1}}}_*\circ\bp^*_{\mathrm{I}_1}\mathcal{G}_{\mathrm{I}_1}{\iota_{E_{\mathrm{I}_2}}}_*\circ\cdots\circ {\iota_{E_{\mathrm{I}_k}}}_*\bp^*_{\mathrm{I}_k}\mathcal{G}_{\mathrm{I}_k}\Res_{E_{\mathrm{I}_k}}(\DDD\alpha)\Big)\wedge \beta\\
					&=\langle \bp T^{\prime}_k,\beta\rangle,
				\end{aligned}
			\end{equation*}
			where 
			$$T^{\prime}_k:=(-1)^{\frac{k(k+1)}{2}} \sum_{\mathrm I_1\subset\cdots\subset \III_{k}} \bp^*\mathcal{G}{\iota_{E_{\mathrm{I}_k}}}_*\circ\bp^*_{\mathrm{I}_1}\mathcal{G}_{\mathrm{I}_1}{\iota_{E_{\mathrm{I}_2}}}_*\circ\cdots\circ {\iota_{E_{\mathrm{I}_k}}}_*\bp^*_{\mathrm{I}_k}\mathcal{G}_{\mathrm{I}_k}\Res_{E_{\mathrm{I}_k}}(\DDD\alpha).$$
			
			Putting everything together, we have, as conic currents on $X$,
			\begin{equation*}
				T_k=2\pi i\cdot(-1)^{p+q-k-1}  {T}_{k+1}+\bp T^{\prime}_k+T_{\mathrm{on}\,E}.
			\end{equation*}
			In other words, 
			\begin{equation*}
				T_k\equiv 2\pi i\cdot(-1)^{p+q-k-1}  {T}_{k+1}+\bp T^{\prime}_k\quad\mathrm{mod}\,\,\mathpzc{co}\mathscr{D}^{\prime p+1,q}_{L^*}(\mathrm{on}\, E,X).
			\end{equation*}
			By induction, the fact  $T_{r+1}=0$ implies
			Lemma \ref{key lemma}.
		\end{proof}
	\end{lemma}

	With these preparations in place, we now turn to the proof of Theorem \ref{n,q-main}.
	\begin{proof}[{Proof of Theorem \ref{n,q-main}}] 
		Thanks to Corollary \ref{coosjcs} and \eqref{important} in Lemma \ref{key lemma}, there exists a conic current $\widehat{T}$ on $X$ such that 
		$$T_{D^\prime_{h^*}\alpha}\equiv\bp\widehat{T}\quad\mathrm{mod}\,\,\mathpzc{co}\mathscr{D}^{\prime p+1,q}_{L^*}(\mathrm{on}\, E,X).$$
		
		Notice that
		$$\codl^{\prime p,q} (\log E, X) \simeq \frac{\codl^{\prime p,q} (X)}{\codl^{\prime p,q}(\textrm{on}\,E, X)}$$
		due to the acyclicity of the sheaf \(\codl^{\prime p,q}(\textrm{on}\,E)\) on $X$, cf. Remark \ref{remark differn co} \ref{<>>>>>}. Consequently, using the isomorphism between the two cohomology groups established in \eqref{!!!!!!!!!} in Proposition \ref{resolution-2}, there exists an $L^*$-valued logarithmic form  $\chi\in A^{0,q-1}(X,\Omega_{X}^{p+1}(\log D)\otimes L^*)$ satisfying \ref{solution 2}.
		Hence, the proof of Theorem \ref{n,q-main} is finalized.
	\end{proof}
	\begin{remark}\label{R}
		\begin{enumerate}[{(a)}]
			\item\label{RR1}As a special case of Theorem \ref{n,q-main}, Theorem \ref{Thm 1.3} can also be obtained via the conic current approach, see Corollary \ref{corooooo} for further details.
			\item\label{RR2}Our method indeed provides an alternative understanding of \cite[Theorem 0.1]{LRW19} where the bundle $L$ is trivial (see \textit{Remark \ref{remark 1.2}}).  Specifically, one can verify that a combination of \cite[Lemma 2.3]{LRW19} with King's original quasi-isomorphism \eqref{KING} will yield the result, without the need to use the bundle-valued Hodge theorem to treat logarithmic  $(p, q)$-forms as bundle-valued  $(0, q)$-forms, as was done in \cite[Proposition 2.4]{LRW19}.
		\end{enumerate}
	\end{remark}
	
	\section{Applications}
In this section, we present several applications of Theorem \ref{n,q-main}. 
We remark that among these, the application concerning the closedness of twisted  logarithmic forms in $\S$ \ref{subsection:clo} relies only on a weaker conclusion obtained in the course of proving Theorem \ref{n,q-main}, while the other applications make full use of Theorem \ref{n,q-main} itself.
	\subsection{Closedness of  twisted logarithmic forms}\label{subsection:clo}
	Our result in this subsection is as follows.
	\begin{theorem}\label{closed}
		Let $X$ be a compact K\"ahler manifold, and let $D = \sum_{i=1}^s D_i$ be a simple normal crossing divisor on $X$. 
		Let $L$ be a holomorphic line bundle over $X$ whose singular metric has a local weight as in \eqref{weight}, such that all fractional numbers $q_i \ge 0$ and $\varphi_{L,0}$ is a smooth plurisubharmonic function. 
		Then, for any $\alpha \in H^0(X, \Omega_X^p(\log D) \otimes L^*)$, we have
		\[
		D'_{h^*} \alpha = 0 \quad \text{pointwise on } X^\circ.
		\]
	\end{theorem}
	\begin{proof}
For any $\alpha\in H^0\!\left(X,\Omega_{X}^p(\log D)\otimes L^*\right)$, we have 
\[
\bp D^\prime_{h^*}\alpha = 0 \quad \text{pointwise on } X^\circ.
\]
Since we are working in a specific bidegree, Corollary \ref{particular} ensures that Corollary \ref{coosjcs} together with \eqref{important} in Lemma \ref{key lemma} remain applicable under the present curvature conditions. Hence we obtain the existence of a conic current $\widehat{T}$ of bidegree $(p,-1)$ (hence trivial) on $X$ such that 
\[
T_{D^\prime_{h^*}\alpha} - \bp \widehat{T}
= T_{D^\prime_{h^*}\alpha} \in \mathpzc{co}\mathscr{D}^{\prime\,p+1,q}_{L^*}(\mathrm{on}\,E,X).
\]
In particular, the conic current $T_{D^\prime_{h^*}\alpha}$ is supported on the simple normal crossing divisor $E$, see Remark \ref{666}. This yields the desired result.
	\end{proof}
	\begin{remark}
	Notice that in Theorem \ref{closed} we may allow $q_i$ to be larger than $1$, since we do not rely on the acyclic resolution of $\Omega_{X}^p(\log D)\otimes L^*$ from Proposition \ref{resolution-2}, which was needed in Theorem \ref{n,q-main}. In that context, the condition $q_i\le 1$ was imposed precisely to exclude the occurrence of higher order poles (for instance, terms such as $dz/z^2$).
	
	\end{remark}
\begin{remark}
	Notice that Theorem \ref{closed} is a special case of a more general result, where $(L, h)$ has analytic singularities and $i\Theta_{h}(L) \ge 0$. This general case was recently proved by Cao--A. H\"oring (see \cite[Proposition 3.1 \& Remark 3.2]{CH24}), using the clever ``integration by parts'' techniques introduced by J.-P. Demailly in \cite{Dem02}, who treated the case $D = \emptyset$ but allowed $L$ to be pseudo-effective with arbitrary singularities. Our approach is different from theirs.
\end{remark}

	\subsection{Degeneracy of spectral sequence}\label{subsection:E1}
	The aim in this subsection is to prove Theorem \ref{Thm 1.4}.
	\begin{proof}
		Thanks to  the logarithmic analogue of the Fr\"olicher spectral sequence
		(cf. \cite[Theorems 1 \& 3]{CFGU97}), one has
		\[
		E_r^{p,q} \simeq \frac{Z_r^{p,q}}{B_r^{p,q}} .
		\]
		
		For $r=1$, the spaces $Z_1^{p,q}$ and $B_1^{p,q}$ are given by
		\[
		Z_1^{p,q}
		=
		\left\{
		\alpha_{p,q} \in A^{0,q}\!\left(X,\Omega_X^p(\log D)\otimes L^*\right)
		\,\middle|\,
		\bp \alpha_{p,q}=0 \text{ on } X^\circ
		\right\},
		\]
		and
		\[
		B_1^{p,q}
		=
		\left\{
		\bp \beta_{p,q-1}
		\,\middle|\,
		\beta_{p,q-1}\in A^{0,q-1}\!\left(X,\Omega_X^p(\log D)\otimes L^*\right)
		\right\}.
		\]
		
		For $r\ge 2$, 
		\begin{equation}\label{Zr-def}
			\begin{aligned}
				Z_r^{p,q}
				=
				\Bigl\{
				\alpha_{p,q}\in &A^{0,q}\!\left(X,\Omega_X^p(\log D)\otimes L^*\right)
				\,\Big|\, \bp\alpha_{p,q}=0 \text{ on } X^\circ,\\
				&\exists\, \alpha_{p+i,q-i}\in
				A^{0,q-i}\!\left(X,\Omega_X^{p+i}(\log D)\otimes L^*\right),
				\ 1\le i\le r-1,\\
				&\text{such that }
				D^\prime_{h^*}\alpha_{p+i-1,q-i+1}
				+ \bp\alpha_{p+i,q-i}=0
				\text{ on } X^\circ
				\Bigr\},
			\end{aligned}
		\end{equation}
		and
		\[
		\begin{aligned}
			B_r^{p,q}
			=
			\Bigl\{
			D^\prime_{h^*}\beta_{p-1,q}
			+ \bp\beta_{p,q-1}
			\,\Big|\, &\exists\, \beta_{p-i,q+i-1}\in
			A^{0,q+i-1}\!\left(X,\Omega_X^{p-i}(\log D)\otimes L^*\right),\\
			&2\le i\le r-1,
			\text{ satisfying}\\
			&D^\prime_{h^*}\beta_{p-i,q+i-1}
			+ \bp\beta_{p-i+1,q+i-2}=0,\\
			&\bp\beta_{p-r+1,q+r-2}=0
			\text{ on } X^\circ
			\Bigr\}.
		\end{aligned}
		\]

	Under these identifications, the differential map of bidegree $(r,1-r)$,
	\begin{equation*}\label{differential-map}
		\zd_r^{p,q} \colon E_r^{p,q} \longrightarrow E_r^{p+r,q-r+1},
	\end{equation*}
	is given by
	\[
	\zd_r^{p,q}\bigl([\alpha_{p,q}]\bigr)
	=
	\bigl[D^\prime_{h^*}\alpha_{p+r-1,q-r+1}\bigr],
	\]
	where $\alpha_{p+r-1,q-r+1}$ is chosen as in \eqref{Zr-def}.
		Furthermore, 
		$$E_{r+1}^{p,q}\simeq \frac{\ker(\zd_r^{p,q}:E_r^{p,q} \rightarrow E_r^{p+r,q-r+1})}{\zd_r^{p-r,q+r-1}(E_r^{p-r,q+r-1})}.$$
		In particular, one has
		$$\zd_1^{p,q}=D^\prime_{h^*}: E^{p,q}_1\simeq H^q(X,\Omega_X^p(\log D)\otimes L^*)\rightarrow E_1^{p+1,q}\simeq  H^q(X,\Omega_X^{p+1}(\log D)\otimes L^*).$$

			By applying Theorem
		\ref{n,q-main} and performing direct computations, we deduce that for every $r\geq 1$ and all $p,q$,
		$$ \zd^{p,q}_r=0.$$ Therefore, in particular,
		\[
		E_1^{p,q}\simeq E_2^{p,q}\simeq \cdots \simeq E_\infty^{p,q}.
		\] 
		We thus get the desired conclusion.
	\end{proof}
	\subsection{Injectivity theorem}\label{subsection:inj}
	In this subsection, we are going to prove Theorem \ref{inj-}. 
	\begin{proof} 
		Define a new holomorphic line bundle \begin{equation}\label{.}
			{L^\prime}	:=\Omega_{X}^n(\log D)\otimes {L}^*\simeq {K_X}\otimes\mathcal{O}_X(D)\otimes L^*,
		\end{equation} 
		which satisfies 
		$$L^\prime=\sum_{i=1}^{s}\mathcal{O}_X(q_iD_i)\in\mathrm{Pic}_{\mathbb{Q}}(X),$$ where $q_i:=1-b_i\in[0,1)$ for all $i$.
		
		\begin{step}
			Injectivity of $j$ \eqref{VVV}.
		\end{step}By \eqref{.}, it is equivalent to show that the restriction homomorphism
		\begin{equation}\label{res}
			H^q(X,K_X\otimes\mathcal{O}_X(D)\otimes {L^\prime}^*)\rightarrow H^q(X^\circ,K_{X^\circ}\otimes  {L^\prime}^*|_{X^\circ})
		\end{equation}
		is injective, for all $q$.
		
		Let  $\iota: X^\circ\hookrightarrow X$ be  the inclusion. Notice that there is an injective map of complexes
		$$\jmath: \Omega_{X}^\bullet(\log D)\otimes {L^\prime}^*\hookrightarrow 
		\iota_*\Omega_{X^\circ}^\bullet\otimes {L^\prime}^*,$$
		which is a quasi-isomorphism by Deligne's result  (cf. \cite[II, Corollary 3.14]{Del70} and \cite[$\S$ 2, 2.11]{EV92} in the analytic case) since all eigenvalues of $\textrm{Res}_{D_i}(\nabla_{h_{{L^\prime}^*}})=q_i$ along $D_i$ lie in $[0,1)$ according to the formula \eqref{`}.
		
		We then consider the  naive filtrations $F$ on these two complexes.   The inclusion $F^n\subseteq F^0$ induces the following commutative digram:
		\[
		\begin{tikzcd}[row sep=4.6em, column sep=4.6em]
			\mathbb{H}^{q+n}(X,F^n(\Omega_{X}^\bullet(\log D)\otimes {L^\prime}^*)) \ar[r, "\alpha"] \ar[d, "\beta^n"'] & \mathbb{H}^{q+n}(X,\Omega_{X}^\bullet(\log D)\otimes {L^\prime}^*) \ar[d, "\beta"] \\
			\mathbb{H}^{q+n}(X,F^n(\iota_*\Omega_{X^\circ}^\bullet\otimes {L^\prime}^*)) \ar[r] & \mathbb{H}^{q+n}(X,\iota_*\Omega_{X^\circ}^\bullet\otimes {L^\prime}^*).
		\end{tikzcd}
		\]
		Since $\jmath$ is a quasi-isomorphism,  it follows that $\beta$ is an isomorphism. On the other hand, Theorem \ref{Thm 1.4} implies that $\alpha$ is injective. We conclude from the commutative diagram that $\beta^n$ is also injective.
		
		\vspace{0.3em}
		But $$F^n(\Omega_{X}^\bullet(\log D)\otimes {L^\prime}^*)=\Omega_{X}^n(\log D)\otimes {L^\prime}^*\,[-n]$$ and $$F^n(\iota_*\Omega_{X^\circ}^\bullet\otimes {L^\prime}^*)=\iota_*\Omega_{X^\circ}^n\otimes {L^\prime}^*\,[-n].$$ Therefore, $\beta^n$ becomes
		$$\beta^n: H^q(X,\Omega_{X}^n(\log D)\otimes {L^\prime}^*)\rightarrow H^q(X,\iota_*\Omega_{X^\circ}^n\otimes {L^\prime}^*).$$
		Moreover, the inclusion $\iota$ is a Stein morphism\footnote{To be more precise we shall call a holomorphic map $p:\Omega\rightarrow X$ a \textit{Stein morphism} of complex manifolds if every point $x\in X$ has a neighborhood $V=V(x)$ such that $p^{-1}(V)$ is Stein. Indeed, if $R$ is an effective Cartier divisor on a complex manifold $X$, then the inclusion $X\setminus\textrm{Supp}\,R\hookrightarrow X$ is a Stein morphism. This is because the property is local on $X$, which allows us to assume that $R$ is defined by an equation in $\mathcal{O}(X)$, where $X$ is Stein. In this case, the assertion is clear, see e.g. \cite[III, Examples 3.3 (6)]{GPR94}.}, so the higher direct images $R^q \iota_*(\Omega_{X^\circ}^n\otimes \iota^*{L^\prime}^*)$ vanish for $q\geq 1$ by Cartan's Theorem B, and hence 
		$$H^q(X,\iota_*\Omega_{X^\circ}^n\otimes {L^\prime}^*)\simeq H^q(X^\circ, \Omega_{X^\circ}^n\otimes {L^\prime}^*|_{X^\circ})$$
		by virtue of the Leray spectral sequence plus the projection formula. As a result, 
		$\beta^n$ becomes the restriction map \eqref{res}	as desired. 
		\begin{step}
			Injectivity of $j^\prime$ \eqref{XXX} and its equivalence to that of  $j$. 
		\end{step}
		We note that the corresponding result in the algebraic setting was first proposed in \cite[Remark 2.6]{Amb14}.  Here we adapt the arguments therein to the analytic case.

		We have 
		$$\mathcal{O}_X\hookrightarrow\mathcal{O}_X(\widehat D)\hookrightarrow\iota_*\mathcal{O}_{X^\circ},$$
		where the first map is given by multiplication with the section defining $\widehat D$. By tensoring this with $L$ and taking the $q$-th cohomology, we obtain 
		\[
		\begin{tikzcd}
			H^q(X,L) \arrow[r, "j'"] \arrow[rrr, "j", bend right] & H^q(X,L\otimes\mathcal{O}_X(\widehat D)) \arrow[r] & H^q(X, \iota_*\mathcal{O}_{X^\circ}\otimes L) \arrow[r, "\simeq"] & H^q(X^\circ, L|_{X^\circ}),
		\end{tikzcd}
		\]
		where the isomorphism follows from the fact that $\iota$ is a Stein morphism. So the injectivity of $j$ implies that of $j^\prime$.

		Conversely, suppose that
		$ H^q(X, L) \to H^q(X, L \otimes \mathcal{O}_X(\widehat{D})) $
		is injective for all divisors \( \widehat{D} \) supported on \( X^\circ \). Then, it follows that the map  
		$ H^q(X, L) \to H^q(X, L \otimes \mathcal{O}_X(mD)) $
		is also injective for every \( m \geq 0 \).
		By Lemma \ref{direc}, this implies the injectivity of the map  
		$$ H^q(X, L) \to \varinjlim_m H^q(X, L \otimes \mathcal{O}_X(mD)). $$  
		
		Now, observe that  
		$$ \varinjlim_m H^q(X, L\otimes \mathcal{O}_X(mD)) \simeq H^q(X, \varinjlim_m\,(L \otimes \mathcal{O}_X(mD))) \simeq H^q(X, \iota_*(L|_{X^\circ})) \simeq H^q(X^\circ, L|_{X^\circ}), $$  
		where the first isomorphism holds by the commutativity of directed limits with cohomology on quasi-compact topological spaces (see e.g. \cite[Tag 01FE, Lemma 20.19.1]{Sta}),
		the  second isomorphism follows from the properties of sheaf operations,
		and the last isomorphism relies on the Stein property of the inclusion map \( \iota \) again.
		We then achieve the injectivity of the map $j$.
		
		The proof of Theorem \ref{inj-} is, therefore, complete.
	\end{proof}
		\begin{remark}\label{difference}
			Our proof differs from Fujino's (\cite[Theorem 1.2]{Fuj17}) in the following way: 
			Fujino's argument relies on constructing a cyclic cover (up to a normalization) and then, 
		proves the $E_1$-degeneration 
			of the Hodge--de Rham spectral sequence associated to the mixed Hodge structure 
			for cohomology with compact support on the covering space. In contrast, our approach 
			is based on an $E_1$-degeneration result that holds directly on the original manifold 
			$X$ (Theorem~\ref{Thm 1.4}), thus bypassing the covering construction and the subsequent 
			Galois-descent arguments.
		\end{remark}
	\begin{remark}
		Theorem \ref{inj-} also applies to the case where $L$ is a holomorphic line bundle over $X$ satisfying $L\sim_{\mb R} K_X+\sum_{i=1}^s b_i D_i$ where $b_i\in (0,1]$ for any $i$ by a standard approximation argument, cf. \cite[Remark 2.9]{Fuj17} for more details.
	\end{remark}
	\subsection{Unobstructed deformations}\label{subsection:un}
	In this subsection, we are devoted to proving Theorem \ref{main thm}.
	\subsubsection{Deformations of pairs}\label{subsection 5.4.1}
	We begin by reviewing basic notations and properties of locally trivial deformations of pairs, following \cite{Kaw78}. For further details, see \cite[$\S$ 4.4 \& Exercise 4.6.7]{Man22} and \cite[$\S$ 3.4.4]{Ser06}.
	\begin{definition}[{\cite[Definition 3]{Kaw78}}]
		A family of \textit{locally trivial (infinitesimal)  deformations}\footnote{Kawamata referred to these as  ``logarithmic deformations''.} of   a pair $(X,D)$, or a family of locally trivial (infinitesimal) deformations of the closed embedding $D\hookrightarrow X,$ is a $7$-tuple $\mathscr F=(\mathcal X^\circ, {\mathcal{X}}, {\mathcal{D}},{\pi}, S, s_0,{\psi})$ that satisfies the following conditions:
		\begin{enumerate}[{a)}]
			\item ${\pi}:\mathcal{X}\rightarrow S$ is a proper smooth morphism between these two complex spaces. 
			\item   ${\mathcal{D}}$ is a closed analytic subset of $\mathcal{X}$ and $\mathcal{X}^\circ=\mathcal{X}-\mathcal{D}.$
			\item $\psi: X\rightarrow\pi^{-1}(s_0)$ is an isomorphism such that $\psi(X-D)=\pi^{-1}(s_0)\cap\mathcal{X}^\circ$.
			\item $\pi$  locally is  a  projection of a product space as well as the restriction of it to $\mathcal{D}$, that is, for each $p\in\mathcal{X},$ there exist an open neighborhood of $\mathpzc U$ of $p$ and an isomorphism  $\zeta: \mathpzc U\stackrel{\simeq}{\longrightarrow} \mathpzc V\times \mathpzc W$ , where $\mathpzc V=\pi (\mathpzc U)$ and $\mathpzc W=\mathpzc U\cap\pi^{-1}(\pi(p))$, such that the following diagram
			\[
			\begin{tikzcd}[row sep=4em, column sep=4em]
				\mathpzc U \arrow[r, "\zeta"] \arrow[dr, "{\pi}"'] & \mathpzc V \times \mathpzc W \arrow[d, "\textrm{pr}_1"] \\
				& \mathpzc V
			\end{tikzcd}
			\]
			is commutative and the restriction $\zeta|_{\mathpzc U\cap \mathcal D} :\mathpzc U\cap \mathcal D\stackrel{\simeq}{\longrightarrow}\mathpzc V\times (\mathpzc W \cap \mathcal{D})$ is also an isomorphism.
		\end{enumerate}
	\end{definition}

	\begin{definition}\label{tangent}
		Let $T_X^1(-\log D)$ be the \textit{logarithmic tangent bundle}, which is the dual bundle of $\Omega_{X}^1(\log D)$. Recall that
		$T_X^1(-\log D)$ is naturally a subsheaf of the \textit{holomorphic tangent bundle} ${T}_X^{1,0}$, and that it is closed
		under the Lie bracket on ${T}_X^{1,0}$. Indeed, under the same settings near \eqref{hhhhhh}, $T_X^1(-\log D)$ is locally
		generated by the $n$ commuting vector fields
		\[
		z^1 \frac{\partial}{\partial z^1}, \ldots, z^d \frac{\partial}{\partial z^d}, \frac{\partial}{\partial z^{d+1}}, \ldots, \frac{\partial}{\partial z^n},
		\]
		and is therefore closed under the Lie bracket.
	\end{definition}
	
	As is well-known (cf. \cite{Gro58}), the set of the  \textit{first-order locally trivial  deformations of a pair} $(X, D)$ (i.e., families of locally trivial deformations over the space $\operatorname{Spec} \mathbb{C}[\epsilon]/(\epsilon^2)$) is the finite-dimensional space $H^1(X, T_X^1(-\log D))$. Furthermore, the  \textit{obstructions} for  locally trivial deformations lie in $H^2(X, T_X^1(-\log D))$. In the usual way, one has a \textit{Kodaira--Spencer map}
	\[
\rho_{s_0}: T_{S, s_0} \rightarrow H^1(X, T_X^1(-\log D)).
	\]
	
	\begin{definition}[{\cite[Definition 5]{Kaw78}}]\label{Kuranishi-}
		A family $\mathscr K=(\mathcal X^\circ, {\mathcal{X}}, {\mathcal{D}},{\pi}, S, s_0,{\psi})$  of  locally trivial deformations of a pair $(X,D)$ is called to be \textit{semi-universal} if for any family $\mathscr K^\prime=({\mathcal X^\circ}^\prime, {\mathcal{X}}^\prime, {\mathcal{D}}^\prime,{\pi}^\prime, S^\prime, {s_0}^\prime,{\psi}^\prime)$ of locally trivial  deformations of the pair $(X,D),$  there exist  an open neighborhood $S^{\prime\prime}$ of ${s_0}^{\prime}$ in $S^\prime$ and a morphism $\alpha: S^{\prime\prime}\rightarrow S$ such that the following conditions are satisfied:
		\begin{enumerate}[{1.}]
			\item the restriction $\mathscr{K}^\prime|_{S^{\prime\prime}}$ of $\mathscr{K}$ over $S^{\prime\prime}$ is isomorphic to the induced family $\alpha^*\mathscr{K}$;
			\item the differential map $T_{S^{\prime\prime}, {s_0}^\prime}\rightarrow T_{S,s_0}$ is unique.
		\end{enumerate}
	\end{definition}
	
	Kawamata then proved the following Kuranishi-type theorem.

	\begin{theorem}[{\cite[Theorem 1]{Kaw78}}]
		There exists a semi-universal (or Kuranishi) family $\mathscr{K}$ of locally trivial  deformations of the pair $(X, D)$.
	\end{theorem}
	Indeed,  as shown in \cite[p. 251]{Kaw78}, the Kuranishi family 
	$\mathscr{K}$ can be obtained as a subset of
	\[
	\Gamma_{\text{real analytic}}(X, T_X^1(-\log D)\otimes \Lambda^{0,1}T^*X)
	\subset A^{0,1}(X,T_X^1(-\log D)),
	\]
	the so-called \emph{space of logarithmic Beltrami differentials}. 
	This space consists of sections $\varphi$ satisfying the integrability condition
	\[
	\bar{\partial}\varphi = \frac{1}{2}[\varphi,\varphi].
	\]
	
	\begin{definition}
		We say that a pair $(X, D)$ has \textit{unobstructed locally trivial  deformations} if its Kuranishi space $B$, 
		that is, the base of the Kuranishi family $\mathscr K$, is smooth. 
		In this case, $B$ can be identified with a sufficiently small open neighborhood of the origin in the 
		$\mathbb{C}$-vector space $H^1(X,T_X^1(-\log D))$.
	\end{definition}
	
	\subsubsection{Differentiable graded Batalin--Vilkovisky algebra}We next introduce the definition and key properties of the differentiable graded Batalin--Vilkovisky algebra, following primarily \cite{Man22} and \cite{KKP08}, to lay the groundwork for establishing the unobstructedness in $\S$ \ref{><}.
	\begin{definition}[{\cite[p. 145]{Man22}}]
		A \textit{differential-graded Lie algebra} (\dgla) is the data of a differential-graded vector space $(\mathrm{\mathbf L},\mathrm{\mathbf d})$ together with a bilinear map $\left[\cdot,\cdot\right]:\LLL\times\LLL\rightarrow\LLL$ (called \textit{bracket}) of degree $0$, such that the following conditions hold:
		\begin{enumerate}[{(a)}]
			\item  (graded skew-symmetry)\,\, $[a,b] = -(-1)^{|a||b|} [b,a]$;
			\item (graded Jacobi identity)\,\, $[a,[b,c]] = [[a,b],c] + (-1)^{|a||b|} [b,[a,c]]$;
			\item (graded Leibniz rule)\,\, $\mathrm{\mathbf{d}}[a,b] = [\mathrm{\mathbf{d}} a, b] + (-1)^{|a|} [a, \mathrm{\mathbf{d}} b]$.
		\end{enumerate}
		 Here $|a|$ denotes the \textit{degree} of $a$.
	\end{definition}
	Particularly, the Leibniz rule implies that the bracket of a \dgla\,\,induces a structure of graded Lie algebra on its cohomology. Furthermore, we have the following terminologies.
	
	\begin{itemize}
		\item[--]A \dgla\,\,is said to be \textit{abelian} if its bracket is trivial.
		\item [--]A \textit{morphism} between two \dgla s \,$\zeta:\LLL\rightarrow\mathrm{\mathbf M}$ is a linear map that commutes with both the brackets and differentials, and also preserves degrees. 
		\item[--]A \textit{quasi-isomorphism} between two DGLAs\, $\zeta:\LLL\rightarrow\mathrm{\mathbf M}$ is a morphism that induces an isomorphism in cohomologies.
		\item[--]Two DGLAs $\LLL$ and $\mathrm{\mathbf M}$ are said to be  \textit{homotopy equivalent}, if there exists a quasi-isomorphism between them. In particular, a \dgla\, is said to be \textit{homotopy abelian} if it is homotopy equivalent to an abelian \dgla.
	\end{itemize}

	\begin{definition}[{\cite[p. 440]{Man22}}]\label{dgbbbbbb}
		Let \( k \) be a fixed odd integer. A \textit{differential graded Batalin--Vilkovisky algebra (\dgbva) of degree \( k \)} over \( \mathbb{C} \) is the data \( (\AAA, \ddd, \bfitDelta) \), where \( (\AAA, \ddd) \) is a differential $\mb Z$-graded commutative  algebra with unit $1\in\AAA$, and \( \bfitDelta\in\mathrm{Hom}_{\mathbb C}^{-k}(\AAA,\AAA)\) is an order $\leq 2$ differential operator of degree \(-k \), such that   the following properties hold:
		\begin{enumerate}
			\item \(\ddd(1)=\bfitDelta(1)= 0 \);
			\item \( \ddd^2 = \bfitDelta^2 = \ddd \bfitDelta + \bfitDelta \ddd = 0 \).
		\end{enumerate}
	\end{definition}
	\begin{remark}
		For the case  $k = 1$, one may also refer to \cite[Chapter III.9]{Man99} or \cite[Definition 4.12]{KKP08} for the definition.
	\end{remark}
	\begin{remark}
		By  e.g.  \cite[Corollary 9.3.5]{Man22}, the  conditions $\bfitDelta\in\mathrm{Diff}^2_{\AAA/\mathbb C}(\AAA,\AAA)$ and $\bfitDelta(1)=0$ in Definition \ref{dgbbbbbb} are equivalent to the \textit{seven-term relation}:
		\[
		\bfitDelta(abc) + \bfitDelta(a)(bc) + (-1)^{|a||b|} \bfitDelta(b)ac + (-1)^{|c|(|a|+|b|)} \bfitDelta(c)ab 
		\]
		\[
		= \bfitDelta(ab)c + (-1)^{|a|(|b|+|c|)} \bfitDelta(bc)a + (-1)^{|b||c|}\bfitDelta(ac)b.
		\]
		
		 Here  for every integer \(l\), we denote by
		$$
		\mathrm{Diff}^l_{\AAA/\mathbb C}(\AAA, \AAA) \subseteq \mathrm{Hom}^\bullet_{\mb C}(\AAA, \AAA)
		$$
		the \textit{graded subspace of differential operators of  order \(\leq l\)}; it is defined recursively by setting 
		$$
		\mathrm{Diff}^l_{\AAA/\mathbb C}(\AAA,\AAA) = 0$$ for  $l<0$,
		and
		$$
		\mathrm{Diff}^l_{\AAA/\mathbb C}(\AAA,\AAA) = \{ f \in \mathrm{Hom}^\bullet_\mb C(\AAA, \AAA) \mid [f, a] \in \mathrm{Diff}^{l-1}_{\AAA/\mathbb C}(\AAA,\AAA), \,\, \text{for all } a \in \AAA \}
		$$
		for \(l \geq 0\).  For \(f \in \mathrm{Diff}^{l_1}_{\AAA/\mathbb C}(\AAA,\AAA)\) and \(g \in \mathrm{Diff}^{l_2}_{\AAA/\mathbb C}(\AAA,\AAA)\) we have (see \cite[$\S$ 2.3]{Man22} for  similar computations):
		$$
		fg \in \mathrm{Diff}^{l_1+l_2}_{\AAA/\mathbb C}(\AAA,\AAA) \quad \text{and} \quad [f, g] \in \mathrm{Diff}^{l_1+l_2-1}_{\AAA/\mathbb C}(\AAA,\AAA).
		$$
		In particular, 
		$$
		\mathrm{Diff}^1_{\AAA/\mathbb C}(\AAA,\AAA) \quad \text{and} \quad \mathrm{Diff}_{\AAA/\mathbb C}(\AAA,\AAA):= \bigcup_l \mathrm{Diff}^l_{\AAA/\mathbb C}(\AAA,\AAA)
		$$
		are graded Lie sub-algebras of \(\mathrm{Hom}^\bullet_\mb C(\AAA,\AAA)\).
	\end{remark}
	\begin{remark}
		As is well-known (see e.g. \cite{Kos85}, \cite[$\S$ 4.2.2]{KKP08} or \cite[Lemma 9.4.3]{Man22}),  given a \dgbva\,\,of degree $k$, a \dgla\,\,associated with it,  is then canonically defined. This DGLA is denoted by $(\mathfrak{g},\ddd_{\mathfrak{g}}),$ where $$\mathfrak{g}:=\AAA[k],\quad \ddd_{\mathfrak{g}}:=-\ddd_{\AAA},$$ and the bracket is defined by
		$$[a,b]:={(-1)}^{|{a}|}\bfitDelta(ab)-{(-1)}^{|{a}|}\bfitDelta(a)b-a\bfitDelta(b).$$
	\end{remark}
	\begin{definition}[{\cite[Definition 13.6.2]{Man22}}]\label{1212}
		A \dgbva\ of degree  $k$  is said to have the \textit{degeneration property} if for every \( a_0 \in \AAA \) such that \( \ddd a_0 = 0 \), there exists a sequence  $\{a_i\}_{i \geq 0},$ with   $|{a_{i+1}}|=|{a_{i}}|-k-1$
		and such that for any $i\geq 0,$
		$$\bfitDelta a_i=\ddd a_{i+1}.$$
	\end{definition}
	\begin{remark}
		For an equivalent definition of Definition \ref{1212} that relies on the $\mathbb{C}[[u]]$-module freeness of \( H^\bullet(\AAA[[u]], \ddd + u\bfitDelta) \), see \cite[Definition 4.1.3]{KKP08}.  Here  $u$  is a formal central variable of (even) degree $k+1$, and \( \AAA[[u]] \) denotes the graded vector space of formal power series with coefficients in \( \AAA \).
	\end{remark}

	The following theorem is of great significance.
	\begin{theorem} \label{()}
		Let $(\AAA,\ddd,\bfitDelta)$ be a \dgbva \,\,with the degeneration property. Then the associated \dgla \,\,$(\mathfrak{g},\ddd_{\mathfrak{g}})$ is homotopy abelian. 
		\begin{proof}
			See \cite[Theorem 4.14]{KKP08} or \cite[Theorem 2]{Ter08} for $k = 1$, \cite[Theorem 7.6]{Iac15} for any odd $k$. See also \cite{BK98} for a slightly weaker result.
		\end{proof}
	\end{theorem}
	\subsubsection{Proof of Theorem \ref{main thm}: algebraic approach}\label{><}
	\begin{proof}The approach follows the spirit of \cite[$\S$ 4.3.3]{KKP08} or \cite[$\S$  3]{Kon}. 
		In our situation, one can verify that the components of the appropriate DGBVA possessing the degeneration property are as follows:
		$$\AAA:=A^{0,\bullet}(X,\wedge^{\bullet}T_X^1(-\log D));$$
		$$\ddd:=\bp;$$
		$$\bfitDelta:={\mathbf{{i}}_{\Upomega}}^{-1}\circ D^\prime_{h^*} \circ \mathbf{{i}}_{\Upomega}.\footnote{The degree of this \dgbva\,\,is then $k=1$.}$$
		 Here $$\Upomega\in A^{0,0}(X,\Omega^n_X(\log D)\otimes L^*)$$ is a nowhere vanishing  section of  logarithmic $(n,0)$-form with values in $L^*,$ where
		\begin{equation*}
			L:=\Omega_{X}^n(\log D)\simeq K_X\otimes\mathcal{O}_X(D),
		\end{equation*}
		thus
		$$L=\sum_{i=1}^s \mathcal{O}_X(\qqq_iD_i)\in\text{Pic}_{\mb Q}(X)\,\,\,\text{with}\,\,\, \qqq_i:=1-\aaa_i\in[0,1];$$
		and 
		$$\mathbf{{i}}_{\Upomega}: \wedge^{\bullet}T_X^1(-\log D)\rightarrow \Omega_{X}^{n-\bullet}(\log D)\otimes L^*$$
		is the isomorphism given by the contraction with $\Upomega$. The degeneration property  follows from Theorem \ref{n,q-main}, or more practically, Theorem \ref{Thm 1.4}, as the contraction \( \mathbf{i}_{\Omega} \) provides an isomorphism of double complexes between the \dgbva\ \( (\AAA, \ddd, \bfitDelta) \) and the logarithmic Dolbeault-type double complex \((A^{0,\bullet}(X, \Omega_{X}^{\bullet}(\log D) \otimes L^*), \bp, D^\prime_{h^*})\).  Theorem \ref{()} then implies that the associated \dgla,
		$$(\mathfrak{g}, \ddd_{\mathfrak{g}})=(A^{0,\bullet}(X,\wedge^{\bullet}T_X^1(-\log D)),\bp)$$ 
		is homotopy abelian.   
		
		Consider the following \dgla,
		$$(\mathfrak{g}^\prime, \ddd_{\mathfrak{g}^\prime})=(A^{0,\bullet}(X, T_X^1(-\log D)),\bp),$$
		which controls  the locally trivial deformations of the pair $(X,D)$. 
		We then have  a natural inclusion of \dgla s
		$$(\mathfrak{g}^\prime, \ddd_{\mathfrak{g}^\prime})\hookrightarrow (\mathfrak{g}, \ddd_{\mathfrak{g}}),$$
		which embeds $(\mathfrak{g}^\prime, \ddd_{\mathfrak{g}^\prime})$ as a direct summand in 
		$(\mathfrak{g}, \ddd_{\mathfrak{g}}),$
		and so induces an embedding
		$$H^\bullet(\mathfrak{g}^\prime,\ddd_{\mathfrak{g}^\prime})\subset H^\bullet(\mathfrak{g},\ddd_{\mathfrak{g}})$$
		in cohomology.
		By \cite[Proposition 4.11 \textbf{(ii)}]{KKP08}, $(\mathfrak{g}^\prime, \ddd_{\mathfrak{g}^\prime})$ is also homotopy abelian.  
		
		Consequently, the formal moduli space associated with the \dgla\
		$(\mathfrak{g}', \ddd_{\mathfrak{g}'})$ is smooth. Combined with Artin's approximation theorem \cite{Art68}
		(see also \cite[Theorem B.1]{RS20}), this implies the desired unobstructedness.
		
		The proof of Theorem \ref{main thm} is, therefore, complete.
	\end{proof}
	
	\appendix
	\renewcommand{\thesection}{\Roman{section}}

		\section{Proof of Theorem \ref{main thm}: analytic approach (joint with Sheng Rao)}\label{jw}
\begin{center}
	\scriptsize SHENG RAO 
	\alphafootnote{School of Mathematics and Statistics, Wuhan University,
		Wuhan 430072, China; Email address: \texttt{likeanyone@whu.edu.cn}
		} \&  RUNZE ZHANG
\end{center}

\medskip

	In this appendix, we aim to reprove Theorem \ref{main thm} from an analytic point of view.
	Under the assumptions of Theorem \ref{main thm}, for any initial cohomology class
	$[\varphi_1]\in H^{1}(X,T_X^1(-\log D))$ and for any parameter
	$t$ in a sufficiently small $\epsilon$-disk $\Delta_\epsilon$ centered at the origin of
	$\mathbb{C}^{d}$, where $d=\dim_{\mathbb{C}} H^{1}(X, T_X^1(-\log D))$,
	we construct a holomorphic family
	\[
	\varphi=\varphi(t)\in A^{0,1}(X, T_X^1(-\log D))
	\]
	satisfying the following integrability and initial conditions:
	\begin{equation}
		\tag{I.I.}\label{log 0.1}
		\bar{\partial}\varphi=\frac{1}{2}[\varphi,\varphi],
		\qquad
		\left.\frac{\partial \varphi}{\partial t}\right|_{t=0}=\varphi_1.
	\end{equation}
	This construction implies the unobstructedness of locally trivial deformations of $(X,D)$.
	
	The proof is divided into three steps.
	\setcounter{step}{0}
\begin{step}
	\label{step 1}
	{\textrm{Explicit solution to the twisted logarithmic $\bar\partial$-equation \ref{solution 2}.}}
\end{step}
Let $(X,D)$ and $L$ be the same as in  Theorem \ref{n,q-main}. Inspired by \cite[$\S$ 2]{LRW19}, we
set
\begin{equation}\label{EEE}
	\mathpzc{E}^p := \Omega_X^p(\log D)\otimes L^*
\end{equation}
as a holomorphic vector bundle. Then (global) $L^*$-valued logarithmic
$(p,q)$-forms can be identified with $\mathpzc{E}^p$-valued $(0,q)$-forms:

\[
A^{0,q}\bigl(X,\Omega_X^p(\log D)\otimes L^*\bigr)
\simeq
A^{0,q}(X,\mathpzc{E}^p).
\]

Fix an arbitrary smooth Hermitian metric $g_p$ on
$\mathpzc E^p$, together with a K\"ahler metric $\omega$ on
$X$. Let $\mathtt D_{g_p} = \mathtt D^{1,0}_{g_p} + \bar{\partial}$ be the Chern connection of
$(\mathpzc E^p, g_p)$, decomposed into its $(1,0)$ and $(0,1)$ parts. Then the following \textit{$\mathpzc E^p$-valued Hodge decomposition} holds:
\begin{equation}\label{eq:Hodge-decomposition}
	\mathpzc I
	=
	\mathtt H^{\prime\prime}_{\mathpzc E^p}
	+\bar\partial\bar\partial^{*}_{\mathpzc E^p}\mathtt G^{\prime\prime}_{\mathpzc E^p}
	+\bar\partial^{*}_{\mathpzc E^p}\bar\partial\mathtt G^{\prime\prime}_{\mathpzc E^p},
\end{equation}
where $\mathpzc I$ denotes the identity operator, and
$\mathtt H^{\prime\prime}_{\mathpzc E^p}$ (resp.\
$\mathtt G^{\prime\prime}_{\mathpzc E^p}$, $\bar\partial^{*}_{\mathpzc E^p}$)
denotes the harmonic projection (resp. Green operator, adjoint operator of
$\bar\partial$) with respect to the metrics $g_p$ and $\omega$. 
 Here the harmonic projection $\mathtt H^{\prime\prime}_{\mathpzc E^p}$ is defined by
\begin{equation}\label{Harmonic pro}
	\mathtt H^{\prime\prime}_{\mathpzc E^p}(\beta) := \sum_i \langle \beta, \zeta_i \rangle_{\mathpzc E^p}\,\zeta_i,
	\qquad \text{for any } \beta \in A^{r,s}(X,\mathpzc E^p),
\end{equation}
where $\{\zeta_i\}_i$ is an orthonormal basis of $\Delta^{\prime\prime}_{\mathpzc E^p}$-harmonic
$(r,s)$-forms with values in $\mathpzc E^p$, and
$\langle \cdot, \cdot \rangle_{\mathpzc E^p}$ is the natural Hermitian
inner product on $A^{r,s}(X,\mathpzc E^p)$ induced by $g_p$ and $\omega$.

With the aid of \eqref{eq:Hodge-decomposition}, and inspired by the proof
of the existence theorem for deformations in
\cite[$\S$ 4.2, Proposition 2.3, (11)]{MK71} (see also \cite{KNS58}),
we are able to single out  a distinguished log solution $\chi$ to
\ref{solution 2}.
\begin{lemma}\label{lemma-exp}
	The log solution in \ref{solution 2} can be chosen as
	\begin{equation*}
		\chi=\bp^{*}_{\mathpzc E^p}\mathtt G^{\prime\prime}_{\mathpzc E^p}(D^\prime_{h^{*}}\alpha).
	\end{equation*}
	\begin{proof}
		Applying the Hodge decomposition \eqref{eq:Hodge-decomposition} to
		$D^\prime_{h^{*}}\alpha \in A^{0,q}(X,\mathpzc E^p)$, we obtain
		\begin{equation}\label{55}
			D^\prime_{h^{*}}\alpha
			=
			\mathtt H^{\prime\prime}_{\mathpzc E^p}(D^\prime_{h^{*}}\alpha)
			+ \bar\partial \bar\partial^{*}_{\mathpzc E^p} \mathtt G^{\prime\prime}_{\mathpzc E^p}(D^\prime_{h^{*}}\alpha)
			+ \bar\partial^{*}_{\mathpzc E^p} \bar\partial \mathtt G^{\prime\prime}_{\mathpzc E^p}(D^\prime_{h^{*}}\alpha).
		\end{equation}
		
		Since $D^\prime_{h^{*}}\alpha$ is $\bar\partial$-exact, as ensured by 
		Theorem \ref{n,q-main}, the first term in \eqref{55} vanishes by the definition of the harmonic projection in 
		\eqref{Harmonic pro}. Moreover, the last term in \eqref{55} also vanishes
		due to the commutativity of $\bar\partial$ and
		$\mathtt G^{\prime\prime}_{\mathpzc E^p}$.
		It then follows that
		\[
		D^\prime_{h^{*}}\alpha = \bar\partial \bar\partial^{*}_{\mathpzc E^p} \mathtt G^{\prime\prime}_{\mathpzc E^p}(D^\prime_{h^{*}}\alpha),
		\]
		which gives the desired conclusion.
	\end{proof}
\end{lemma}
We have thus completed Step \ref{step 1}.

\vspace{1em}
With Lemma \ref{lemma-exp} at hand, we are now in a good position to prove
Theorem \ref{main thm}. The proof is essentially the same as that of
\cite[Theorem 4.12]{LRW19}. For the convenience of the reader, we nevertheless
provide a complete proof in the remaining two steps,
Steps \ref{step 2} and \ref{step 3}.

\begin{step}
	\label{step 2}
	{\textrm{Solve the  equations \eqref{log 0.1} recursively.}}
\end{step}

Let
\begin{equation}\label{LLL}
	L:=\Omega_X^n(\log D)\simeq K_X\otimes\mathcal{O}_X(D)
\end{equation}
be the associated holomorphic line bundle. By \eqref{KKP's setting}, we may write
\[
L=\sum_{i=1}^s \mathcal{O}_X({q}_i D_i)\in \mathrm{Pic}_{\mathbb{Q}}(X),
\quad\text{where }
{q}_i:=1-a_i\in[0,1]\cap\mb Q.
\]

Let 
\[
\Upomega \in H^0\bigl(X, \Omega_X^n(\log D) \otimes L^*\bigr)
\]
be a nowhere vanishing \emph{holomorphic} logarithmic $(n,0)$-form with values in $L^*$. 

\vspace{1em}
We will need a logarithmic version of the Tian--Todorov lemma adapted to our setting.

\begin{lemma}\label{Tian-Todorov}
	For any $\varphi,\psi\in A^{0,1}\bigl(X, T_X^1(-\log D)\bigr)$ and any
	$\Omega\in A^{0,q}\bigl(X,\Omega_X^n(\log D)\otimes L^*\bigr)$, one has
	\[
	[\varphi,\psi]\lrcorner\, \Omega
	=
	- D^\prime_{h^{*}}\bigl(\psi\lrcorner(\varphi\lrcorner\Omega)\bigr)
	+ \psi\lrcorner D^\prime_{h^{*}}\bigl(\varphi\lrcorner\Omega\bigr)
	+ \varphi\lrcorner D^\prime_{h^{*}}\bigl(\psi\lrcorner\Omega\bigr),
	\]
	where $$\varphi\l=i_{\varphi}: A^{0,q}(X,\Omega_X^p(\log D)\otimes L^*)\rightarrow A^{0,q+1}(X,\Omega_X^{p-1}(\log D)\otimes L^*)$$ is the contraction map.
\end{lemma}

\begin{proof}
	Locally, the operator $D^\prime_{h^{*}}$ admits the expression
	\begin{equation*}
		D^\prime_{h^{*}}
		=
		\partial+\sum_{i=1}^s \qqq_i\,\frac{dz^i}{z^i}\wedge \bullet,
	\end{equation*}
	where $D_i=\{z^i=0\}$ for $i=1,\ldots,s$.
	Comparing with \cite[Lemma 4.4]{LRW19} (which in turn follows from
	\cite[Proposition 3.2]{LR12}), it suffices to verify that
	\[
	\sum_{i=1}^s \qqq_i\,\frac{dz^i}{z^i}\wedge
	\bigl(\psi\lrcorner(\varphi\lrcorner\Omega)\bigr)
	=
	\psi\lrcorner\Bigl(\sum_{i=1}^s \qqq_i\,\frac{dz^i}{z^i}\Bigr)\wedge
	(\varphi\lrcorner\Omega)
	+
	\varphi\lrcorner\Bigl(\sum_{i=1}^s \qqq_i\,\frac{dz^i}{z^i}\Bigr)\wedge
	(\psi\lrcorner\Omega).
	\]
	This identity can be checked directly.
\end{proof}
\begin{proposition}\label{prop-1}
	Suppose that there exist two smooth families
	\[
	\varphi(t)\in A^{0,1}\bigl(X,T_X^1(-\log D)\bigr)
	\quad\text{and}\quad
	\Omega(t)\in A^{0,0}\bigl(X,\Omega_X^{n}(\log D)\otimes L^*\bigr)
	\]
	satisfying the following system of equations:
	\begin{equation}\tag*{($\diamondsuit$)}\label{log 0.2}
		\begin{cases}
			\bigl(\bar{\partial}+\dfrac{1}{2}D^\prime_{h^{*}}\circ i_{\varphi}\bigr)
			\bigl(i_{\varphi}\Omega(t)\bigr)=0,\\[0.6em]
			\bigl(\bar{\partial}+D^\prime_{h^{*}}\circ i_{\varphi}\bigr)\Omega(t)=0,\\
			\Omega(0)=\Upomega,
		\end{cases}
	\end{equation}
	Then $\varphi(t)$ satisfies the integrability condition in  \eqref{log 0.1} for any sufficiently small $t$.
\end{proposition}
\begin{proof}			Combining \ref{log 0.2} with Lemma \ref{Tian-Todorov}, we compute
	\begin{align*}\label{log 0.3}
		\begin{split}
			\bar{\partial}\bigl(\varphi\lrcorner \Omega(t)\bigr)
			&=-\frac{1}{2}D^\prime_{h^{*}}\bigl(i_{\varphi}\circ i_{\varphi}\Omega(t)\bigr)\\
			&=\frac{1}{2}[\varphi,\varphi]\lrcorner \Omega(t)
			- i_{\varphi}D^\prime_{h^{*}}\bigl(i_{\varphi}\Omega(t)\bigr)\\
			&=\frac{1}{2}[\varphi,\varphi]\lrcorner \Omega(t)
			+ i_{\varphi}\bar{\partial}\Omega(t).
		\end{split}
	\end{align*}
	Hence
	\[
	(\bar{\partial}\varphi)\lrcorner \Omega(t)
	=
	\bar{\partial}\bigl(\varphi\lrcorner \Omega(t)\bigr)
	- i_{\varphi}\bar{\partial}\Omega(t)
	=
	\frac{1}{2}[\varphi,\varphi]\lrcorner \Omega(t).
	\]
	
	Since $\Omega(t)$ depends smoothly on $t$ and satisfies
	$\Omega(0)=\Omega_0=\Upomega$, which is nowhere vanishing, it follows that
	$\Omega(t)$ has no zeros for any sufficiently small $t$.
	This completes the proof.
\end{proof}

We denote by \ref{log 0.2}$_k$ the subsystem of \ref{log 0.2} consisting of the
equations of order $(k+1)$ from the first equation and of order $k$ from the second.
We now solve \ref{log 0.2} using an \emph{iterative approach}.

For $k=1$, the system   \ref{log 0.2}$_1$ reads
\begin{equation}\label{123}
	\begin{cases}
		\bar{\partial}\bigl(i_{\varphi_2}\Omega_0 + i_{\varphi_1}\Omega_1\bigr)
		+ \frac{1}{2} D^{\prime}_{h^*} \circ i_{\varphi_1}(i_{\varphi_1}\Omega_0) = 0,\\[0.3em]
		\bar{\partial}\Omega_1 + D^{\prime}_{h^*} \circ i_{\varphi_1}\Omega_0 = 0.
	\end{cases}
\end{equation}
Examining the second equation in \eqref{123}, we observe that
\[
\bar{\partial} (D^\prime_{h^*} \circ i_{\varphi_1}\Omega_0)
= - D^\prime_{h^*} \bar{\partial}(\varphi_1 \lrcorner \Omega_0)
= - D^\prime_{h^*} (\bar{\partial} \varphi_1 \lrcorner \Omega_0 + \varphi_1 \lrcorner \bar{\partial} \Omega_0) = 0,
\]
since both $\Omega_0 = \Upomega$ and $\varphi_1$ are $\bar{\partial}$-closed.  
Hence, the conditions of Theorem \ref{n,q-main} are satisfied. By Lemma \ref{lemma-exp}, we can select
\[
\Omega_1 = \bar{\partial}_{\mathpzc{E}^n}^* \, \mathtt G''_{\mathpzc{E}^n}
\bigl( - D^\prime_{h^*} \circ i_{\varphi_1} \Omega_0 \bigr)
\in A^{0,0}(X, \Omega_X^n(\log D) \otimes L^*),
\]
where $$\mathpzc{E}^n= \Omega_X^n(\log D) \otimes L^*\simeq \mathcal{O}_X$$ by  \eqref{EEE} together with \eqref{LLL}.  
Similarly, the first equation in \eqref{123} can be solved in the same manner. 
More precisely, there exists
\[
\varphi_2 \in A^{0,1}(X, T_X^1(-\log D))
\]
such that
\[
i_{\varphi_2}\Omega_0
=
-\, i_{\varphi_1}\Omega_1
+ \bar{\partial}_{\mathpzc E^{n-1}}^{\,*}\,\mathtt G''_{\mathpzc E^{n-1}}
\Bigl(
- \frac{1}{2} D^\prime_{h^*} \circ i_{\varphi_1}\bigl(i_{\varphi_1}\Omega_0\bigr)
\Bigr)
\in A^{0,1}\bigl(X, \Omega_X^{n-1}(\log D)\otimes L^*\bigr),
\]
where
\[
\mathpzc E^{n-1}= \Omega_X^{n-1}(\log D)\otimes L^* .
\]
Here we use the fact that the contraction map
\begin{equation}\label{a}
	\bullet \l  \Upomega :
	A^{0,1}\bigl(X, T_X^1(-\log D)\bigr)
	\longrightarrow
	A^{0,1}\bigl(X, \Omega_X^{n-1}(\log D)\otimes L^*\bigr)
\end{equation}
is an isomorphism, since the logarithmic volume form $\Upomega$ has no zeroes, see \cite[Lemma 4.9]{LRW19}.

Consequently, $\varphi_2$ can be written explicitly as
\begin{equation}\label{hg}
	\varphi_2
	=
	\Omega_0^* \l
	\Bigl(
	- i_{\varphi_1}\Omega_1
	+ \bar{\partial}_{\mathpzc E^{n-1}}^{\,*}\,\mathtt G''_{\mathpzc E^{n-1}}
	\Bigl(
	- \frac{1}{2} D^\prime_{h^*} \circ i_{\varphi_1}\bigl(i_{\varphi_1}\Omega_0\bigr)
	\Bigr)
	\Bigr),
\end{equation}
where
\[
\Upomega^* \l \bullet :
A^{0,1}\bigl(X, \Omega_X^{n-1}(\log D)\bigr)
\longrightarrow
A^{0,1}\bigl(X, T_X^1(-\log D)\bigr)
\]
is the inverse of the map in \eqref{a}. For an explicit construction of this inverse, see \cite[(4.18)]{LRW19}.  Notice that we have the nice properties:
for any 
\(\psi \in A^{0,1}(X, \Omega_X^{n-1}(\log D))\) and 
\(\phi \in A^{0,1}(X, T_X^1(-\log D))\), 
\begin{equation}\label{66666}
	(\Upomega^* \lrcorner \, \psi) \lrcorner \, \Upomega = \psi, 
	\qquad 
	\Upomega^* \lrcorner \, (\phi \lrcorner \, \Upomega) = \phi.
\end{equation}

By induction, assume that \ref{log 0.2}  is solvable for all $k \le N$, and that we have already found
\[
\Omega_i, \, 0 \le i \le N, \quad \text{and} \quad \varphi_i, \, 1 \le i \le N+1.
\]
Thus, Proposition \ref{prop-1} implies
\begin{equation}\label{log 0.4}
	\bar{\partial} \varphi_k = \frac{1}{2} \sum_{i+j=k} [\varphi_i, \varphi_j], \quad k \le N+1,
\end{equation}
and then \cite[Lemma 4.11]{LRW19} gives
\begin{equation}\label{log 0.5}
	\bar{\partial} \Bigl( \sum_{i+j=k} [\varphi_i, \varphi_j] \Bigr) = 0, \quad k \le N+2.
\end{equation}

For the $(N+1)$-th step, \ref{log 0.2}$_{N+1}$  reads
\begin{equation}\label{log 0.6}
	\begin{cases}
		\sum_{i+j=N+2} \bar{\partial} i_{\varphi_i} \Omega_j
		+ \sum_{i+j+k=N+2} \frac{1}{2} D^\prime_{h^*} \circ i_{\varphi_i} \circ i_{\varphi_j} \Omega_k = 0,\\[0.3em]
		\bar{\partial} \Omega_{N+1} + \sum_{i+j=N+1} D^\prime_{h^*} \circ i_{\varphi_i} \Omega_j = 0.
	\end{cases}
\end{equation}

For the second equation  in \eqref{log 0.6}, we have 
\allowdisplaybreaks
\begin{align*}
	-\bar\partial \sum_{i+j=N+1} D^\prime_{h^{*}}(\varphi_i \lrcorner \Omega_j)
	&= D^\prime_{h^{*}}\Bigg(
	\sum_{i=1}^{N+1} \bar\partial \varphi_i \lrcorner \Omega_{N+1-i}
	+ \sum_{i=1}^{N+1} \varphi_i \lrcorner \bar\partial \Omega_{N+1-i}
	\Bigg) \\[0.3em]
	&= D^\prime_{h^{*}}\Bigg(
	\frac12 \sum_{i=1}^{N+1} \sum_{j=1}^{i-1}
	[\varphi_j,\varphi_{i-j}] \lrcorner \Omega_{N+1-i} \\
	&\qquad\qquad
	- \sum_{i=1}^{N+1} \varphi_i \lrcorner
	D^\prime_{h^{*}}\Bigg(
	\sum_{j=1}^{N+1-i}
	\varphi_j \lrcorner \Omega_{N+1-i-j}
	\Bigg)
	\Bigg) \\[0.3em]
	&= D^\prime_{h^{*}}\Bigg(
	\frac12 \sum_{i=1}^{N+1} \sum_{j=1}^{i-1}
	\Bigg(
	- D^\prime_{h^{*}}\Big(
	\varphi_j \lrcorner
	(\varphi_{i-j} \lrcorner \Omega_{N+1-i})
	\Big) \\
	&\qquad\qquad\qquad
	+ \varphi_j \lrcorner
	D^\prime_{h^{*}}(\varphi_{i-j} \lrcorner \Omega_{N+1-i})
	+ \varphi_{i-j} \lrcorner
	D^\prime_{h^{*}}(\varphi_j \lrcorner \Omega_{N+1-i})
	\Bigg) \\
	&\qquad
	- \sum_{i=1}^{N+1} \varphi_i \lrcorner
	D^\prime_{h^{*}}\Bigg(
	\sum_{j=1}^{N+1-i}
	\varphi_j \lrcorner \Omega_{N+1-i-j}
	\Bigg)
	\Bigg) \\[0.3em]
	&= D^\prime_{h^{*}}\Bigg(
	\sum_{1 \le j < i \le N+1}
	\varphi_j \lrcorner
	D^\prime_{h^{*}}(\varphi_{i-j} \lrcorner \Omega_{N+1-i}) \\
	&\qquad
	- \sum_{i=1}^{N+1} \varphi_i \lrcorner
	D^\prime_{h^{*}}\Bigg(
	\sum_{j=1}^{N+1-i}
	\varphi_j \lrcorner \Omega_{N+1-i-j}
	\Bigg)
	\Bigg) \\
	&= 0,
\end{align*}
where the third equality holds thanks to Lemma \ref{Tian-Todorov}.
We then can set 
\begin{equation}\label{Omega}
	\Omega_{N+1} = \sum_{i+j=N+1} \bar{\partial}_{\mathpzc E^n}^* \, \mathtt G''_{\mathpzc E^n}
	\bigl( - D^\prime_{h^*} \circ i_{\varphi_i} \Omega_j \bigr)
	\in A^{0,0}(X, \Omega_X^n(\log D) \otimes L^*)
\end{equation}
by using Lemma \ref{lemma-exp}.

Moreover, for the first equation  in \eqref{log 0.6}, one has
\begin{align*}
	&\quad \bar{\partial} \sum_{i+j+k=N+2} \frac{1}{2} D^\prime_{h^*} \circ i_{\varphi_i}\circ i_{\varphi_j} \Omega_k \\
	&= \sum_{i+j+k=N+2} \bar{\partial} \Bigl(
	- [\varphi_j, \varphi_i] \lrcorner \Omega_k
	+ \varphi_i \lrcorner D^\prime_{h^*} (\varphi_j \lrcorner \Omega_k)
	+ \varphi_j \lrcorner D^\prime_{h^*} (\varphi_i \lrcorner \Omega_k)
	\Bigr) \\
	&= - \sum_{i+j+k=N+2} \bar{\partial} ([\varphi_j, \varphi_i] \lrcorner \Omega_k)
	- 2 \sum_{j+k=N+2,\, k \ge 1} \bar{\partial} (\varphi_j \lrcorner \bar{\partial} \Omega_k) \\
	&= - \sum_{k=0}^{N+2} \Bigl( \bar{\partial} \sum_{i+j=N+2-k} [\varphi_j, \varphi_i] \Bigr) \lrcorner \Omega_k
	+ 2 \sum_{k=1}^{N+2} \Bigl( \frac{1}{2} \sum_{i+j=N+2-k} [\varphi_j, \varphi_i]
	- \bar{\partial} \varphi_{N+2-k} \Bigr) \lrcorner \bar{\partial} \Omega_k \\
	&= 0,
\end{align*}
where the first equality is due to   Lemma \ref{Tian-Todorov},
and the last equality uses \eqref{log 0.4} and \eqref{log 0.5}.
Consequently, we can also choose $\varphi_{N+2} \in A^{0,1}(X, T_X^1(-\log D))$ such that
\begin{align}\label{log varphi}
	i_{\varphi_{N+2}} \Omega_0
	&= - \sum_{j=1}^{N+2} i_{\varphi_{N+2-j}} \Omega_j
	+ \bar{\partial}_{\mathpzc E^{n-1}}^* \, \mathtt G''_{\mathpzc E^{n-1}}
	\Bigl( - \sum_{i+j+k=N+2} \frac{1}{2} D^\prime_{h^*} \circ i_{\varphi_i} \circ i_{\varphi_j} \Omega_k \Bigr) \\
	&\in A^{0,1}(X, \Omega_X^{n-1}(\log D) \otimes L^*). \notag
\end{align}

In conclusion, we have constructed
$$
\Omega_i \in A^{0,0}(X, \Omega_X^n(\log D) \otimes L^*), \, 0 \le i \le N+1,
\,\, \text{and} \,\,
\varphi_i \in A^{0,1}(X, T_X^1(-\log D)), \, 1 \le i \le N+2,
$$
which solve \ref{log 0.2}$_{N+1}$. Therefore, for sufficiently small $t$, \ref{log 0.2} can be solved recursively.
This completes Step \ref{step 2}.
\begin{step}
	\label{step 3}
	{\textrm{Regularity procedures for $\varphi(t)$ and $\Omega(t)$.}}	
\end{step}
Fix an integer $k \ge 2$ and a real number $\alpha \in (0,1)$. 
We denote by $\|\cdot\|_{k,\alpha}$ the H\"older norm and by $\|\cdot\|_{0}$ 
the $C^{0}$-norm of bundle-valued $(0,q)$-forms on $X$; 
see \cite[p. 275]{Kod86}.
We first prove the convergence of the formal power series $\varphi(t)$ and $\Omega(t)$ under the $C^{k,\alpha}$-norm.

\vspace{0.3em}
 For simplicity, we will always assume \( d =  \dim_{\mathbb{C}} H^{1}(X, T_X^1(-\log D)) \)=1 in what follows; the case \( d > 1 \) can be treated in a similar way.

\vspace{0.3em}
Consider the following power series, which plays an important role in the convergence analysis in deformation theory of complex structures:
\[
\mathrm {A}(t) = \frac{\mathrm b}{16\mathrm c} \sum_{m=1}^{\infty} \frac{(\mathrm c t)^m}{m^2} = \sum_{m=1}^{\infty} \mathrm a_m t^m, \quad \mathrm a_m = \frac{\mathrm b \mathrm c^{m-1}}{16 m^2},
\]
where $\mathrm b,\mathrm c>0$.  This series satisfies
\[
\mathrm A^n(t) \le \left(\frac{\mathrm b}{\mathrm c}\right)^{n-1} \mathrm A(t),
\]
and converges for $|t|<1/\mathrm c$.

Suppose that the coefficients $\varphi_i$ and $\Omega_i$ have been chosen so that
\[
\|\varphi_i\|_{k,\alpha} \le \mathrm a_i, \quad 1 \le i \le N+1, \qquad
\|\Omega_i\|_{k,\alpha} \le \left(\frac{\mathrm c}{\mathrm b}\right)^{1/2} \mathrm a_i, \quad 1 \le i \le l.
\]
Then from \eqref{Omega}, we obtain
\begin{align}\label{Holder-Omega}
	\begin{split}
		\|\Omega_{N+1}\|_{k,\alpha} 
		&\le C \sum_{i+j=N+1} \|\varphi_i\|_{k,\alpha} \|\Omega_j\|_{k,\alpha} \\
		&\le C \left( \left(\frac{\mathrm c}{\mathrm b}\right)^{1/2} \sum_{i+j=N+1} \mathrm a_i \mathrm a_j + \|\Omega_0\|_{k,\alpha} \mathrm a_{N+1} \right) \\
		&\le C \left( \frac{\mathrm b}{\mathrm c} + \left(\frac{\mathrm b}{\mathrm c}\right)^{1/2} \|\Omega_0\|_{k,\alpha} \right) \left( \frac{\mathrm c}{\mathrm b} \right)^{1/2} \mathrm a_{N+1}.
	\end{split}
\end{align}
Similarly to \eqref{hg}, from \eqref{log varphi}, we can write
\[
\varphi_{N+2} = \Omega_0^* \l 
\Biggl(
- \sum_{j=1}^{N+2} i_{\varphi_{N+2-j}} \Omega_j
+ \bar{\partial}_{\mathpzc E^{n-1}}^* \, \mathtt G''_{\mathpzc E^{n-1}}
\Bigl( - \sum_{i+j+k=N+2} \frac{1}{2} D^\prime_{h^*} \circ i_{\varphi_i} \circ i_{\varphi_j} \Omega_k \Bigr)
\Biggr).
\]
It then follows that 
\begin{align}\label{Holder-varphi}
	\begin{split}
		\|\varphi_{N+2}\|_{k,\alpha} 
		&\le C \Biggl( 
		\frac{\mathrm b}{\mathrm c} \Bigl(\frac{\mathrm c}{\mathrm b}\Bigr)^{1/2} 
		+ \frac{\mathrm b}{\mathrm c} \|\Omega_0\|_{k,\alpha} 
		+ \Bigl(\frac{\mathrm c}{\mathrm b}\Bigr)^{1/2} \Bigl(\frac{\mathrm b}{\mathrm c}\Bigr)^2
		\Biggr)\mathrm  a_{N+2} \\
		&\le C \Biggl(
		\Bigl(\frac{\mathrm b}{\mathrm c}\Bigr)^{1/2} 
		+ \frac{\mathrm b}{\mathrm c} \|\Omega_0\|_{k,\alpha} 
		+ \Bigl(\frac{\mathrm b}{\mathrm c}\Bigr)^{3/2}
		\Biggr) \mathrm a_{N+2}.
	\end{split}
\end{align}

We may now choose $\mathrm b/\mathrm c$ sufficiently small so that
\[
C \left( \frac{\mathrm b}{\mathrm c} + \left(\frac{\mathrm b}{\mathrm c}\right)^{1/2} \|\Omega_0\|_{k,\alpha} \right) \le 1, \quad
C \left( \left(\frac{\mathrm b}{\mathrm c}\right)^{1/2} + \frac{\mathrm b}{\mathrm c} \|\Omega_0\|_{k,\alpha} + \left(\frac{\mathrm b}{\mathrm c}\right)^{3/2} \right) \le 1.
\]
Then from \eqref{Holder-Omega} and \eqref{Holder-varphi}, we deduce
\[
\|\Omega_{N+1}\|_{k,\alpha} \le \left(\frac{\mathrm c}{\mathrm b}\right)^{1/2} \mathrm a_{N+1}, \quad
\|\varphi_{N+2}\|_{k,\alpha} \le \mathrm a_{N+2}.
\]

As  a result,  the two power series
\[
\varphi(t) = \sum_{i=1}^\infty \varphi_i t^i, \quad
\Omega(t) = \Omega_0 + \sum_{i=1}^\infty \Omega_i t^i
\]
converge in the $C^{k,\alpha}$-norm for $|t|<1/\mathrm c$, with
\[
\|\varphi(t)\|_{k,\alpha} \le \mathrm A(t), \quad
\|\Omega(t)\|_{k,\alpha} \le \|\Omega_0\|_{k,\alpha} + \left(\frac{\mathrm c}{\mathrm b}\right)^{1/2} \mathrm A(t).
\]

\vspace{2em}
Finally, we address the regularity of $\varphi=\varphi(t)$ and $\Omega(t)$. 
From \eqref{Omega} and \eqref{log varphi}, we obtain
\begin{equation}\label{log 1.3}
	\begin{cases}
		\bigl(\mathpzc I+ \bar{\p}_{\mathpzc E^n}^{\,*}\mathtt G''_{\mathpzc E^n} D^\prime_{h^*} i_{\varphi}\bigr)\Omega(t)=\Omega_0,\\[0.3em]
		\bigl(\mathpzc I+\dfrac{1}{2}\bar{\p}_{\mathpzc E^{n-1}}^{\,*}\mathtt G''_{\mathpzc E^{n-1}} D^\prime_{h^*} i_{\varphi}\bigr)i_{\varphi}\Omega(t)
		=i_{\varphi_1}\Omega_0.
	\end{cases}
\end{equation}

From the first equation of \eqref{log 1.3}, we deduce that
\[
\frac{\p \Omega(t)}{\p \bar t}=0.
\]
Consequently, we obtain
\begin{equation}\label{log 1.5}
	\begin{cases}
		\dfrac{\p^2}{\p t \p \bar t}\Omega(t)
		+\Delta^{\prime\prime}_{\mathpzc E^n}\Omega(t)
		=
		\bar{\p}\,\bar{\p}_{\mathpzc E^n}^{\,*}\Omega_0
		-\bar{\p}_{\mathpzc E^n}^{\,*} D^\prime_{h^*} i_{\varphi}\Omega(t),\\[0.4em]
		\dfrac{\p^2}{\p t \p \bar t}\bigl(i_{\varphi}\Omega(t)\bigr)
		+\Delta^{\prime\prime}_{\mathpzc E^{n-1}}\bigl(i_{\varphi}\Omega(t)\bigr)
		=
		\Delta^{\prime\prime}_{\mathpzc E^{n-1}}\bigl(i_{\varphi_1}\Omega_0\bigr)
		-\dfrac{1}{2}\bar{\p}_{\mathpzc E^{n-1}}^{\,*} D^\prime_{h^*} i_{\varphi}\circ i_{\varphi}\Omega(t).
	\end{cases}
\end{equation}

Since both $\varphi$ and $\Omega(t)$ converge in the $C^{k,\alpha}$-norm, it follows that
\begin{equation}\label{log 1.8}
	\| i_{\varphi}\Omega(t) \|_{k,\alpha} \le C.
\end{equation}
From the second equation of \eqref{log 1.5}, we further obtain
\begin{equation}\label{log 1.6}
	\bar{\p}_{\mathpzc E^{n-1}}^{\,*}\bigl(i_{\varphi}\Omega(t)\bigr)
	=
	\bar{\p}_{\mathpzc E^{n-1}}^{\,*}\bigl(i_{\varphi_1}\Omega_0\bigr).
\end{equation}
Noting that $[\bar{\p}_{\mathpzc E^{n-1}}^{\,*}, D^\prime_{h^*}]$ is a first-order differential operator, we deduce from \eqref{log 1.6} that
\begin{align*}
	\bigl\| \bar{\p}_{\mathpzc E^{n-1}}^{\,*} D^\prime_{h^*} i_{\varphi}\Omega(t) \bigr\|_{k-1,\alpha}
	&=
	\bigl\| [\bar{\p}_{\mathpzc E^{n-1}}^{\,*}, D^\prime_{h^*}]\bigl(i_{\varphi}\Omega(t)\bigr)
	- D^\prime_{h^*} \bar{\p}_{\mathpzc E^{n-1}}^{\,*}\bigl(i_{\varphi}\Omega(t)\bigr) \bigr\|_{k-1,\alpha} \\
	&\le C \| i_{\varphi}\Omega(t) \|_{k,\alpha}
	+ \| D^\prime_{h^*} \bar{\p}_{\mathpzc E^{n-1}}^{\,*}(i_{\varphi_1}\Omega_0) \|_{k-1,\alpha}
	\le C.
\end{align*}

By the first equation of \eqref{log 1.5}, we conclude that
\begin{equation} \label{sq}
	\| \Omega(t) \|_{k+1,\alpha} \le C.
\end{equation}
Recall that $\Omega(t)$ has the local expression
\[
\Omega(t)= v(t)\, \frac{dz^1\wedge\cdots\wedge dz^n}{z^1\cdots z^s},
\]
where $v(t)$ is a smooth function  and $D=\{z^1\cdots z^s=0\}$ locally.
Then the above estimate \eqref{sq} implies that, for each coordinate polydisk $U_i$, one has
$|v(t)|^{U_i}_{k+1,\alpha} \le C$.
After possibly shrinking the disc, we may assume that $|v(t)|^{U_i}_0 \ge c^\prime>0$, since $v(0)$ has no zeroes. A direct computation then yields
\begin{equation}\label{sg}
	\left| \frac{1}{v(t)} \right|^{U_i}_{k+1,\alpha} \le C.
\end{equation}

Consequently, by \eqref{66666} and \eqref{sg}, for such $t$ we may regard
\begin{equation*}\label{log 1.7}
	-\frac{1}{2}\bar{\p}_{\mathpzc E^{n-1}}^{\,*} D^\prime_{h^*} i_{\varphi}\circ i_{\varphi}\Omega(t)
	=
	-\frac{1}{2}\bar{\p}_{\mathpzc E^{n-1}}^{\,*} D^\prime_{h^*}
	\bigl(\Omega(t)^* \l i_{\varphi}\Omega(t)\bigr)
	\l i_{\varphi}\Omega(t)
\end{equation*}
as a linear second-order differential operator in $i_{\varphi}\Omega(t)$ with $C^{k-1,\alpha}$ coefficients.

Since $\varphi(0)=0$, the second equation of \eqref{log 1.5} can therefore be viewed, for $t$ sufficiently small, as a linear elliptic equation with $C^{k-1,\alpha}$ coefficients. By \cite[Theorem 6.17]{GT01}, we obtain
\begin{equation}\label{log 1.9}
	\| i_{\varphi}\Omega(t) \|_{k+1,\alpha} \le C.
\end{equation}

Iterating the above argument from \eqref{log 1.8} to \eqref{log 1.9}, we conclude that both $\Omega(t)$ and $i_{\varphi}\Omega(t)$ are smooth. Consequently,
\[
\varphi(t)=\Omega(t)^* \l \bigl(i_{\varphi}\Omega(t)\bigr)
\]
is smooth on $X\times\{|t|<\epsilon\}$ for some small $\epsilon$.
This completes the proof of Step \ref{step 3}, and hence of Theorem \ref{main thm}.

			\section{Proof of King's quasi-isomorphisms}\label{KINGKING}
		We provide here a detailed proof of Proposition \ref{PPP} for the comfort of the readers.
		\begin{proposition}[{=Proposition \ref{PPP}}]\label{Q}Let $X$ be a complex manifold of complex dimension $n$, and let $Y$ be a simple normal crossing divisor on $X$. Then we have the following properties.
			\begin{enumerate}[{(a)}]
				\item\label{Ap1}The sheaf $\mathscr{D}^{\prime \bullet,\bullet}_X (\log Y)$ is a bigraded  $\mathscr{A}^{0,\bullet} (\Omega_{X}^\bullet(\log Y))$-module.
				\item \label{Ap2}
				This  module structure from  \ref{Ap1} defines a canonical isomorphism: $$\Omega_{X}^p(\log Y)\otimes_{\mathcal{O}_X} \mathscr{D}_X^{\prime 0,q}\stackrel{\simeq}{\longrightarrow} \mathscr{D}_X^{\prime p,q}(\log Y).$$ 
				\item\label{Ap3}There exist quasi-isomorphisms of complexes of sheaves:
				\begin{equation*}
					\Omega_{X}^p(\log Y)\hookrightarrow\mathscr{A}^{0,\bullet}(\Omega_{X}^p(\log Y))\hookrightarrow \mathscr{D}^{\prime p,\bullet}_X (\log Y).
				\end{equation*}
				So,  $( \mathscr{D}^{\prime p,\bullet}_X (\log Y),\bp)$ is an acyclic resolution of $ \Omega_{X}^p(\log Y).$
			\end{enumerate}
		\end{proposition}

		The above proposition is purely a local problem. Hence, we may assume  $Y=\left\{z^1\cdots z^k=0\right\}$ and   set $h=z^1\cdots z^k$.  We also denote  the coordinates as $(z^1,\ldots,z^k, w^{k+1},\ldots,w^n)$.
		\begin{proof}[{Proof of Proposition \ref{Q} \ref{Ap1}}]
			Suppose that $T$ is a current and  $g$ is  a holomorphic function such that $h$ can be divided by $g$, then it would be desirable to define $T\wedge \frac{dg}{g}$ to be $\frac{T}{g}\wedge dg$ where $\frac{T}{g}$ is a suitably chosen current $S$ satisfying $gS=T$. In fact, by virtue of the work of Schwartz \cite{Sch55}, it is always possible to find such an $S$. However, $S$ is not unique, since it can be modified by adding a  current whose product with $g$ equals  zero. Conversely, if $gS=gS^\prime$, then $g(S-S^\prime)=0.$ Sometimes, it is possible to define explicit procedures to obtain a preferred  $S$  under various transversality conditions, but in general, this is impossible, and the  current $\frac{T}{g}\wedge dg$ is not well-defined. However, by passing to the complex of log currents through the quotient map, the ambiguity disappears.
			We will define the following maps to form a commutative diagram:
			\[
			\begin{tikzcd}[row sep=5.2em, column sep=5.2em] 
				\mathscr{D}_X^{\prime\bullet,\bullet} \otimes_{\mathscr{A}_X} \mathscr{A}^{0,\bullet}(\Omega_{X}^\bullet(\log Y)) 
				\arrow[r, "\alpha"] 
				\arrow[d, "\gamma \otimes\mathrm{id}"'] 
				&	\mathscr{D}_X^{\prime\bullet,\bullet}(\log Y).  \\
				\mathscr{D}_X^{\prime\bullet,\bullet}(\log Y)\otimes_{\mathscr{A}_X} \mathscr{A}^{0,\bullet}(\Omega_{X}^\bullet(\log Y))	
				\arrow[ur, "\overline{\alpha}"']
			\end{tikzcd}
			\]
			
			Here the map $\gamma$ is the quotient map $\mathscr{D}^{\prime \bullet,\bullet}_X\rightarrow \mathscr{D}^{\prime \bullet,\bullet}_X (\log Y)$.  Let $T$ be a current and $\mu$ be a logarithmic form. Choose a current $S$ such that $hS=T$. 
			The form $h\mu$ is then smooth, therefore the product $S\wedge(h\mu)$ is defined. This current depends on the choice of $S$ but $\gamma(S\wedge(h\mu)$) just depends on $T$ and $\mu$. Indeed, if $S^\prime$ is another current satisfying $hS^\prime=T$, then $hQ=0$, where $Q:=S-S^\prime$. For any null-$Y$ form $\varphi$ with suitable bidegree,
			as currents, $$(Q\wedge(h\mu))\wedge\varphi=Q\wedge(h\mu\wedge \varphi)=hQ\wedge(\mu\wedge\varphi)=0,$$ since $\mu\wedge\varphi$ is smooth. Therefore, $Q\wedge(h\mu)$ is an on-$Y$ current by the equivalent Definition  \eqref{kkk}. In other words, $$\gamma(S\wedge(h\mu))=\gamma(S^\prime\wedge(h\mu))$$  is  a well-defined log current, which we denote by
			$$\alpha(T\otimes\mu).$$
			
			To define the map $\overline{\alpha},$ it  suffices to show that $$	\mathscr{D}_X^{\prime\bullet,\bullet}(\mathrm{on}\, Y) \otimes_{\mathscr{A}_X} \mathscr{A}^{0,\bullet}(\Omega_{X}^\bullet(\log Y))\subset\mathrm{Ker}\,\alpha.$$ If $T$ above is an on-$Y$ current and the forms are the same as before, we then have $$(S\wedge(h\mu))\wedge\varphi=hS\wedge(\mu\wedge\varphi)=T\wedge(\mu\wedge\varphi)=0$$
			since $\mu\wedge\varphi$ is a null-$Y$ form according to the ``ideal'' property  in Proposition \ref{null sheaf} \ref{null sheaf 2}. Therefore, $S\wedge(h\mu)$ is also an on-$Y$ current,  which completes the proof.
		\end{proof}
		\begin{proof}[{Proof of Proposition \ref{Q} \ref{Ap2}}]
			If $\alpha$ and $\gamma$ are the same maps in the proof of \ref{Ap1}, one then has a commutative diagram:
			\[
			\begin{tikzcd}[row sep=4.4em, column sep=4.4em]
				& {\mathscr{D}_X^{\prime p,q}} \arrow[ld, "\beta"'] \arrow[rd, "\gamma"] &                                                                                       \\
				{\mathscr{D}_X^{\prime 0,q}\otimes_{\mathcal{O}_X}\Omega_{X}^p(\log Y)}  \arrow[rr, "\alpha"] &                                                                         & {\mathscr{D}_X^{\prime p,q}(\log Y)}. \arrow[ll, "\overline{\beta}", dashed, bend left]
			\end{tikzcd}
			\]
			 Here if $T=\sum_{|I|+|J|=p}T_{IJ} dz^I\wedge dw^J,$ where $T_{IJ}$ is a $(0,q)$-current,
			then we set $$\beta(T)=\sum_{|I|+|J|=p}(z^IT_{IJ})\otimes\frac{dz^I}{z^I}\wedge dw^J.$$
			We should also observe that in this case,  $\alpha=\overline{\alpha},$ since $\mathscr{D}^{\prime 0,q}_X (\log Y)=\mathscr{D}^{\prime 0,q}_X$ (see Remark \ref{rmk 1} \ref{p=0}). 
			
			Through the rest of the proof, we will adopt the convention that for a current $T$ and a holomorphic function $g$, $\frac{T}{g}$ refers to  any \textit{fixed} current $S$ such that $gS=T.$
			
			\vspace{0.4em}
			\noindent(1) \textit{$\beta$ is surjective.}\quad 
			
			\vspace{0.4em}
			This follows because
			$$\sum_{|I|+|J|=p} R_{IJ}\otimes\frac{dz^I}{z^I}\wedge dw^J=\beta\Bigg(\sum_{|I|+|J|=p} \Big(\frac{R_{IJ}}{z^I}\Big)\,dz^I\wedge dw^J\Bigg).$$

			\noindent(2) \textit{$\overline{\beta}$ is surjective.}\quad 
			
			\vspace{0.4em}
			We first show that $\overline\beta$ is well-defined.  Then,  it suffices to show that ${\mathscr{D}_X^{\prime p,q}}(\mathrm{on}\, Y)\subset\mathrm{Ker}\,\beta.$
			For any null-$Y$ form $\varphi$ of  $(n-p,n-q)$-type,  by \eqref{868686} we can write
			\begin{equation*}
				\begin{aligned}
					\varphi&=\sum_{|K|+|L|=n-p}h\cdot\varphi_{KL}\frac{dz^K}{z^K}\wedge dw^L\\
					&=\sum_{|K|+|L|=n-p}\varphi_{KL} \,z^{K^\prime}dz^K\wedge dw^L,
				\end{aligned}
			\end{equation*}
			where $\varphi_{KL}$ is a $(0,n-q)$-smooth  form, and 
			$K^\prime=\left\{1,\ldots,k\right\}-K.$
			Therefore, $T$ is an on-$Y$ current if and only if 
			$$z^I T_{IJ}=0$$
			as currents. Then in this case, $\beta(T)=0$.  
			
			As a result, $\overline{\beta}$ is also surjective thanks to (1) above.

			\vspace{0.6em}
			\noindent(3) \textit{$\alpha\circ\beta=\gamma.$}
			\begin{equation*}
				\begin{aligned}
					\alpha(\beta(T))&=
					\gamma\Bigg(\sum\Big(\frac{z^IT_{IJ}}{h}\Big)\wedge h\frac{dz^I}{z^I}\wedge dw^J\Bigg)\\
					&=	\gamma\Bigg(\sum\Big(\frac{T_{IJ}}{z^{I^\prime}}\Big)\wedge z^{I^\prime}dz^I\wedge dw^J\Bigg)\\
					&=
					\gamma\Big(\sum{T_{IJ}}dz^I\wedge dw^J\Big)=\gamma(T),
				\end{aligned}
			\end{equation*}
			where $I^\prime=\left\{1,\ldots,k\right\}-I.$
			Consequently, $\alpha\circ\overline{\beta}=\mathrm{id}$ so $\overline{\beta}$ is injective and $\alpha={\overline{\beta}}^{\,-1}$.

			\vspace{1em}
			Putting (1), (2) and (3) together, we get 
			$$\Omega_{X}^p(\log Y)\otimes_{\mathcal{O}_X} \mathscr{D}_X^{\prime 0,q}\stackrel{\simeq}{\longrightarrow} \mathscr{D}_X^{\prime p,q}(\log Y).$$ 
		\end{proof}

		\begin{proof}[{Proof of Proposition \ref{Q} \ref{Ap3}}]
			Thanks to \ref{Ap2}, the sequence
			$$	\Omega_{X}^p(\log Y)\hookrightarrow\mathscr{A}^{0,\bullet}(\Omega_{X}^p(\log Y))\hookrightarrow \mathscr{D}^{\prime p,\bullet}_X (\log Y)$$
			can be obtained by tensoring the sequence
			$$\mathcal{O}_X\hookrightarrow\mathscr{A}_X^{0,\bullet}\hookrightarrow\mathscr{D}_X^{\prime 0,\bullet}$$ with the locally free $\mathcal{O}_X$-sheaf $\Omega_{X}^p(\log Y).$ The desired quasi-isomorphisms in  the former sequence follows since the latter sequence consists of quasi-isomorphisms, as established by the celebrated Dolbeault--Grothendieck lemma.
		\end{proof}

		\section{Directed limits}
		We give some basic knowledge of the directed limit of a directed family of abelian groups to be used in the proof of Theorem \ref{inj-}.
		
		\begin{definition}
			A \textit{directed family} of abelian groups $\left\{A_m\right\}_{m \in \mathbb{Z}}$ is a collection of abelian groups equipped with homomorphisms $\phi_{mn}: A_m \to A_n$ for all $m \leq n$, satisfying the following conditions:
			\begin{itemize}
				\item[--]$\phi_{mm} = \mathrm{id}_{A_m}$ for all $m$ (identity property);
				\item[--] $\phi_{np} \circ \phi_{mn} = \phi_{mp}$ for all $m \leq n \leq p$ (compatibility property).
			\end{itemize}
		\end{definition}
		
		\begin{definition}
			The \textit{directed limit} $\varinjlim\limits_m A_m$ of $\left\{A_m\right\}_{m \in \mathbb{Z}}$  is defined as the quotient of $\bigoplus\limits_{m\in\mb Z}A_m$ modulo the subgroup generated by $$x_m-\phi_{mn}(x_m)$$
			for all $m\leq n$ and $x_m\in A_m$.
		\end{definition}
		One  obtains from this definition the homomorphisms 
		$$\mu_m:A_m\rightarrow \varinjlim\limits_n A_n, \quad a_m\mapsto [a_m],$$
		which are compatible with $\phi_{mn}$, and satisfy the following \textit{universal property}:
		if $G$ is an abelian group and $g_m: A_m\rightarrow G$ are homomorphisms compatible with $\phi_{mn}$, then there exists a unique homomorphism 
		$$g:\varinjlim_n A_n\rightarrow G$$
		such that $g_m=g\circ\mu_m$ 
		for all $m$.
		
		From the explicit description of the directed limit, we have the following properties:
		\begin{enumerate}[{(i)}]
			\item $\varinjlim\limits_n A_n=\bigcup\limits_m \mu_m(A_m)$;
			\item $\mathrm{Ker}\,(A_m\rightarrow\varinjlim\limits_n A_n)=\bigcup\limits_{m\leq n}\mathrm{Ker}\,(A_m\rightarrow A_n).$
		\end{enumerate}
		In particular, one obtains the following lemma:
		\begin{lemma}\label{direc}
			Let $\left\{A_m\right\}_{m \in \mathbb{Z}}$ be a directed family of abelian groups.
			\begin{enumerate}
				\item We have that $A_m\rightarrow\varinjlim\limits_n A_n$ is injective if and only if $A_m\rightarrow A_n$ is injective for all $m\leq n$.
				\item 	Let $\left\{B_m\right\}_{m \in \mathbb{Z}}$ be another  directed family of abelian groups. Let $f_m:A_m\rightarrow B_m$ be a sequence of compatible homomorphisms, which induce a homomorphism $$f:\varinjlim\limits_m A_m\rightarrow\varinjlim\limits_m B_m.$$ If $f_m$ is injective for some $m\geq m_0$, then $f$ is injective.
			\end{enumerate}
		\end{lemma}

		\vspace{2em}
	\noindent\textbf{\large Declarations}
	
	\vspace{1em}
	\noindent\textbf{Conflict of interest}\quad The author stated that there is no conflict of interest.

\end{document}